%% file: main.tex
\documentclass[a4paper,english]{article}

\usepackage[english]{babel}
\usepackage{bibgerm}
\usepackage[utf8]{inputenc}
\usepackage[hyphens]{url}
\usepackage{hyperref}
\usepackage{caption}
\usepackage{enumerate}
\usepackage{xcolor}
\usepackage{xspace}
\usepackage{amsmath,amssymb,amsfonts, amsthm}
\usepackage[a4paper, total={6.5in, 9in}]{geometry}
\usepackage{rotating}
\usepackage{fancyvrb}
\usepackage{subfigure}
\usepackage{caption}
\usepackage{subcaption}
\usepackage{graphicx}
\usepackage{mathtools}
\DeclarePairedDelimiter\ceil{\lceil}{\rceil}
\DeclarePairedDelimiter\floor{\lfloor}{\rfloor}
\usepackage{pgfplots}
\usepackage{pifont}
\usepackage{pdflscape}
\usepackage{multirow}

\usepackage{tikz}
\usetikzlibrary{patterns}
\usetikzlibrary{positioning}

\usetikzlibrary{shapes}
\usetikzlibrary{shapes.geometric}
\usetikzlibrary{backgrounds}
\usetikzlibrary{arrows.meta}

\usepackage{zibtitlepage}
\usepackage{authblk}

\usepackage{natbib}
 \bibpunct[, ]{(}{)}{,}{a}{}{,}

\newtheorem{theorem}{Theorem}[section]
\newtheorem{lemma}[theorem]{Lemma}
\newtheorem{definition}[theorem]{Definition}
\newtheorem{corollary}[theorem]{Corollary}
\newtheorem{remark}[theorem]{Remark}

\newtheorem{example}[theorem]{Example}

\newcommand\conv{\operatorname{conv }}
\newcommand\pesp{\text{PESP}}
\newcommand\xpesp{\text{XPESP}}
\newcommand\cxpesp{\text{cXPESP}}
\newcommand\ttp{\text{TTP}}
\newcommand\xttp{\text{XTTP}}
\newcommand\cxttp{\text{cXTTP}}

\newcommand\cxpespw{\text{cXPESP}^{w}}

\newcommand{\BFx}{\mathbf{x}}
\newcommand{\BFa}{\mathbf{a}}

\begin{document}
\ZTPAuthor{
    \ZTPHasOrcid{Stephanie Riedmüller}{0009-0006-4508-4262},
    \ZTPHasOrcid{Niels Lindner}{0000-0002-8337-4387},
}
\ZTPNumber{xx-xx}
\ZTPMonth{October}
\ZTPYear{2024}

\title{Column Generation for Periodic Timetabling}

\author[1]{Stephanie Riedmüller\footnote{riedmueller@zib.de, \ZTPOrcid{0009-0006-4508-4262}}}
\author[1,2]{Niels Lindner\footnote{lindner@zib.de, \ZTPOrcid{0000-0002-8337-4387}}}

\affil[1]{Zuse Institute Berlin, Berlin, Germany}
\affil[2]{Freie Universtität Berlin, Berlin, Germany}

\maketitle

\begin{abstract}
Periodic timetabling for public transportation networks is typically modelled as a Periodic Event Scheduling Problem (PESP). Solving instances of the benchmark library PESPlib to optimality continues to pose a challenge. As a further approach towards this goal, we remodel the problem by a time discretization of the underlying graph and consider arc-based as well as path-based integer programming formulations. For the path-based case, we use cycles on the graph expansion of the operational lines as variables and, therefore, include more of the problem inherent structure into the model. A consequence is the validity of several known inequalities and a lower bound on the LP-relaxation, that is the best known to date. As an extension we integrate passenger routing into the new model. 
The proposed models have an advantage in the linear programming relaxation, on the one hand, but have an increased problem size, on the other hand. We define the corresponding pricing problems for the use of column generation to handle the size. Both models are practically tested on different problem instances. 
\end{abstract}

\section{Introduction}

The timetable is the backbone of a public transportation network.
A careful timetable design is hence key to offer attractive service, to enable efficient operation, and to contribute towards more sustainable mobility.
The optimization of timetables however remains to be a challenging planning step. The goal is often a periodic timetable, i.e., a repeating schedule for a specific time period. 
There are several techniques available to optimize periodic timetables, all of which struggle when applied to large data sets.
It is therefore an ongoing topic of interest to find new computational perspectives.

The most common model for periodic timetabling is the Periodic Event Scheduling Problem (PESP), first described for a periodic scheduling problem in terms of job shops by \cite{SerafiniUkovich1989}.
A variety of integer programming formulations are available for PESP \citep{Liebchen2006PeriodicTO}, leading to success stories, e.g., the first mathematically optimized timetable in practice \citep{liebchen_first_2008}.
On the downside, computing optimal timetables for large instances continues to pose challenges:
The instances of the benchmark library PESPlib (\url{https://timpasslib.aalto.fi/pesplib.html}) remain to be unsolved for more than a decade despite several attempts.

It is therefore natural to ask for reformulations.
Our contribution is influenced by two approaches:
In aperiodic railway timetabling, it is common to work on time-expanded networks \citep{brannlund_railway_1998,caprara_modeling_2002,schlechte_railway_2012}. 
Another theme, which is also common practice independent of timetabling, is to transform arc-based models into path-based formulations in connection with column generation techniques \citep{Barnhart1998,borndorfer_column-generation_2007,cacchiani_column_2008}.

We transfer these ideas to the Periodic Event Scheduling Problem. On the basis of the \xpesp{} model introduced by \cite{Kinder2008}, we apply a time discretization to the underlying event-activity network of PESP and contribute new integer programming models based both on arcs and on paths (\cxpesp). 
Due to the periodicity of the problem, the occurring paths are mostly cycles.  
The increase of variables, resulting from the large number of possible cycles, is compensated by inclusion of the inherent problem structure, which is not present in traditional methods, and allows for stronger dual bounds of the linear programming relaxation.
Moreover, we cope with the problem size by analyzing column generation methods, and show that the resulting pricing problems can be solved via shortest path computations on directed acyclic graphs.

Furthermore, timetabling and passenger routing affect each other in the process of optimizing public transport. 
Most models are based on the assumption that the chosen routes of passengers are independent of the timetable (\cite{Liebchen2006PeriodicTO}, \cite{Lindner2000}).
We will follow instead an integrated approach, extending \cite{BorndoerferHoppmannKarbstein2017}, and show that passenger routing integrates seamlessly into our time-expanded timetabling model (\cxttp) and the column generation process. 

Finally, we evaluate our models and techniques computationally.
We start with solving the linear programming relaxation of cXPESP. First, we investigate the integrality gap between its optimal solution and the optimal solution of the PESP integer program.
Next, we compare different variants of cXPESP in terms of the number of pricing rounds, generated variables, and computation time, including the portion spent on pricing. 
Furthermore, we solve the integer program of the best-performing cXPESP variant and compare its performance to state-of-the-art PESP formulations. In addition to the mentioned pricing statistics, we also evaluate the number of nodes explored in the branch-and-bound tree.
For the case of integrated passenger routing, we analyze the linear programming relaxation of cXTTP regarding its integrality gap, and the pricing behavior for cycles and passenger paths.

The paper is structured as follows. We review related literature in Section~\ref{section:literature}.
In Section~\ref{section:PESP}, we recall PESP and its typical underlying event-activity networks. We describe in Section~\ref{section:expansion} an expansion of the event-activity network based on a time discretization and a reformulation of PESP on this graph based on arcs. 
In Section~\ref{section:cxpesp}, we introduce the novel path-based formulation. After a description of the models and their linear programming formulation, we investigate the polyhedral structure of the solution space and compare it to PESP. Due to the increasing model size, we look into column generation and its model specific pricing problem.
In Section~\ref{section:integration}, we integrate passenger routing into the model and state its corresponding pricing problem. 
Section~\ref{section:computation} is dedicated to a computational investigation, before we conclude the paper in Section~\ref{section:conclusions}.

\section{Literature Review}
\label{section:literature}

The basis of our considerations is the \emph{Periodic Event Scheduling Problem} (PESP), as introduced by \cite{SerafiniUkovich1989}. The literature on PESP is rich, and we refer to \cite{nachtigall_periodic_1998} and \cite{Liebchen2006PeriodicTO} for a general overview.

As finding feasible solutions is already NP-hard for a fixed period time \citep{Odijk1994} or on series-parallel graphs \citep{lindner_analysis_2022}, starting heuristics based on SAT techniques (\citealt{grossmann2012}) and mimicking the Phase I of the simplex method have been developed (\citealt{goerigk_phase_2021}).
A plenty of improving heuristics are available as well: The polyhedral structure of PESP is exploited by the modulo network simplex heuristic (\citealt{nachtigall_solving_2008}, \citealt{goerigk_improving_2013}) and tropical neighborhood search \citep{BortolettoLindnerMasing2022}. Further approaches include, e.g., maximum cuts \citep{lindner_new_2019}, multi-agent reinforcement learning \citep{matos_solving_2021}, solution merging \citep{LindnerLiebchen2022}, and line-based decomposition approaches \citep{patzold_matching_2016,lindner_incremental_2023}. 
The currently best known primal bounds on the PESPlib instances have been achieved by a combination of many of these approaches \citep{BorndoerferLindnerRoth2020}.

Providing good dual bounds is notoriously hard. Several types of cutting planes are known, e.g., the cycle inequalities \cite{Odijk1994}, and the change-cycle inequalities by \cite{Nachtigall1996}. The best-performing cutting plane technique relies on the separation of the more general class of flip inequalities, that are equivalent to split cuts \citep{lindner_split_2025}.

A time-expanded version of PESP has been developed by \cite{Kinder2008}, leading to an arc-based integer programming formulation. One advantage is that this expansion provides a natural linearization when integrating passenger routing, as observed by \cite{BorndoerferHoppmannKarbstein2017}, while using the standard PESP formulation leads to a more compact, but quadratic problem \citep{lubbe_passagierrouting_2009}. The integrated periodic timetabling and passenger routing problem is an active research topic \citep{schmidt_timetabling_2015,borndorfer_passenger_2017,schiewe2020,lobel_restricted_2020,loebel_geometric_2025}. Recently, the benchmarking library TimPassLib\footnote{\url{https://timpasslib.aalto.fi}} has been established \citep{schiewe_introducing_2023}.

\cite{martin-iradi_column-generation-based_2022} describe a time-expanded path-based integer programming formulation for periodic timetabling, but with a narrower scope than PESP, and targeted at conflict-free symmetric railway timetables. Their model includes passenger routing, column generation, and Benders decomposition.

Concerning aperiodic timetabling, time-space networks are the foundation for the influential integer programming models by \cite{brannlund_railway_1998} and \cite{caprara_modeling_2002}. The transformation of arc-based models to path-based formulations, often in conjunction with column generation, has been investigated not only in railway timetabling \citep{cacchiani_column_2008,schlechte_railway_2012,min_appraisal_2011}, but also in other fields of public transport optimization, e.g., line planning \citep{borndorfer_column-generation_2007}, vehicle scheduling \citep{ribeiro_column_1994}, and crew scheduling \citep{Barnhart1998}.

\section{The Periodic Event Scheduling Problem}\label{section:PESP}

The Periodic Event Scheduling Problem (PESP) is based on an \textit{event-activity network}. In this section, we recall the definition both of the graph structure and of PESP, which will function as a basis for further introductions and as a reference model. To simplify notation, we will assume that all graphs under consideration are simple.

\subsection{The Event-Activity Network}

The underlying graph in PESP is an \textit{event-activity network} $N$, which is a directed graph, whose vertices are called \textit{events} and whose edges are called \textit{activities}. A \emph{line network} is an undirected graph $G$, together with a set $\mathcal{L}$ of \emph{lines}, where each line is a path in $G$. In public transportation, determining the line network is typically preceding the timetabling phase (see, e.g., \citealt{bussiek1997}).
We will consider event-activity networks $N$ constructed as follows, as in \cite{masing_forward_2023}:
\begin{itemize}
    \item For each line $L \in \mathcal{L}$ and each edge $\{v,w\} \in E(L)$ add \textit{departure events} $(v,L,dep, +)$ and $(w,L,dep, -)$ and \textit{arrival events} $(w,L,arr, +)$ and $(v,L,arr, -)$ to $V(N)$. Furthermore, add \textit{driving activities} $((v,L,dep,+),(w,L,arr,+))$ and $((w,L,dep,-),(v,L,arr,-))$ to $E(N)$.
    \item For each line $L \in \mathcal{L}$ and each stop $v \in V(L)$ add \textit{waiting activities} $((v,L,arr,+),(v,L,dep,+))$ and $((v,L,arr,-),(v,L,dep,-))$ to $E(N)$ if the corresponding events exist.
    \item For each line $L \in \mathcal{L}$ and for the first and last stop $v_s, v_t \in V(L)$ of $L$ add \textit{turnaround activities} $((v_s,L,arr,-),(v_s,L,dep,+))$ and $((v_t,L,arr,+),(v_t,L,dep,-))$ to $E(N)$.
    \item For each $v \in V(G)$ and each pair $(L,L')$ of distinct lines with $v \in V(L) \cap V(L')$ and for $s_1, s_2 \in \{ +, - \}$ add a \textit{transfer activity} $((v,L,arr, s_1),(v,L', dep, s_2))$ to $E(N)$ if both events exist.
\end{itemize} 
For an undirected line $L= (v_1, \ldots, v_n) \in \mathcal{L}$ in $G$, we define the directed path
\begin{equation*}
    ((v_1,L,dep, +), (v_2,L,arr, +), (v_2,L,dep, +), \ldots, (v_{n-1},L,dep, +), (v_n,L,arr, +))
\end{equation*}
as the \textit{forward direction} of $L$ in $N$ and 
\begin{equation*}
    ((v_n,L,dep, -), (v_{n-1},L,arr, -), (v_{n-1},L,dep, -), \ldots, (v_{n-1},L,dep, -), (v_1,L,arr, -))
\end{equation*}
as the \textit{backward direction} of $L$ in $N$.
For a given line $L$ in the event-activity network $N$, we call the closed path consisting of the forward and backward direction of $L$ together with its turnaround activities the \textit{line cycle} of $L$.
We consider the \textit{frequency} of each line, that is, how often a line is served by a vehicle in a given time frame, to be fixed to one for the purpose of simplification.

\subsection{Problem Definition}

Given an event-activity network $N$, lower and upper bounds $l,u \in \mathbb{Z}^{A(N)}, l \leq u$, weights $\omega \in \mathbb{Q}^{A(N)}$
and a \textit{period time} $T \in \mathbb{N}$ with $[T] = \{ 0, \ldots, T-1\}$, the \textit{Periodic Event Scheduling Problem (PESP)} is to find 
$\pi \in [T]^{V(N)}$ and  $\BFx \in \mathbb{Z}^{A(N)}$ such that
\begin{itemize}
\item $\forall (v,w) \in A(N): \BFx_{vw} \equiv \pi_w - \pi_v \bmod T,$
\item $l \leq \BFx \leq u,$
\item $\omega^T\BFx$ is minimal,
\end{itemize}
or decide that such $\pi$ does not exist.
The resulting vector $\pi \in [T]^{V(N)}$ is called a \textit{periodic timetable}. The periodic timetable assigns to each event $v \in V(N)$ a \textit{potential} $\pi_v \in [T]$, that can be interpreted as the time at which a vehicle arrives at or departs from a station in the given period.
The \textit{periodic tension} $\BFx \in \mathbb{Z}^{A(N)}$ represents the duration of the activities.
We assume the periodic tension and periodic timetable to be integral. 

There are several integer programming formulations for PESP
available. The following formulation models the modulo operator with a vector $p \in \mathbb{Z}^{A(N)}$:
\begin{align} 
	\min & \sum_{\alpha = (v,w) \in A(N)} \omega_\alpha (\pi_w - \pi_v + T p_\alpha)  &\text{PESP} \notag\\
	&\pi_w - \pi_v + T p_\alpha \geq l_\alpha   &\forall \alpha \in A(N)\label{eq:pesp1}\\
	&\pi_w - \pi_v + T p_\alpha \leq u_\alpha  &\forall \alpha \in A(N)  \label{eq:pesp2}\\
	&0 \leq \pi_v   \leq T - 1& \forall  v \in V(N) \label{eq:pesp3}\\ 
	& \pi_v \in \mathbb{Z}  & \forall v \in V(N) \label{eq:pesp4} \\
	& p_\alpha \in \mathbb{Z}  & \forall \alpha \in A(N). \label{eq:pesp5}
\end{align}
The periodic tension is included only implicitly: For $\alpha = (v, w) \in A(N)$, we have $x_\alpha = \pi_w - \pi_v + T p_\alpha$.
For the purpose of comparison to other models, we define the following solution polytopes:
\begin{definition} \label{def:polypesp}
Denote by
\begin{align*}
	& P_{IP}(\text{PESP}) \coloneqq \conv \{ (\pi, p) \in \mathbb{Z}^{V(N)} \times \mathbb{Z}^{A(N)} \mid (\pi, p) \text{ satisfies } (\ref{eq:pesp1}) - (\ref{eq:pesp5})\},\\ 
	& P_{LP}(\text{PESP}) \coloneqq \{ (\pi, p) \in \mathbb{Q}^{V(N)} \times \mathbb{Q}^{A(N)} \mid (\pi, p) \text{ satisfies } (\ref{eq:pesp1}) - (\ref{eq:pesp3})\}
\end{align*}
the polyhedra associated to the integer program and linear programming relaxation for PESP, respectively.
\end{definition}

\section{The Expanded Periodic Event Scheduling Problem} \label{section:expansion}

Requiring timetables to attain integer values, the event-activity network $N$ can be expanded by a time discretization, which will be the underlying graph of further models in this work.

\subsection{Expansion of the Event-activity Network}

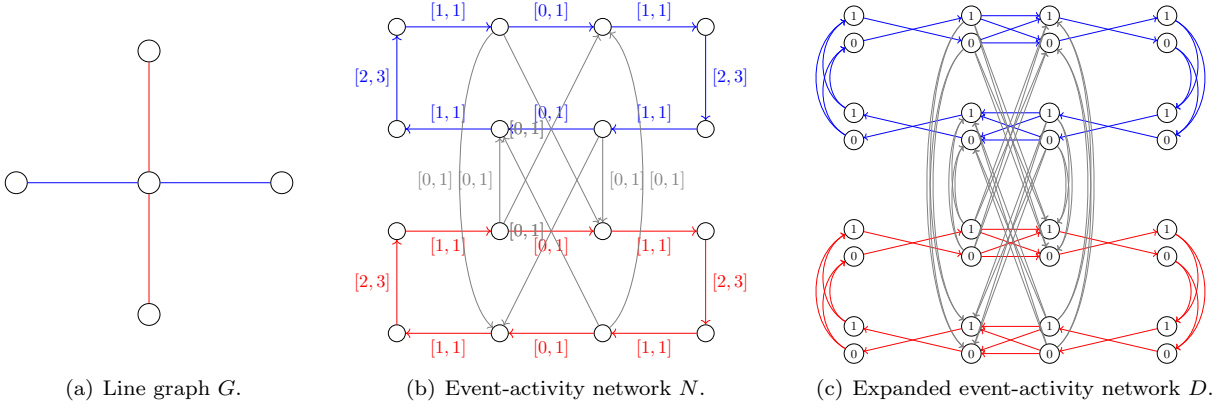
\begin{figure}
\hfill
\centering
\subfigure[Line graph $G$.]{\raisebox{0.5cm}
{\resizebox{4cm}{!}{
  \input{graphics/figure1}}
}}
\hfill
\subfigure[Event-activity network $N$.]
{\resizebox{5.5cm}{!}{
  \input{graphics/figure2}}
}
\hfill
\subfigure[Expanded event-activity network $D$.]
{\resizebox{5.5cm}{!}{
  \input{graphics/figure3}}
}
\caption{Exemplary graph expansion. (a) The line network $G$ of a two-line example with driving activity bounds fixed to 1. (b) The corresponding event-activity network $N$ with exemplary bounds for waiting, transfer, and turnaround activities. (c) The corresponding expanded event-activity network $D$ for period $T=2$.
For each event, the different time steps are marked by the numbers in the nodes.}
\label{fig:graph-comparison}
\end{figure}

Let $T \in \mathbb{N}$ be a fixed period time. 
The \textit{expanded event-activity network} $D$ for a time period $T$ is constructed by \citep{Kinder2008}:
\begin{itemize}
\item For each event $v \in V(N)$ and for each time step $t \in [T]$
add a node $v[t]$ to $V(D).$
\item For each activity $\alpha = (v,w) \in A(N)$ add to $A(D)$ all arcs of the set
$$\mathcal{A}(\alpha) = \{(v[t], w[t']) \in V(D) \times V(D) \mid t, t' \in [T], (t'-t-l_{\alpha}) \bmod T \leq u_{\alpha} - l_{\alpha} \}.$$
\end{itemize}
Each arc $a \in A(D)$ has a fixed \textit{duration} given by $\tau_a \coloneqq (t' - t -l_\alpha) \bmod T + l_\alpha$.
Consequently, all arcs in the expanded event-activity network $D$ represent feasible activities, since they obey the lower and upper bounds of the input data. 
Hence, a feasible solution to a given PESP instance corresponds to a subgraph of the expanded event-activity network: The tension of a given activity $\alpha$ in PESP is represented by the choice of exactly one arc in $\mathcal{A}(\alpha)$. The value of the timetable at a given event $v \in V(N)$ is determined by the choice of exactly one $v[t] \in V(D)$. An illustration of an exemplary graph expansion can be found in Figure~\ref{fig:graph-comparison}.

For a line cycle $(v_1, \ldots, v_n)$ in $N$, we call $D$ restricted to the node set
$$\{v[t] \mid v \in  \{v_1, \ldots, v_n\}, t \in [T]\}$$
an \textit{expanded line cycle}. See, for example, the blue colored expanded line cycle in Figure~\ref{fig:graph-comparison}c.

We denote by $X(N) \subseteq A(N)$ the set of \textit{vehicle activities} (driving, waiting, and turnaround activities), and by $Z(N) \subseteq A(N)$ the set of transfer activities, so that $A(N) = X(N) \cup Z(N)$.
Analogously, we denote by $X(D) \subseteq A(D)$  the set of arcs belonging to expanded line cycles
and by $Z(D) \subseteq A(D)$ the set of transfer arcs, so that $A(D) = X(D) \cup Z(D)$.

The construction of the expanded event-activity network results in $|V(D)| = |V(N)| \cdot T$ nodes.
For each activity $\alpha \in A(N)$ there are $(u_{\alpha} - l_{\alpha} +1  \bmod T) \cdot T$ arcs in $A(D)$, and, therefore, 
\begin{equation}
\label{eq:expanded_arcs_count}
  |A(D)| = T \cdot \sum_{\alpha \in A(N)} (u_{\alpha} - l_{\alpha} +1  \bmod T) \leq T \cdot \sum_{\alpha \in A(N)} (T-1) \in \mathcal{O}(T^2 \,|A(N)|) .  
\end{equation}

\subsection{XPESP: An Arc-based Model}
Several integer programming formulations for arc-based models on an expanded event-activity network are provided by \cite{Kinder2008}, one of them is the following:
For each transfer arc $a \in Z(D)$, we introduce a binary variable $z_a \in \{0,1\}$ and for each arc $a \in X(D)$ a binary variable $x_a \in \{0,1\}$ that models if the arc is present in the chosen subgraph. 
For the objective we consider for each arc $a \in A(D)$ its weight $\omega_a \in \mathbb{Q}$ multiplied by its duration $\tau_a \in \mathbb{Q}$. 
Denote by $\delta^{+/-}_S(v)$ the outgoing/ingoing arcs of a node $v \in V(D)$, respectively, restricted to the set $S \subseteq A(D)$.
For each activity $\alpha \in A(N)$ and each arc $a \in \mathcal{A}(\alpha)$, the lower and upper bound is given by $[l_a, u_a] = [l_{\alpha}, u_{\alpha}]$, and we assume $\omega_a = \omega_\alpha$ for $a \in \mathcal{A}(\alpha)$.
\begin{align} 
	\min & \sum_{a \in X(D)} \omega_a \tau_a x_a + \sum_{a \in Z(D)} \omega_a \tau_a z_a  &\text{XPESP} \notag\\
	& \sum_{a \in \mathcal{A}(\alpha)} x_a= 1 & \forall \alpha \in X(N) \label{eqn:cut}\\ 
	& \sum_{a \in \delta_{X(D)}^-(v) } x_{a} - \sum_{a \in \delta_{X(D)}^+(v) } x_{a}  = 0 & \forall v \in V(D)\label{eqn:flowconsxpesp}\\ 
	& \sum_{a \in \delta_{X(D)}^-(v[t]) } x_{a} - \sum_{a \in \delta_{\mathcal{A}(\alpha)}^+(v[t]) } z_{a}  = 0 & \forall \alpha = (v,w) \in Z(N), t \in [T]\label{eqn:tr1}\\
	& \sum_{a \in \delta_{X(D)}^+(w[t']) } x_{a} - \sum_{a \in \delta_{\mathcal{A}(\alpha)}^-(w[t']) } z_{a}  = 0 & \forall \alpha = (v,w) \in Z(N), t' \in [T]\label{eqn:tr2}\\
	&0 \leq x_a   \leq 1& \forall  a\in X(D) \label{eqn:bound1}\\
	&0 \leq z_a   \leq 1& \forall  a\in Z(D) \label{eqn:bound2}\\
	& x_a \in \mathbb{Z}  & \forall a \in X(D) \label{eqn:int1}\\
	& z_a \in \mathbb{Z}  & \forall a \in Z(D) \label{eqn:int2}
\end{align}

XPESP is formulated similar to a min cost flow problem.
Constraint~\eqref{eqn:cut} partitions the flow over $\mathcal{A}(\alpha)$ to exactly one arc $a \in \mathcal{A}(\alpha)$, 
Constraint~\eqref{eqn:flowconsxpesp} ensures flow conservation on nodes of an expanded line cycle and Constraints~\eqref{eqn:tr1} and \eqref{eqn:tr2} are coupling constraints for the transfers between lines.

For the remainder of this section, we include the following definition and lemma, which we will encounter in later chapters:
\begin{definition} \label{def:polyxpesp}
Denote by
\begin{align*}
	& P_{IP}(\text{XPESP}) = \conv \{ (x,z) \in \mathbb{Z}^X \times \mathbb{Z}^{Z} \mid (x,z) \text{ satisfies } (\ref{eqn:cut}) - (\ref{eqn:int2})\},\\
	& P_{LP}(\text{XPESP}) = \{ (x,z) \in \mathbb{Q}^X \times \mathbb{Q}^{Z} \mid (x,z) \text{ satisfies } (\ref{eqn:cut}) - (\ref{eqn:bound2})\}
\end{align*}
the polyhedra associated to the integer and linear program for XPESP, respectively. 
\end{definition}

\begin{lemma} \label{lem:xpesplowerbounds}
The optimal solution value to the linear programming relaxation of XPESP
is the weighted sum of lower bounds on the activities
$$\sum_{\alpha \in A(N)} \omega_\alpha l_\alpha.$$
\end{lemma}

\cite{Kinder2008} proves the lower bound for XPESP to be zero, if the objective function minimizes the slack, i.e., the difference of the tensions to the lower bounds on the activities. Lemma~\ref{lem:xpesplowerbounds} follows directly from that proof.

\begin{remark}
\label{rem:non-linear}
We use $\omega_a = \omega_\alpha$ for $a \in \mathcal A(\alpha)$ only for convenience. The XPESP model and our later developments do work with arbitrary arc weights $\omega_a$ in the expanded network $D$, and can hence model non-linear objective functions in terms of the periodic tensions of the activities in $N$ as well.
\end{remark}

\section{cXPESP: A Path-based Model} \label{section:cxpesp}

\begin{figure}
\centering
  \input{graphics/figure4}
\caption{A cycle (blue) in an expanded line cycle. The closed walk colored in red does not fulfill the definition of a cycle due to the cardinality of its node set.}
\label{fig:cycle}
\end{figure}
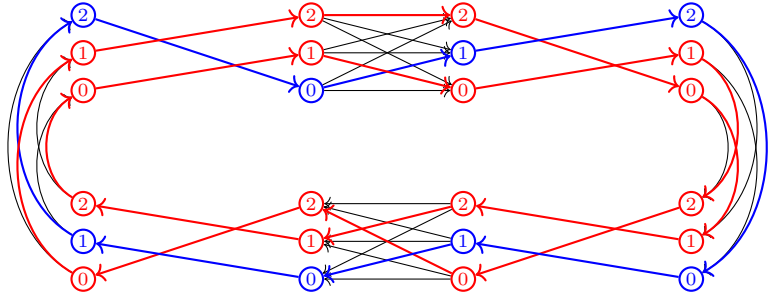

In contrast to Kinder's arc-based model, we introduce the new path-based model cXPESP, which provably includes more of the problem-inherent structure than PESP and XPESP by exploiting the expanded line cycles. This arises to be especially beneficial when dealing with its linear programming relaxation. To that end, we introduce the notion of a \textit{cycle} in an expanded line as illustrated in Figure~\ref{fig:cycle}.
\begin{definition}
Let $\gamma_L$ be the line cycle of a line $L \in \mathcal{L}$ in the event-activity network $N$ and let $\textbf{c}_L$ be the corresponding expanded line cycle in the expanded event-activity network $D$. We call a closed path $c$ in $\textbf{c}_L$ a \textit{cycle} if $|V(c)| = |V(\gamma_L)|$.
We denote by $C_L $ the set of cycles in the expanded line cycle $\textbf{c}_L$ for line $L \in \mathcal{L}$ and by $C$ the set of all cycles $C = \bigcup_{L \in \mathcal{L}} C_L$.
\end{definition}  

The goal is again to determine an optimal subgraph of the expanded event-activity network. For each line $L$ and for each possible cycle $c \in C_L$,
we introduce a variable
$$x_c = \begin{cases} 
		1 &\text{if line cycle }c\text{ is selected,}\\
		0 &\text{otherwise.} \\
	\end{cases}$$
The model inherits the variables $z_a \in \{0,1\}$ for transfer arcs $a \in Z(D)$ from XPESP.
Further, $\tau_a$ and $\tau_c$ denote the duration for each arc $a = (v[t], w[t']) \in A(D)$ and for each cycle $c \in C$, respectively. The durations are computed by $$\tau_a \coloneqq (t' - t - l_a ) \bmod T + l_a \qquad \text{ and }  \qquad \tau_c \coloneqq \sum_{a \in c} \tau_a.$$
Note that $\tau_c$ is necessarily an integer multiple of $T$. As before, we denote the weight of an arc $a \in A(D)$ as $\omega_a$. We further write
$$\vartheta_c \coloneqq \sum_{a \in A(c)} \omega_a \tau_a$$
for the weighted duration of a cycle $c \in C$ with arc set $A(c)$.

\begin{figure}
\centering
  \input{graphics/figure5}
\caption{Illustration of the coupling constraints~\eqref{eq:c1new} and \eqref{eq:c2new} in cXPESP. The summed values of the identically colored arc variables equal the summed values of the correspondingly colored cycle variables.}
\label{fig:couplingconst}
\end{figure}
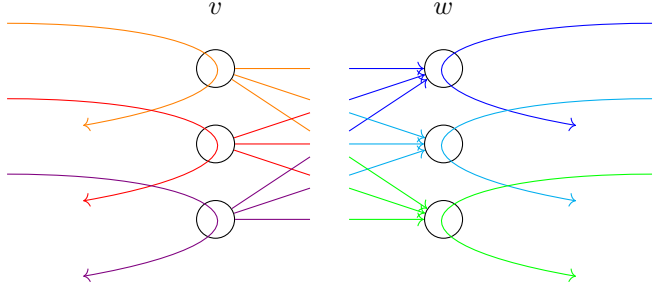

We define cXPESP to be formulated as
\begin{align}
	\min & \sum_{c \in C} \vartheta_c x_c + \sum_{a \in Z(D)} \omega_a \tau_a z_a  & \text{cXPESP} \notag\\
	& \sum_{c \in C_L} x_c = 1 & \forall L \in \mathcal{L} \label{eq:partition}\\
	& \sum_{c \in C: v[t] \in c} x_{c} - \sum_{a \in \delta^+_{\mathcal{A}(\alpha)}(v[t])} z_a  =  0 & \forall \alpha = (v,w) \in Z(N), t \in [T] \label{eq:c1new}\\
	&  \sum_{c \in C: w[t'] \in c} x_{c} - \sum_{a \in \delta^-_{\mathcal{A}(\alpha)}(w[t'])} z_a  =  0 & \forall \alpha = (v,w) \in Z(N), t' \in [T] \label{eq:c2new}\\
	& x_c   \geq 0& \forall c \in C	\label{eq:b1}\\
	& z_a   \geq 0& \forall a \in Z(D) \label{eq:b2}\\
	& x_c \in \mathbb{Z}  & \forall c \in C	\label{eq:int1}\\
	& z_a \in \mathbb{Z}  & \forall a \in Z(D). \label{eq:int2}
\end{align}
The partitioning constraint~\eqref{eq:partition} ensures the resulting subgraph to include exactly one cycle per line.
The coupling constraints~\eqref{eq:c1new} and \eqref{eq:c2new} describe flow conservation at each node of an expanded line cycle, where at least one outgoing arc is a transfer arc, as illustrated in Figure~\ref{fig:couplingconst}. 
Note that a coupling constraint is necessary for each node and pair of lines with a possible transfer.
All variables are implicitly binary, since each $x_c$ is at most 1 due to Constraint~\eqref{eq:partition} and each $z_a$ is at most 1 due to
$$z_a \leq \sum_{a \in \delta^+_{\mathcal{A}(\alpha)}(v[t])} z_a \overset{\eqref{eq:c1new}}{=} \sum_{c \in C: v[t] \in c} x_{c} \leq  \sum_{c \in C_L} x_c \overset{\eqref{eq:partition}}{=} 1,$$
where $a = (v[t],w[t])$ and $L$ is the line through $v$.

\begin{definition} 
Denote by
\begin{align*}
	& P_{IP}(\text{cXPESP}) = \conv \{ (x,z) \in \mathbb{Z}^C \times \mathbb{Z}^{Z} | (x,z) \text{ satisfies } \eqref{eq:partition} - \eqref{eq:int2}\},\\ 
	& P_{LP}(\text{cXPESP}) = \{ (x,z) \in \mathbb{Q}^C \times \mathbb{Q}^{Z} | (x,z) \text{ satisfies } \eqref{eq:partition}-\eqref{eq:b2}\}
\end{align*}
the polyhedra associated to the integer program and linear programming relaxation for cXPESP, respectively.
\end{definition}

\subsection{Comparison of Problem Size}

The maximal number of variables and constraints across the different models is described in Table~\ref{tb:numbervariablesconstraints} dependent on the size of the event-activity network $N$, the period $T$, and the line set $\mathcal{L}$.
We do not count upper and lower bound constraints.
The number of nodes in the longest line cycle is denoted by $k \coloneqq \max_{L \in \mathcal L} |V(C_L)|$.

In PESP, there are exactly one variable and two constraints for each activity in $A(N)$, plus the potential variables for each event in $V(N)$.
XPESP involves one variable for each arc in $A(D)$ as well as one constraint for each activity $\alpha \in X(N)$, $2T$ constraints for each activity $\alpha \in Z(N)$,
and one constraint for each node $v[t] \in V(D)$.
In total, we obtain $|A(D)| \in \mathcal O(T^2 |A(N)|)$ variables (see \eqref{eq:expanded_arcs_count}), and
$$|X(N)| + 2T \cdot | Z(N) | + T \cdot |V(N)| \in \mathcal{O}\left(T \cdot (|A(N)| + |V(N)|)\right)$$ constraints.

In cXPESP, there are at most $T^k$ possible cycles in each expanded line cycle for each line in $\mathcal{L}$.
We have further one variable for each transfer arc in $|Z(D)|$ and $|Z(D)|\leq T^2 \cdot |Z(N)|$.
There is one constraint for each line in $\mathcal{L}$ and $2T$ constraints for each activity in $Z(N)$. 

\begin{table}
\centering
\caption{Comparison of problem sizes.}
\label{tb:numbervariablesconstraints}
\begin{tabular}{ l l l}
  \hline
  & no. variables & no. constraints \\ 
 \hline 
 PESP & $|A(N)| + |V(N)|$ &  $\mathcal{O}( |A(N)|)$ \\  
 XPESP& $\mathcal{O}(T^2 \,|A(N)|)$ &  $ \mathcal{O}\left(T \, (|A(N)| + |V(N)|)\right)$  \\
 cXPESP & $\mathcal{O}( T^k  \,  |\mathcal{L}| + T^2 \, |Z (N)|) $   &  $\mathcal{O}(|\mathcal{L}| + T \,|Z (N)|)$ \\
  \hline
\end{tabular}
\end{table}

Table~\ref{tb:numbervariablesconstraints} shows that XPESP and cXPESP have significantly larger size than PESP.
While XPESP has the largest number of constraints, cXPESP shows the largest number of variables due to the exponential number of possible cycles in an expanded line cycle.

\subsection{Comparison of Solution Polytopes}

The increased number of variables in cXPESP provides a richer structure and we therefore investigate the relationship between the linear programming relaxations of PESP, XPESP and cXPESP. Since the solution spaces of these three relaxations differ in dimension, statements must be given under linear transformation.

\begin{definition}\label{def:transformation}
Define $\psi$ to be the linear transformation 
\begin{align*}
	\psi : \mathbb{Q}^{X(D)} \times \mathbb{Q}^{Z(D)} & \rightarrow \mathbb{Q}^{V(N)} \times \mathbb{Q}^{A(N)}\\
 	\begin{pmatrix}{c} x \\ z \end{pmatrix} & \mapsto \begin{pmatrix}{c} \Pi \quad 0\\  P \end{pmatrix} \cdot \begin{pmatrix}{c} x \\ z \end{pmatrix},
\end{align*}
where
\begin{align*} 
	\Pi &\in \mathbb{Q}^{V(N) \times X(D)} && \Pi = SR, \\
	R &\in \mathbb{Q}^{V(D) \times X(D)} && R_{v,a} = \begin{cases} 
		t &\text{if $a = (v[t], w[t']) \in \delta^+_X(v)$,}\\
		0 &\text{otherwise,} \\
	\end{cases}\\
	S &\in \{0,1\}^{V(N) \times V(D)} && S_{v,w} = \begin{cases} 
		1 &\text{if $w = v[t]$ for some $t$,}\\
		0 &\text{otherwise,} \\
	\end{cases}\\
	P &\in \mathbb{Q}^{A(N) \times A(D)} && P_{\alpha,a} = \begin{cases} 
		0 &\text{if $a = (v[t], w[t']) \in \mathcal{A}(\alpha)$ and $t \leq t'$,}\\
		1 &\text{otherwise.} \\
	\end{cases}
\end{align*} 
and define $\varphi$ to be the linear transformation
		\begin{align*}
			\varphi : \mathbb{Q}^C \times \mathbb{Q}^{Z(D)} & \rightarrow \mathbb{Q}^{X(D)} \times \mathbb{Q}^{Z(D)}\\
			 \begin{pmatrix}{c} x \\ z \end{pmatrix} & \mapsto \begin{pmatrix}{cc} M_C & 0 \\  0 &I \end{pmatrix} \cdot \begin{pmatrix}{c} x \\ z \end{pmatrix}, 
		\end{align*}
		where 
		\begin{align*}
			&M_C = (m_{ac})_{a \in X(D), c \in C},  &m_{ac} = \begin{cases} 1 &\text{if $a \in c$,}\\ 0 &\text{otherwise.} \end{cases}
		\end{align*}
\end{definition}

\begin{theorem} \label{thm_transformation_LP}
The linear transformations $\psi$ and $\varphi$ have the property:
\begin{equation*}
    \psi(\varphi(P_{LP}(\cxpesp ))) \subseteq \psi(P_{LP}( \xpesp)) \subseteq P_{LP}( \pesp).
\end{equation*}
\end{theorem}

\proof{Proof}
$\varphi(P_{LP}(\cxpesp )) \subseteq P_{LP}( \xpesp)$:\\
Let $(x^{\circ} , z^{\circ})  \in P_{LP}(\cxpesp )$. 
We show $(\bar{x}, \bar{z}) \coloneqq \varphi (x^{\circ} , z^{\circ}) \in P_{LP}(\xpesp )$.	
Note that for any $a \in X(D)$ holds
\begin{align}
     \bar{x}_a = \sum_{c \in C:a \in c} x^{\circ}_c.	\label{eq:def_phi}
\end{align}
\begin{itemize}
    \item (\ref{eqn:cut}) Let $\alpha \in X(N)$ be an activity of line $L \in \mathcal{L}$. Then  
    $$\sum_{a \in \mathcal{A}(\alpha)} \bar{x}_a \overset{\eqref{eq:def_phi}}{=} \sum_{a \in \mathcal{A}(\alpha)} \sum_{c \in C:a \in c} x^{\circ}_c = \sum_{c \in C_L} x^{\circ}_c \overset{\eqref{eq:partition}}{=} 1.$$
    \item (\ref{eqn:flowconsxpesp})  Let $v \in V(D)$. Then $$\sum_{a \in \delta^-_{X(D)}(v)} \bar{x}_a \overset{\eqref{eq:def_phi}}{=}  \sum_{a \in \delta^-_{X(D)}(v)}\sum_{c \in C: a \in c} x^{\circ}_c =  \sum_{c \in C : v\in c} x^{\circ}_c   = \sum_{a \in \delta^+_{X(D)}}\sum_{c \in C : a \in c} x^{\circ}_c \overset{\eqref{eq:def_phi}}{=} \sum_{a \in \delta^+_{X(D)}(v)} \bar{x}_a.$$
    \item (\ref{eqn:tr1}) and analogously (\ref{eqn:tr2}) Let $\alpha = (v,w) \in Z(N)$ and $t \in [T]$. Then
    $$\sum_{a \in \delta^+_{\mathcal{A}(\alpha)}(v[t])} \bar{z}_a = \sum_{a \in \delta^+_{\mathcal{A}(\alpha)}(v[t])} z^{\circ}_a \overset{\eqref{eq:c1new}}{=} \sum_{c \in C: v[t] \in c} x^{\circ}_{c}=  \sum_{a \in \delta^-_{X}(v[t])}\sum_{c \in C: a \in c} x^{\circ}_c \overset{\eqref{eq:def_phi}}{=} \sum_{a \in \delta^-_{X(D)}(v[t])} \bar{x}_a. $$
    \item (\ref{eqn:bound1}) Let $a\in X(D)$. Then, for a unique $L \in \mathcal{L}$,  \begin{align*} 
    0 \leq \bar{x}_a =  \sum_{c \in C:a\in c} x^{\circ}_c  \leq \sum_{c \in C_L} x^{\circ}_c \overset{\eqref{eq:partition}}{=} 1. \end{align*}
    \item (\ref{eqn:bound2}) Let $\BFa = (v[t],w[t']) \in Z(N)$. Then, for an unique $L \in \mathcal{L}$,  
    \begin{align*}
    0 \leq \bar{z}_{\BFa} = z^{\circ}_{\BFa} \leq   \sum_{a \in \delta^+_{\mathcal{A}((v,w))}(v[t])} z^{\circ}_a \overset{\eqref{eq:c1new}}{=}  \sum_{c \in C: v[t]\in c} x^{\circ}_c \leq \sum_{c \in C_L} x^{\circ}_c \overset{\eqref{eq:partition}}{=} 1. \end{align*}
\end{itemize}
The inclusion $\psi(P_{LP}( \xpesp)) \subseteq P_{LP}( \pesp)$ has been proven by \cite{Kinder2008} assuming $l_\alpha < T$ for all $\alpha \in A(N)$.
\endproof

\begin{theorem} \label{thm_transformation_IP}
The linear transformations $\psi$ and $\varphi$ have the property:
\begin{equation*}
    \psi(\varphi(P_{IP}(\cxpesp ))) = \psi(P_{IP}( \xpesp)) = P_{IP}( \pesp).
\end{equation*}
\end{theorem}

\proof{Proof}
$\varphi(P_{IP}(\cxpesp )) = P_{IP}( \xpesp)$:\\
	Notice that an integral solution $(x^{\circ} , z^{\circ}) \in \mathbb{Z}^C \times \mathbb{Z}^{Z(D)}$ is mapped to an integral solution  $(\bar{x}, \bar{z}) = \varphi (x^{\circ} , z^{\circ}) \in P_{LP}(\xpesp ) \in \mathbb{Z}^{X(D)} \times \mathbb{Z}^{Z(D)}$, since the defining matrices have exclusively integer entries. As linear maps preserve convex combinations, $\varphi(P_{IP}(\cxpesp )) \subseteq P_{IP}( \xpesp)$ holds. For the reverse inclusion, we show that the restricted map $$\varphi_{IP}: P_{IP}(\cxpesp ) \rightarrow P_{IP}(\xpesp )$$ is surjective on integer points.
	Let $(\bar{x}, \bar{z}) \in P_{IP}(\xpesp ) \cap (\mathbb Z^{X(D)} \times \mathbb Z^{Z(D)})$ be an integer point. We construct $(x^{\circ} , z^{\circ}) \in P_{IP}(\cxpesp)$ such that $(\bar{x}, \bar{z}) = \varphi_{IP} (x^{\circ} , z^{\circ})$ by setting
	\begin{align*}
		&x^{\circ}_c \coloneqq \begin{cases} 1 & \text{if } \bar{x}_a = 1 \text{ for all } a \in c, \\ 0 & \text{otherwise,} \end{cases} &&\text{and} &&&z^{\circ}_a \coloneqq \bar{z}_a
	\end{align*}
	for all $c \in C$ and $a \in Z$, respectively.  Intuitively, $\bar{x}$ is supported on exactly one cycle per expanded line cycle, and we select precisely those cycles for $x^\circ$. It remains to check the formal details.
    
    First, we show that $(x^{\circ}, z^{\circ})$ is indeed an element of $P_{IP}(\cxpesp )$: 
    \begin{itemize}
        \item Constraint~\eqref{eq:partition} can be    validated by a proof of contradiction using case distinction to contradict Constraint~\eqref{eqn:cut} and Constraint~\eqref{eqn:flowconsxpesp}.
        \item By case distinction, $\sum_{c \in C:v[t] \in c} x^{\circ}_c = \sum_{a \in \delta^+_X(v[t])} \bar{x}_a$ holds for $t \in [T]$ and, therefore, Constraints~\eqref{eq:c1new} and \eqref{eq:c2new}.
	\end{itemize}
	Second, we show that $(x^{\circ}, z^{\circ})$ is indeed mapped to $(\bar{x}, \bar{z})$ by case distinction. Let $a \in Z(D)$.
    \begin{itemize}
        \item Let $\bar{x}_a = 0$. Since $x^{\circ}_c = 0$ for all $c$ containing $a$, we 
        obtain $$\varphi_a(x^{\circ} , z^{\circ})_a = \sum_{c \in C: a \in c} x^{\circ}_c = 0 = \bar{x}_a.$$
        \item Let $\bar{x}_a = 1$. If $x^{\circ}_c = 0$ for all $c$ containing $a$, we obtain a contradiction to the flow conservation constraints in \xpesp. If $x^{\circ}_c = 1$ for some $c$ containing $a$, then $c$ is uniquely determined. Otherwise, Constraint~\eqref{eq:partition} would be violated. Therefore, $$\varphi_a (x^{\circ} , z^{\circ})_a = \sum_{c \in C : a \in c} x^{\circ}_c = 1 = \bar{x}_a.$$ 
    \end{itemize}

The identity $\psi(P_{IP}( \xpesp)) = P_{IP}( \pesp)$ has been proven by \cite{Kinder2008} assuming $l_\alpha < T$ for all $\alpha \in A(N)$.	
\endproof

\begin{corollary}
    There is a one-to-one correspondence between the integer solutions of PESP, XPESP and cXPESP. In particular, XPESP and cXPESP are correct.
\end{corollary}

\proof{Proof} 
While the one-to-one correspondence between the integer solutions of PESP and XPESP has been proven by \cite{Kinder2008}, the one-to-one correspondence between the integer solutions of XPESP and cXPESP follows from Theorem~\ref{thm_transformation_IP} by injectivity of the restricted map $\varphi_{IP}$ on integer points.
Assume for the sake of contradiction $(x_1^{\circ} , z_1^{\circ}) \ne (x_2^{\circ} , z_2^{\circ}) \in P_{IP}(\cxpesp)$ such that 
$(\bar{x}, \bar{z}) = \varphi_{IP} (x_1^{\circ} , z_1^{\circ})= \varphi_{IP} (x_2^{\circ} , z_2^{\circ})$. For all activities $\alpha \in X(N)$, there is exactly one arc $a \in \mathcal{A}(\alpha)$ with $\bar{x}_a = 1$. Thus, for all line cycles $C$, there is exactly one cycle $c \in C$ including those non-zero arcs and hence $x^{\circ}_{1c}=x^{\circ}_{2c} = 1$. As we must have that $z^\circ_1 = z^\circ_2$ by definition of $\varphi_{IP}$, we conclude that $(x_1^{\circ} , z_1^{\circ}) = (x_2^{\circ} , z_2^{\circ})$.

\endproof

\begin{figure}
\centering
\subfigure[PESP: Line graph G.]{\raisebox{0.8cm}{{\resizebox{5cm}{!}{
  \input{graphics/figure7}}
}}}
\hfill
\subfigure[XPESP: Solution to LP relaxation.]{{\resizebox{5.3cm}{!}{
  \input{graphics/figure8}}
}}
\hfill
\subfigure[cXPESP: Possible cycles.]{{\resizebox{5.3cm}{!}{
  \input{graphics/figure9}}
}}
\caption{Illustration of counterexample in Remark~\ref{rem:incl}.}
\label{fig:ce}
\end{figure}
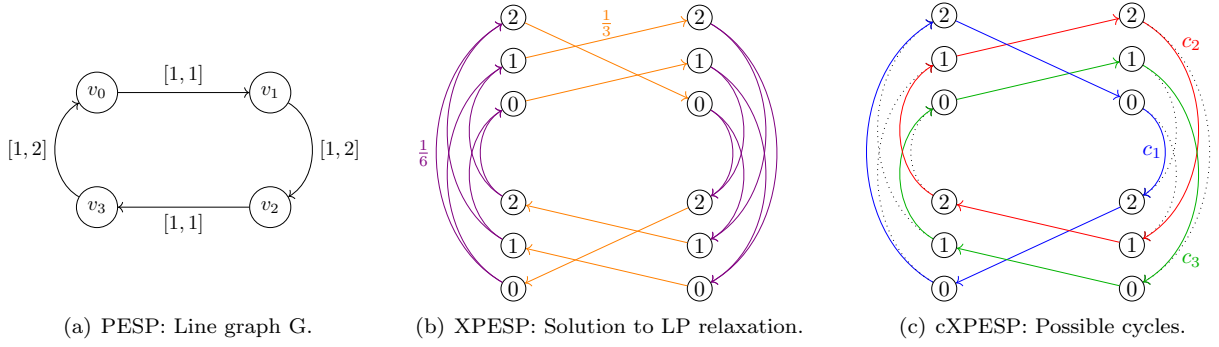

\begin{remark} \label{rem:incl}
In general, $\varphi(P_{LP}(\cxpesp )) = P_{LP}( \xpesp)$ and $\psi(P_{LP}( \xpesp)) = P_{LP}( \pesp)$ do not hold.
By counterexample: For the first equation, consider the expanded event-activity network in Figure~\ref{fig:ce} consisting of a single expanded line cycle with a fixed driving time of $1$, a turnaround time in $[1,2]$ and period $T=3$.	
Consider $(\bar{x}, \bar{z}) \in P_{LP}( \xpesp)$ such that 
$$\bar{x} = (\bar{x}_a)_{a \in A(D)}, \qquad \bar{x}_a \coloneqq \begin{cases}  \frac{1}{3} &\text{if $a \in \mathcal{A}(\alpha), \alpha \in X(N)$  driving activity,}\\  \frac{1}{6} &\text{if $a \in \mathcal{A}(\alpha), \alpha \in X(N)$ turnaround activity.} \end{cases}$$

Let $c_1, c_2, c_3$ denote the possible cycles in the expanded line cycle.
We show that there is no $x^{\circ} \in P_{LP}(\cxpesp )$ such that $\varphi(x^{\circ}) = \bar{x}$. Assume there is such an element $x^{\circ}= (x^{\circ}_{c_1},x^{\circ}_{c_2},x^{\circ}_{c_3})^\top \in P_{LP}(\cxpesp )$.
Then the following equation must hold, considering the rows corresponding to the arcs $(v_0[0], v_1[1])$ and $(v_1[1], v_2[0])$:

\begin{equation*}
    \varphi \begin{pmatrix}{c} x^{\circ}_{c_1}\\x^{\circ}_{c_2} \\ x^{\circ}_{c_3} \end{pmatrix}  = \begin{pmatrix}{ccc}  
    &\vdots&  \\ 0 & 0 & 1 \\  0 & 0 & 1 \\   &\vdots& \end{pmatrix}  \cdot \begin{pmatrix}{c} x^{\circ}_{c_1}\\x^{\circ}_{c_2} \\ x^{\circ}_{c_3} \end{pmatrix}=  \begin{pmatrix}{c}  
    \vdots   \\x^{\circ}_{c_3} \\ x^{\circ}_{c_3} \\  \vdots  \end{pmatrix} = \begin{pmatrix}{c}   
    \vdots  \\1/3 \\ 1/6 \\  \vdots \end{pmatrix} = 
    \begin{pmatrix}{ccc} 
    \vdots \\ \bar{x}_{(v_0[0], v_1[1])} \\  \bar{x}_{(v_1[1], v_2[0])}\\   \vdots \end{pmatrix}
    =\bar{x}.
\end{equation*}\
		This is a contradiction, for example by $1/6 = x^{\circ}_{c_3} = 1/3$. Hence, $\varphi(P_{LP}(\cxpesp )) \ne P_{LP}( \xpesp)$.
A counterexample for the second identity is given by \cite{Kinder2008}.
\end{remark}

This result shows that the linear programming relaxation of cXPESP can indeed be stronger than the linear programming relaxation of XPESP and PESP, while having the same integer solution space.
We claim that the relaxation of cXPESP, furthermore, gives an advantage for the optimal value. Recall from Lemma~\ref{lem:xpesplowerbounds}
that the optimal value to XPESP is the weighted sum of lower bounds on the activities. This is not necessarily the case for cXPESP and we can construct a counterexample. 

\begin{example} \label{ex:notzero}
Again, consider the example given in Figure~\ref{fig:ce}.
There are three different possible cycles in the expanded line cycle and all of them have duration $6 > \sum_{a \in A(N)} l_a = 4$.  
\end{example}

\begin{remark}
	The concatenation  $\varphi \circ \psi$ of the transformations given in Definition~\ref{def:transformation} can be used to project a solution to cXPESP to a corresponding solution to PESP. 
	If one is solely interested in the tensions of the projection onto PESP, a single periodic tension for some $\alpha \in X(N)$ is given by
	$$\textbf{\textup{x}}_{\alpha} = \sum_{a \in \mathcal{A}(\alpha)} \tau_{a} \cdot \left( \sum_{c \in C_L: a \in c} x^{\circ}_c \right)$$
	and for some $\alpha \in Z(N)$ by 
	$$\textbf{\textup{x}}_{\alpha} = \sum_{a \in \mathcal{A}(\alpha)} \tau_{a}  z^{\circ}_a.$$
\end{remark}

\subsection{Valid Inequalities}

Since cXPESP includes more of the problem structure than other PESP variants, we show that there exist additional valid cuts.
The following theorem holds for solutions to cXPESP with a focus on single line cycles. 
\begin{theorem} \label{thm:splitclosure}
Let $(x^{\circ},z^{\circ}) \in P_{LP}(\cxpesp)$ and $(\pi, p) = \psi \circ \varphi((x^{\circ},z^{\circ})) \in P_{LP}(\pesp)$. Let $T \in \mathbb{N}$ be a fixed period, let $l,u \in \mathbb{Z}^A$ be the lower and upper bounds and denote by $\Gamma$ the set of line cycles in the event-activity network $N$. If $\textbf{\textup{x}}$ is the tension corresponding to $(\pi, p)$, then 
	$$\textbf{\textup{x}} \in \bigcap_{\gamma \in \Gamma} \conv \left\{ \textbf{\textup{y}} \in \mathbb{Z}^{A(N)}\mid l \leq \textbf{\textup{y}} \leq u,  \frac{\gamma^T \textbf{\textup{y}}}{T} \in \mathbb{Z} \right\}.$$
\end{theorem}

\proof{Proof}

Let $\alpha=(v, w) \in A(N)$ be an arbitrary activity in the event-activity network $N$. 
Let $(x^{\circ},z^{\circ}) \in P_{LP}(\cxpesp)$,  $(\bar{x}, \bar{z}) \in P_{LP}(\xpesp)$ and $(\pi, p) \in P_{LP}(\pesp)$ such that 
$$(\pi, p) = \psi \left((\bar{x}, \bar{z}) \right)= \psi \circ \varphi \left((x^{\circ},z^{\circ}) \right).$$
Using the definition of the transformation $\psi$, the periodic tension $\BFx$ for an activity $\alpha$ in the given solution $(\pi, p)$ is computed by 
\begin{align}
			\BFx_{\alpha} &= \pi_w - \pi_v + Tp_{\alpha}  \notag \\
				&= \sum_{a = (v[t],w[t']) \in \mathcal{A}(\alpha)} t' \bar{x}_a - \sum_{a = (v[t],w[t'])\in \mathcal{A}(\alpha)} t \bar{x}_a + T  \cdot \sum_{a = (v[t],w[t'])\in \mathcal{A}(\alpha), t > t'} \bar{x}_a  \notag \\ 
				& = \sum_{a \in \mathcal{A}(\alpha)} \beta_a \bar{x}_a 	\label{eq:tension}
\end{align}
with $$\beta_a = \beta_{(v[t],w[t'])} \coloneqq \begin{cases}
								t' - t + T & \text{if } t>t', \\
								t' - t 		& \text{otherwise.}
							\end{cases}$$
Consider a line cycle $\gamma \in \Gamma$ and denote by $C$ the expansion of $\gamma$ in the expanded event-activity network $D$. Then 
\begin{align*}
			(\BFx_{\alpha})_{\alpha \in \gamma} &\overset{\eqref{eq:tension}}{=} \left( \sum_{a \in \mathcal{A}(\alpha)} \beta_a \bar{x}_a \right)_{\alpha \in \gamma} \overset{\eqref{eq:c1new}}{=}  \left( \sum_{a \in \mathcal{A}(\alpha)} \beta_a \sum_{c \in C: a\in c} x^{\circ}_c \right)_{\alpha \in \gamma}\\
			&  = \left( \sum_{c \in C} \left(\sum_{a\in \mathcal{A}(\alpha) \cap c}  \beta_a \right)  x^{\circ}_c \right)_{\alpha \in \gamma}= \left( \sum_{c \in C} \textbf{\textup{y}}_{\alpha}  x^{\circ}_c \right)_{\alpha\in \gamma} =  \sum_{c \in C}  x^{\circ}_c \left(\textbf{\textup{y}}_{\alpha}\right)_{\alpha \in \gamma} , 
\end{align*}
where $\textbf{\textup{y}}_\alpha \coloneqq \sum_{a \in A(\alpha) \cap c} \beta_a \in \mathbb{Z}_{\geq 0}$. Since $\sum_{c \in C}  x^{\circ}_c = 1$ due to Constraint~\eqref{eqn:cut}, it remains to check that the vector $\textbf{\textup{y}} $ is a feasible periodic tension. Indeed, $ \mathcal{A}(\alpha) \cap c$ contains a unique arc, and the collection of these arcs for $\alpha \in \gamma$ is precisely the arc set of $c$. Hence $\textbf{\textup{y}}_\alpha \in [l_\alpha, u_\alpha]$ and
$$ \frac{\gamma^\top \textbf{\textup{y}}}{T} = \frac{1}{T} \sum_{\alpha \in \gamma} \sum_{a \in \mathcal{A}(\alpha) \cap c} \beta_a = \frac{1}{T}  \sum_{a \in c}  \beta_a  \in \mathbb Z,$$
resolving the telescoping sum in the definition of $\beta_a$.

\endproof

We use this result to see that some known inequalities are valid for the projection of a cXPESP solution to a PESP solution. 
To that end, we recall the following established results, where each $\gamma$ is considered to be a vector in $\{0, -1, 1\}^{X(N)}$. The entries in the vector represent if an activity is present in the oriented cycle and determine its direction. Decompose the cycle into positive and negative directions $\gamma = \gamma^+ - \gamma^-$ and note that  $\gamma = \gamma^+$ if $\gamma$ is a line cycle.
\begin{lemma}[\citealt{Odijk1994}] \label{lem:cycleinequality} 
	Let $\gamma$ be an oriented cycle, $(\pi, p) \in P_{IP}(\pesp)$ and \textbf{\textup{x}} the corresponding periodic tension. Then, the \textit{cycle inequality} 
			$$\ceil*{\frac{\gamma^T_+l - \gamma_-^T u}{T}} \leq \frac{\gamma^T \textbf{\textup{x}}}{T}  \leq \floor*{\frac{\gamma^T_+ u - \gamma_-^T l}{T}} $$ 
   is valid.
\end{lemma}
\begin{lemma}[\citealt{Nachtigall1996}] \label{lem:changecycleinequality} 
	Let $\gamma$ be an oriented cycle, $(\pi, p) \in P_{IP}(\pesp)$ and $\textbf{\textup{x}}$ the corresponding periodic tension. Then the \textit{change-cycle-inequality} 
			$$(T- \xi) \gamma_+^T (\textbf{\textup{x}}-l) + \xi \gamma_-^T (\textbf{\textup{x}}-l) \geq \xi (T-\xi), \quad \xi = [- \gamma^Tl]_T$$ is valid. 
\end{lemma}
\begin{lemma}[\citealt{lindner_determining_2020}]\label{lem:flipcycleinequality} 
	Let $F \subseteq A$, $\gamma$ be an oriented cycle, $(\pi, p) \in P_{IP}(\pesp)$ and $\textbf{\textup{x}}$ the corresponding periodic tension. Then the \textit{flip-cycle inequality}
		\begin{align*}
			 &(T- \xi_F) \sum_{a \in A \setminus F: \gamma_a=1} (\textbf{\textup{x}}_a - l_a) + \xi_F \sum_{a \in A \setminus F: \gamma_a=-1} (\textbf{\textup{x}}_a - l_a) \\
			& + \xi_F \sum_{a \in F: \gamma_a = 1} (u_a - \textbf{\textup{x}}_a) + (T - \xi_F) \sum_{a \in F: \gamma_a = -1} (u_a - \textbf{\textup{x}}_a) \geq \xi_F(T- \xi_F),
		\end{align*}
		where $$ \xi_F = \left[- \sum_{a \in A \setminus F} \gamma_a l_a - \sum_{a \in F} \gamma_a u_a \right]_T$$ is valid.
\end{lemma}

The stated Lemmata and Theorem~\ref{thm:splitclosure} yield the following theorem.

\begin{theorem} \label{thm:inequalities}
 	Let $(x^{\circ},z^{\circ}) \in P_{LP}(\cxpesp)$ and $(\pi, p) = \psi \circ \varphi((x^{\circ},z^{\circ})) \in P_{LP}(\pesp)$. Denote by $\Gamma$ the set  of line cycles in the event-activity network $N$. If $\textbf{\textup{x}}$ is the periodic tension corresponding to $(\pi, p)$, then
	\begin{enumerate}
		\item the cycle inequality,  
		\item the change-cycle-inequality, and
		\item the flip-cycle inequality
	\end{enumerate}
	hold for each $\gamma \in \Gamma$.
\end{theorem}

\proof{Proof}
    Let $\textbf{\textup{x}}$ be the periodic tension corresponding to $(\pi, p)$.    Theorem~\ref{thm:splitclosure} has an interpretation in terms of split inequalities \citep{cook_chvatal_1990}: It implies that $\textbf{\textup{x}}$ satisfies all split inequalities for $P_{LP}(\pesp)$ with respect to all split disjunctions given by cycles in $\Gamma$. However, split inequalities are the same as flip-cycle inequalities \cite[Theorem 3.1]{lindner_split_2025}, of which the cycle and change-cycle inequalities are special cases \citep{lindner_determining_2020}. 
\endproof

\begin{remark}
Notice that there are exponentially many flip-cycle inequalities due to the choice of $F\subseteq A(N)$.
\end{remark}

Theorem~\ref{thm:inequalities} shows that cXPESP includes an exponential number of cuts and therefore has a large benefit over PESP and XPESP.
However, the inequalities are only true for line cycles and do not hold for arbitrary cycles in general, as is illustrated by the following example.

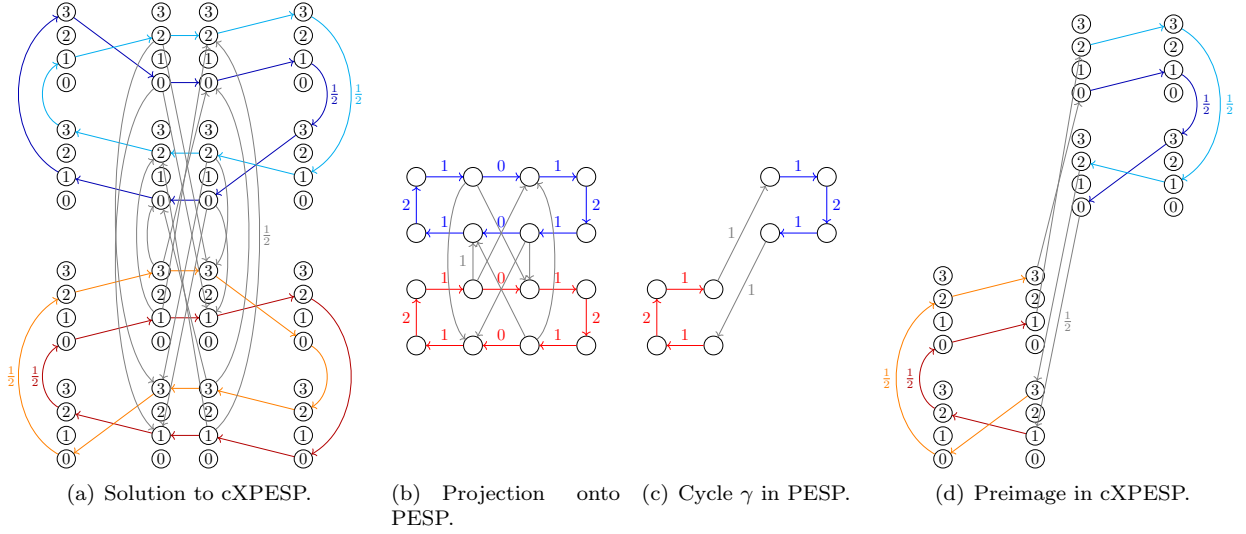
\begin{figure}
\centering
\subfigure[Solution to cXPESP.]{{\resizebox{4.9cm}{!}{
  \input{graphics/figure10}}
}}
\hfill
\subfigure[Projection onto PESP.]{\raisebox{1.5cm}{{\resizebox{2.9cm}{!}{
  \input{graphics/figure12}}
}}}
\hfill
\subfigure[Cycle $\gamma$ in PESP.]{\raisebox{1.5cm}{{\resizebox{2.9cm}{!}{
  \input{graphics/figure13}}
}}}
\hfill
\subfigure[Preimage in cXPESP.]{{\resizebox{4.9cm}{!}{
  \input{graphics/figure11}}
}}
\caption{Counterexample for Theorem~\ref{thm:splitclosure} for arbitrary cycles. }
\label{fig:example_projection}
\end{figure}

\begin{example}
	Consider the ongoing two-line example for period $T=4$. 
	Figure~\ref{fig:example_projection}a shows a possible solution to the linear programming relaxation of cXPESP and Figure~\ref{fig:example_projection}b shows the corresponding projection onto PESP with its tensions $\textbf{\textup{x}}$ . 
	In Figure~\ref{fig:example_projection}c we have chosen a cycle $\gamma$ in the event-activity network $N$ that is not a line cycle and Figure~\ref{fig:example_projection}d shows the preimage of that cycle under the projection.

    Assume that 
    $$ \gamma^\top \textbf{\textup{x}} \in \conv \left\{ \textbf{\textup{y}} \in \mathbb{Z}^{A(N)}\mid l \leq \textbf{\textup{y}} \leq u,  \frac{\gamma^T \textbf{\textup{y}}}{T} \in \mathbb{Z} \right\}.  $$
    Due to the choice of bounds (see Figure~\ref{fig:graph-comparison}b) and $T = 4$, the only feasible periodic tensions are $\textbf{\textup{y}} = l$ with $\frac{\gamma^\top l}{T} = 2$ and $\textbf{\textup{y}} = u$ with $\frac{\gamma^\top u}{T} = 3$. However, $\textbf{\textup{x}}_{\alpha} = l_{\alpha} = 2 < 3 = u_{\alpha}$ for the turnaround activities $\alpha$, and $\textbf{\textup{x}}_{\alpha} = u_{\alpha} = 1 > 0 = l_{\alpha}$ for the gray transfer activities $\alpha$, so that $\textbf{\textup{x}}$ cannot be a convex combination of $l$ and $u$. We therefore conclude that $\textbf{\textup{x}}$ violates a split inequality for $\gamma$, which translates to violating a flip-cycle inequality for $\gamma$ \citep[Theorem~3.1]{lindner_split_2025}.
    
    In fact, consider the set $F$ consisting of the two gray turnaround activities, at which  $\textbf{\textup{x}}$ is the upper bound $u$. For the flip-cycle inequality with respect to $\gamma$ and $F$, we find
    $ \xi_F = 2 $, $\textbf{\textup{x}}_{\alpha} = l_{\alpha}$ for all $\alpha \in A(\gamma) \setminus F$ and $\textbf{\textup{x}}_{\alpha} = u_{\alpha}$ for all $\alpha \in F$. Note that also $\gamma_{\alpha} = 1$ for all $\alpha \in A(\gamma)$. The flip-cycle inequality (see Lemma~\ref{lem:flipcycleinequality}) then reads as
    $$0 = (4- 2) \cdot 0 + 2 \cdot 0 \geq 2 (4 - 2) = 4,$$
    and is clearly violated.
  
    A look at the corresponding cXPESP solution in the expanded event-activity network $D$ reveals that the reason is a cycle in the cXPESP solution that includes nodes from the same event at more than one time step.
\end{example}

For the remainder of this section, the focus lies on computing a lower bound for the linear programming relaxation of cXPESP. To that end, we consider single line cycles with the help of the so-called \textit{flip polytope}.

\begin{definition}
Denote by $\BFx$ the (fractional) periodic tension corresponding to $(\pi, p) \in P_{LP}(\pesp)$ and by $\Gamma$ the set of oriented cycles.
The \textit{flip polytope} is defined as
$$P_{flip} = \{ (\pi, p) \in P_{LP}(\pesp) \mid  \text{ $\BFx$ satisfies the flip inequality for $F \subseteq A$,$\gamma \in \Gamma$} \}.$$
\end{definition}

\begin{theorem}[\citealt{lindner_determining_2020}] \label{thm:flip}
	Suppose that each activity $\alpha \in A(N)$ is contained in at most one (undirected) cycle. Then $P_{flip} = P_{IP}(\pesp)$.
\end{theorem}

Theorem~\ref{thm:flip} yields the following result for cXPESP:
\begin{lemma} \label{lem:onelinecycle}
If the event-activity network $N$ consists of exactly one line cycle, then $P_{IP}(\pesp) = \psi(\varphi(P_{LP}(\cxpesp ))) $.
\end{lemma}

\proof{Proof}
For the first inclusion, notice that 
$$P_{IP}(\pesp) = \psi(\varphi(P_{IP}(\cxpesp ))) \subseteq  \psi(\varphi(P_{LP}(\cxpesp )))$$
due to Theorem~\ref{thm_transformation_IP}.
For the other inclusion, let $(\pi, p) \in \psi(\varphi(P_{LP}(\cxpesp ))) \subseteq P_{LP}(\pesp)$ and let $\BFx$ be the corresponding, (possibly fractional) periodic tension. Denote the unique line cycle of $N$ by $\gamma$. Then the flip inequalities hold for $\BFx$ and $\gamma$ due to Theorem~\ref{thm:inequalities}. Hence, $(\pi, p) \in P_{flip}$ per definition and the result is a direct consequence of Theorem~\ref{thm:flip}.
\endproof

\begin{figure}
\centering
\subfigure[PESP IP.]{\raisebox{0.5cm}{{\resizebox{5cm}{!}{
  \input{graphics/figure14}}
}}}
\hfill
\subfigure[PESP LP.]{\raisebox{0.5cm}{{\resizebox{5cm}{!}{
  \input{graphics/figure15}}
}}}
\hfill
\subfigure[cXPESP LP.]{{\resizebox{5cm}{!}{
  \input{graphics/figure16}}
}}
\caption{cXPESP for a single line cycle. For a single line cycle, the objective value of the linear programming relaxation of cXPESP equals the objective value of PESP.}
\label{fig:singlecycle}
\end{figure}
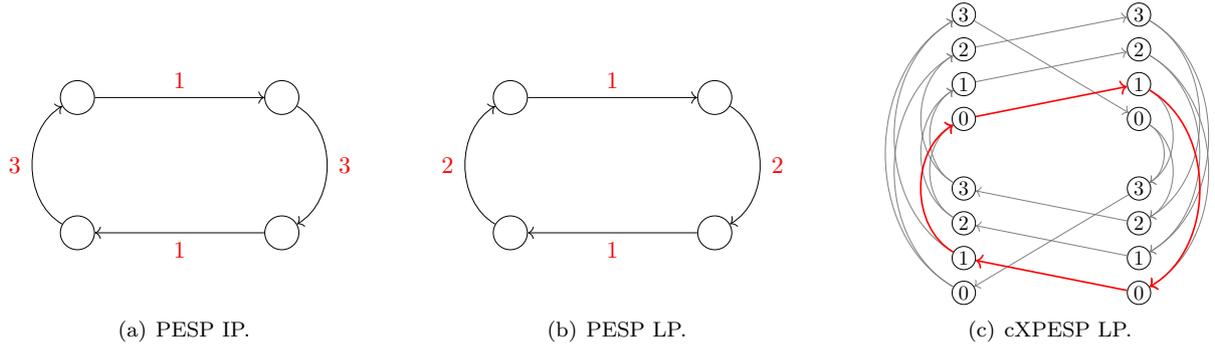

\begin{example}
Figure~\ref{fig:singlecycle} shows an examples for the statement of Lemma~\ref{lem:onelinecycle}, where the event-activity network $N$ consists of exactly one line cycle. We set $T=4$, the bounds for driving activities are fixed to 1 and the bounds of turnaround activities are $[l_\alpha,u_\alpha] = [2,3]$. An optimal solution to PESP (Figure~\ref{fig:singlecycle}a), its linear programming relaxation (Figure~\ref{fig:singlecycle}b) and the linear programming relaxation of cXPESP (Figure~\ref{fig:singlecycle}c) are colored in red. 
We denote them by $o(\text{PESP IP})$, $o(\text{PESP LP})$ and $o(\text{cXPESP LP})$, respectively.
It is
$$6 = o(\text{PESP LP}) < o(\text{PESP IP}) = o(\text{cXPESP LP}) = 8.$$

\end{example}

Lemma~\ref{lem:onelinecycle} gives an indication for a bound on the integrality gap for cXPESP:
\begin{remark}
Note that the optimal objective value of the linear programming relaxation of cXPESP 
decomposes into the contribution of the cycle variables and of the transfer arc variables, and the cycle variables can be grouped into the lines.
Since the transfer arc variables in cXPESP are modeled identically to XPESP, Lemma~\ref{lem:xpesplowerbounds} yields
$$o^\star (\cxpesp_{LP}) = \sum_{c \in C} \vartheta_c x_c + \sum_{a \in Z} \omega_a \tau_a z_a \geq \sum_{L \in \mathcal{L}} \sum_{c \in C_L} \vartheta_c x_c + \sum_{\alpha \in Z(N)} \omega_\alpha l_\alpha.$$
By Lemma~\ref{lem:onelinecycle}, the optimal value to the problem restricted to $L\in \mathcal{L}$ coincides with the optimal value of its integer solution.
Denote by $o(\cxpesp_{IP})_L$ the optimal value of the restricted problem.
Note that $o^\star (\cxpesp_{LP})$ does not have to include the optimal solution for each individual line. 
Instead $o^\star (\cxpesp_{LP})$ is bounded by them from below such that
$$o^\star (\cxpesp_{IP}) \geq o^\star (\cxpesp_{LP}) \geq  \sum_{L \in \mathcal{L}} o(\cxpesp_{IP})_L + \sum_{\alpha \in Z(N)} \omega_\alpha l_\alpha.$$
\end{remark}

\subsection{Column Generation}

Solving an integer program typically starts by solving the linear programming relaxation. While we have shown in this section that the linear programming relaxation is stronger in cXPESP in comparion to PESP and XPESP, the disadvantage of the introduced model lies in the increased number of variables. 
We deal with this increased size through the use of column generation. In the following, we discuss the pricing problems for the line cycle variables $x_c$ and the transfer arc variables $z_a$.

By the \textit{primal program}, we denote the linear programming relaxation of cXPESP, that is, cXPESP without the integer constraints \eqref{eq:int1} and \eqref{eq:int2}.
For the \textit{dual program}, we introduce a dual variable $\mu_l$ for each partitioning constraint~\eqref{eq:partition}, that is, for each $l \in \mathcal{L}$.
For each coupling constraint~\eqref{eq:c1new} at a transfer activity $(v,w) \in Z(N)$ and $t\in [T]$ introduce
a dual variable $\nu_{v[t],w}$.
For each coupling constraint~\eqref{eq:c2new} at a transfer activity $(v,w) \in Z(N)$ and $t'\in [T]$ introduce a dual variable $\rho_{v, w[t']}$.
The dual program then is
\begin{align}
	\max \quad &\sum_{L \in \mathcal{L}} \mu_L && dual \notag\\
	& \mu_L + \sum_{u \in \delta^+_{Z(N)}(v)} \sum_{v[t] \in c} \nu_{v[t], u} + \sum_{u \in \delta^-_{Z(N)}(w)} \sum_{w[t'] \in c} \rho_{u, w[t']} \leq \vartheta_c && \forall c \in C_L, \forall L \in \mathcal{L} \label{dual:c11}\\
	& - \nu_{v[t],w} - \rho_{v,w[t']}  \leq \omega_a \tau_a && \forall a = (v[t],w[t'] )  \label{dual:c22}.
\end{align}

\subsubsection{Pricing cycle variables} 
\label{sss:pricing_cycle}

The aim is to find a cycle that violates constraint~\eqref{dual:c11}, i.e., to find $c \in C$ such that

\begin{align*}
	&\mu_L > \vartheta_c - \sum_{u \in \delta^+_{Z(N)}(v)} \sum_{v[t] \in c} \nu_{v[t], u} -\sum_{u \in \delta^-_{Z(N)}(w)} \sum_{w[t'] \in c} \rho_{u, w[t']}.
\end{align*}
This could be solved for each $L \in \mathcal{L}$ individually by solving
\begin{align}
	& \min_{c \in C_L} \quad \vartheta_c - \sum_{u \in \delta^+_{Z(N)}(v)} \sum_{v[t] \in c} \nu_{v[t], u} -\sum_{u \in \delta^-_{Z(N)}(w)} \sum_{w[t'] \in c} \rho_{u, w[t']}\label{eq:reccost1}
\end{align}
and checking if the optimal value is smaller than $\mu_L$. 

\begin{lemma} \label{cycle_pricing}
The pricing problem for a cycle variable $x_c$ in the expanded line cycle of a given line $L$ in cXPESP is a set of $T$ shortest path problems in an acyclic graph.
\end{lemma}

\begin{figure}
\centering
\subfigure[]{{\resizebox{8cm}{!}{
  \input{graphics/figure17}}
}}
\hfill
\subfigure[]{{\resizebox{8cm}{!}{
  \input{graphics/figure18}}
}}
\caption{Cutting a graph at an arbitrary event.}
 \label{fig:cut}
\end{figure}
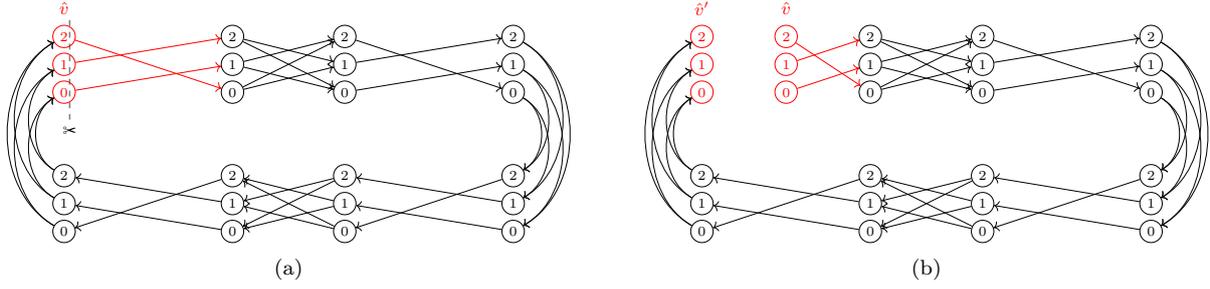

\proof{Proof}
Let $L$ be a line in $G$. 
The aim is to find new cycles in the expanded line cycle of $L$.
Recall cycles to be of fixed length to avoid closed paths that pass the same event more than once, compare Figure~\ref{fig:cycle}.
We avoid the longer closed paths by cutting the expanded line cycles:
Choose an arbitrary event $\hat{v} \in V(N)$, duplicate all nodes $\hat{v}[t]$ for $t \in [T]$, denote them by $\hat{v}[t]'$, and thereby cut the expanded line cycle as depicted in Figure~\ref{fig:cut}. The resulting graph $C_{L, \hat{v}}$ is acyclic.

For each arc $a = (v[t], w[t'])$ in the subgraph $C_L$ of the expanded event-activity network $D$, define the reduced arc cost as
\begin{align*}
	\overline{c}_a &\coloneqq \omega_{a} \tau_{a} -  \sum_{u \in \delta^+_{Z(N)}(v)} \nu_{v[t],u} - \sum_{u \in \delta^-_{Z(N)}(w)} \rho_{u, w[t']}, 
\end{align*}
which represents the objective in \eqref{eq:reccost1} restricted to each arc of a cycle.
Using the reduced arc costs $\overline{c}_a $ for each $a \in X(C_L)$ as weight, we can find a cycle passing through $v[t]$ for some $t\in [T]$ by solving a shortest path problem in $C_{L, \hat{v}}$ with $\hat{v}[t]$ as source and $\hat{v}'[t]$ as target.
Since we want to check the cycles for each $t\in [T]$, we need to solve $T$ shortest path problems per line.
\endproof

Since  $C_{L, \hat{v}}$ is acyclic, we can apply topological search to find shortest paths, and hence solve the pricing problem for line cycles of line $L$ within a time complexity of
$$\mathcal O(T \cdot (|V(C_{L, \hat{v}})| + |A(C_{L, \hat{v}})| )) = \mathcal O(T \cdot |A(C_L)|) . $$

\subsubsection{Pricing for transfer variables}

The aim is to find a transfer arc, that violates Constraint~\eqref{dual:c22}, i.e. find $a = (v[t], w[t']) \in Z$ such that
$$- \nu_{v[t],w} - \rho_{v,w[t']} > \omega_a \tau_a.$$
This could be solved by 
$$\min_{a=(v[t],w[t'])} \quad \omega_a \tau_a + \nu_{v[t],w} + \rho_{v,w[t']},$$
which could be simply approached by enumerating the $T \cdot (u_\alpha - l_\alpha \mod T)$ arcs per transfer activity $\alpha \in Z(N)$. 
This sums up to $T \cdot \sum_{\alpha \in Z(N)} (u_\alpha - l_\alpha \mod T)$, so that there are in total $\mathcal{O}(T^2)$ many arcs to enumerate for each transfer activity.

\section{$\cxpespw$: Linearizing the Number of Expanded Transfer Arcs}

While cXPESP has beneficial theoretical properties, one of its drawbacks is that each transfer activity $\alpha \in Z(N)$ produces up to $T^2$ expanded transfer arcs in $Z(D)$. We therefore formulate a variation of cXPESP, called $\cxpespw$, that only introduces $T$ transfer arcs, at the cost of introducing $T$ additional \emph{waiting arcs}, that link departure nodes. This has already been suggested by \citet{BorndoerferHoppmannKarbstein2017}.

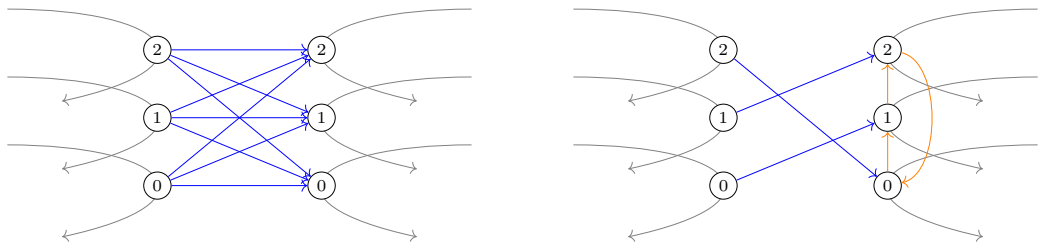
\begin{figure}[h!]
\centering
\subfigure[The standard expanded event-activity network at a transfer activity, with $T^2 = 9$ expanded transfer arcs in blue.]{
 \scalebox{0.9}{
    \begin{tikzpicture}[xscale=0.8, node distance={10mm}, font=\scriptsize,  inner sep=1pt, minimum size=4mm, time/.style = {}, invisi/.style = {}]
    \tikzstyle{main} = [draw, circle, fill=white]
    
     \node[invisi] (i1) at (7,1.2) {};
     \node[invisi] (i2) at (7,0.2) {};
     \node[invisi] (i3) at (7,-0.8) {};
     \node[invisi] (i4) at (6,2.6) {};
     \node[invisi] (i5) at (6,1.6) {};
     \node[invisi] (i6) at (6,0.6) {};
     
     \node[invisi] (j1) at (14,1.2) {};
     \node[invisi] (j2) at (14,0.2) {};
     \node[invisi] (j3) at (14,-0.8) {};
     \node[invisi] (j4) at (15,2.6) {};
     \node[invisi] (j5) at (15,1.6) {};
     \node[invisi] (j6) at (15,0.6) {};   
    
    \draw[gray, ->] (i4) to [out=0, in=10, looseness=4.5] (i1);
    \draw[gray, ->] (i5) to [out=0, in=10, looseness=4.5] (i2);
    \draw[gray, ->] (i6) to [out=0, in=10, looseness=4.5] (i3);
    
    \draw[gray, ->] (j4) to [out=180, in=170, looseness=4.5] (j1);
    \draw[gray, ->] (j5) to [out=180, in=170, looseness=4.5] (j2);
    \draw[gray, ->] (j6) to [out=180, in=170, looseness=4.5] (j3);
    
    \node[main] (4) at (9,2) {$2$};
    \node[main] (5) [below of=4] {$1$}; 
    \node[main] (6) [below of=5]  {$0$};
    \node[main] (10) at (12,2) {$2$};
    \node[main] (11) [below of=10]  {$1$};
    \node[main] (12) [below of=11]  {$0$};
    
     \draw[blue, ->] (4) --  (10);
    \draw[blue, ->] (4) --  (11);
    \draw[blue, ->] (4) --  (12);
    \draw[blue, ->] (5) --  (10);
    \draw[blue, ->] (5) --  (11);
    \draw[blue, ->] (5) --  (12);
    \draw[blue, ->] (6) --  (10);
    \draw[blue, ->] (6) --  (11);
    \draw[blue, ->] (6) --  (12);
    
    \end{tikzpicture}}}
\hspace{1ex}
\subfigure[The network at a transfer activity $\alpha$ with $l_\alpha = 1$ for the waiting arc transfer model. The $T = 3$ transfer arcs are blue, the $T=3$ waiting arcs are orange.]{
     \scalebox{0.9}{
    \begin{tikzpicture}[xscale=0.8, node distance={10mm}, font=\scriptsize, inner sep=1pt, minimum size=4mm, time/.style = {}, invisi/.style = {}]
    \tikzstyle{main} = [draw, circle, fill=white]    
     \node[invisi] (i1) at (7,1.2) {};
     \node[invisi] (i2) at (7,0.2) {};
     \node[invisi] (i3) at (7,-0.8) {};
     \node[invisi] (i4) at (6,2.6) {};
     \node[invisi] (i5) at (6,1.6) {};
     \node[invisi] (i6) at (6,0.6) {};
     
     \node[invisi] (j1) at (14,1.2) {};
     \node[invisi] (j2) at (14,0.2) {};
     \node[invisi] (j3) at (14,-0.8) {};
     \node[invisi] (j4) at (15,2.6) {};
     \node[invisi] (j5) at (15,1.6) {};
     \node[invisi] (j6) at (15,0.6) {};

     \draw[gray, ->] (i4) to [out=0, in=10, looseness=4.5] (i1);
    \draw[gray, ->] (i5) to [out=0, in=10, looseness=4.5] (i2);
    \draw[gray, ->] (i6) to [out=0, in=10, looseness=4.5] (i3);
    
    \draw[gray, ->] (j4) to [out=180, in=170, looseness=4.5] (j1);
    \draw[gray, ->] (j5) to [out=180, in=170, looseness=4.5] (j2);
    \draw[gray, ->] (j6) to [out=180, in=170, looseness=4.5] (j3);

    \node[main] (4) at (9,2) {$2$};
    \node[main] (5) [below of=4] {$1$}; 
    \node[main] (6) [below of=5]  {$0$};
    \node[main] (10) at (12,2) {$2$};
    \node[main] (11) [below of=10]  {$1$};
    \node[main] (12) [below of=11]  {$0$};

    \draw[blue, ->] (4) --  (12);
    \draw[blue, ->] (5) --  (10);
    \draw[blue, ->] (6) --  (11);
    \draw[orange, ->] (12) --  (11);
    \draw[orange, ->] (11) --  (10);
    \draw[orange, ->] (10) to [out=350, in=10, looseness=1]  (12);
    
    \end{tikzpicture}}}
    \caption{Graph structure for different transfer models, visualized for $T=3$.}
    \label{fig:cxpespw-transfers}
\end{figure}

The idea is visualized in Figure~\ref{fig:cxpespw-transfers}. For each transfer activity $\alpha = (v, w) \in Z(N)$ in the PESP instance, we only construct the $T$ expanded arcs $(v[t], w[(t + l_\alpha) \bmod T])$ for $t \in [T]$. To model that a transfer activity might take longer than its lower bound, we add waiting arcs $(w[t], w[(t+1) \bmod T])$ for $t \in [T]$ that allow waiting at the expanded departure nodes. We set $\tau_a \coloneqq 1$ and $\omega_a \coloneqq \tau_\alpha$ for each waiting arc $a$ arising from a transfer activity $\alpha$. This way, we can model a transfer from $v$ to $w$ at time $t$ of duration $\ell_\alpha + k$ by a path comprised of the expanded transfer arc $(v[t], w[t + l_\alpha])$, and the $k$ waiting arcs $$(w[t + l_\alpha], w[t + l_\alpha+1]), \dots, (w[t + l_\alpha+k-1], w[t + l_\alpha+k]).$$
We omit the modulo $T$ expressions in the time indices for better readability. We denote the union of the set of reduced transfer activities and the set of waiting arcs by $Z^w(D)$. The modified model, that we call $\cxpespw$, is as follows:
\begin{align} 
	\min & \sum_{c \in C} \omega_c \tau_c x_c + \sum_{a \in Z^w(D)} \omega_a \tau_a z_a  & \text{cXPESP}^{w} \notag\\ 
	& \sum_{c \in C_L} x_c = 1 & \forall L \in \mathcal{L} \label{eq:bwpart}\\
	& \sum_{c: v[t] \in c} x_{c} -  z_{(v[t], w[t])}  =  0 & \forall \alpha = (v,w) \in Z(N), t \in [T] \label{eq:bwc1} \\
	&  \sum_{c: w[t] \in c} x_{c} + z_{(w[t], w[t+1])} -   z_{(w[t-1], w[t])}   -  z_{(v[t-l_{\alpha}], w[t])} =  0 & \forall \alpha = (v,w) \in Z(N), t \in [T] \label{eq:bwc2}\\
	& \sum_{t = 0} ^{T-1} z_{w[t],w[t+1]} \leq u_\alpha - l_\alpha &  \forall \alpha = (v,w) \in Z(N) \label{eq:bwbw}\\
	& x_c   \geq 0& \forall c \in C \label{eq:bwb1}\\
	& z_a   \geq 0& \forall a \in Z^w(D) \label{eq:bwb2}\\
	& x_c \in \mathbb{Z}  & \forall c \in C\\
	& z_a \in \mathbb{Z}  & \forall a \in Z^w(D) \label{eq:bwi2}.
\end{align}
Comparing to the $\cxpesp$ model, we can replace the sum over the transfer arcs in \eqref{eq:c1new} by a single variable to obtain \eqref{eq:bwc1}. Constraint \eqref{eq:bwc2} ensures flow conservation, replacing \eqref{eq:c2new}. The constraint \eqref{eq:bwbw} ensures that transferring and waiting do not exceed the upper bound $u_\alpha$ of the original transfer activity $\alpha$. This constraint can be removed whenever $\alpha$ is free, i.e., $u_\alpha - l_\alpha \geq T-1$. In this case, we obtain an equally strong linear programming relaxation:

\begin{theorem}
\label{thm:cxpespw-relaxation}
    Let $o(\cxpesp\text{ LP})$ and $o(\cxpespw\text{ LP})$ denote the optimal values of the LP relaxations of $\cxpesp$ and $\cxpespw$, respectively. Then $o(\cxpesp\text{ LP}) \geq o(\cxpespw\text{ LP})$ and equality holds if all transfer activities are free.
\end{theorem}

\proof{Proof}
For a given transfer activity $\alpha = (v,w) \in Z(N)$, we consider the expanded transfer subgraph as a flow network: Each $v[t]$ represents a source and $w[t']$ a sink, so that the $z$-variables describe a multi-commodity flow.

Now let $(x, z) \in P_\text{LP}(\cxpesp)$ be optimal. To construct a point of $P_\text{LP}(\cxpespw)$, the polyhedron associated to the LP relaxation of $\cxpespw$, we proceed as follows: We leave $x$ unchanged. For all transfer arcs $a = (v[t], w[t']) \in Z(D)\setminus Z^w(D)$, we add $z_a$ units of flow on the $v[t]$-$w[t']$-path consisting of the transfer arc $(v[t], w[t+l_\alpha])$ and the $w[t+l_\alpha]$-$w[t']$-path along $\tau_a - l_\alpha$ waiting arcs. This procedure leaves the objective value unchanged, conserves the flow \eqref{eq:bwc1} and \eqref{eq:bwc2}, and adheres to upper bounds \eqref{eq:bwbw}. We conclude $o(\cxpesp\text{ LP}) \geq o(\cxpespw\text{ LP})$.

To prove equality when $\alpha$ is free, let $(x^w,z^w) \in P_\text{LP}(\cxpespw)$ be optimal. We construct $(x,z) \in P_\text{LP}(\cxpesp)$ with the same objective value. We apply flow decomposition to $z^w$. As $(x^w,z^w)$ is optimal, we can assume that this decomposition contains no cycles, but only $v[t]$-$w[t']$ paths. For each of these paths, we add its amount of flow to the value of $z_{v[t],w[t']}$. We can guarantee that the expanded transfer arc $(v[t],w[t'])$ exists, since $\alpha$ is free. Again, the objective value is not modified by this procedure.

\endproof

\begin{remark}
    In the context of column generation, the pricing problem for the cycle variables remains unchanged, and can be solved with the same strategies as described in Section~\ref{sss:pricing_cycle}. 
\end{remark}

\section{cXTTP: Integrated Periodic Timetabling and Passenger Routing} \label{section:integration}

Since timetabling and passenger routing affect each other in the process of optimizing public transport, we will apply the integration of timetabling and passenger routing to cXPESP in this section.

\subsection{Model Description}

We first recall PESP and XPESP with integrated passenger routing. 
Given a graph $(V,A)$ an \textit{origin-destination-matrix (OD-matrix)} is a $V \times V$-matrix $(d_{st})_{(s,t) \in V \times V}$, where $d_{st}$ is a non-negative integer that describes the demand from node $s$ to node $t$. 
An \textit{OD-pair} is a tuple of nodes $(s,t) \in V \times V$ such that $d_{st} >0$.
We denote the set of OD-pairs by $\mathcal{OD}$.
Furthermore, $P_{st}$ denotes the set of possible paths between nodes $s$ and $t$.
Note that it is not always clear which line is taken by a passenger traveling from a starting location $s$ of an OD-pair. 
It is, therefore, recommended to add artificial nodes $v_{s}$ and $v_{t}$ and corresponding arcs to the underlying network for possible OD-pair nodes. 

In addition to the variables in PESP and XPESP, we introduce for each $(s,t) \in \mathcal{OD}$ and for each $s$-$t$-path $p \in P_{st}$ the variable $y_p \in \mathbb{Q}^{+}$, representing the passenger flow on path $p$.
For PESP with integrated passenger routing, we use a mixed integer programming formulation inspired by \cite{BorndoerferHoppmannKarbstein2017}:
\begin{align}
	\min & \sum_{st \in \mathcal{OD}} \sum_{p \in P_{st}} \sum_{a = (v,w) \in p} d_{st} y_p \omega_a (\pi_w - \pi_v + T p_a) & \text{TTP}	\notag\\
	&\pi_w - \pi_v + T p_a \geq l_a  & \forall a= (v,w) \in A \notag\\
	&\pi_w - \pi_v + T p_a \leq u_a & \forall a= (v,w) \in A \notag\\
	&0 \leq \pi_v   \leq T - 1& \forall  v \in V \notag\\
	& \pi_v \in \mathbb{Z}  & \forall v \in V \notag\\
	& p_a \in \mathbb{Z}  & \forall a \in A \notag\\
	& \sum_{p \in P_{st}} y_p = 1 & \forall (s,t) \in \mathcal{OD} \label{eq:ttp6}\\
	& y_p \geq 0 & \forall p \in P_{st}, \forall (s,t) \in \mathcal{OD}. \notag
\end{align}

This mixed integer program is called TTP as suggested in \cite{BorndoerferHoppmannKarbstein2017}, where also another formulation based on cycle bases is introduced.
Here we use a formulation that is more similar to the PESP version we used before.
In addition to the known constraints from PESP, we add Constraint~\eqref{eq:ttp6} to model a total passenger flow of one from node $s$ to node $t$. 
Note that the objective of TTP is not linear.

Now denote for each path $p \in P_{st}$ by $\tau_{p} = \sum_{a \in p} \tau_a$ its duration.
For XPESP with integrated passenger routing (XTTP), we use the formulation based on \cite{BorndoerferHoppmannKarbstein2017}:
\begin{align} 
	\min & \sum_{st \in \mathcal{OD}} \sum_{p \in P_{st}} d_{st} \tau_p y_p  + \sum_{a \in X(D)} \omega_a \tau_a x_a &\text{XTTP}\notag \\
	& \sum_{a \in \mathcal{A}(\alpha)} x_a= 1 & \forall \alpha \in X(D) \label{eq:xttp1}\\
	& \sum_{a \in \delta_X^-(v) } x_{a} - \sum_{a \in \delta_X^+(v) } x_{a}  = 0 & \forall v \in V(D) \label{eq:xttp2} \\ 
	& x_a - \sum_{p \in P_{st}: a \in p} y_p \geq 0 & \forall (s,t) \in \mathcal{OD}, \forall a \in X(D)  \label{eq:xttp3}\\
	& \sum_{p \in P_{st}} y_p = 1 & \forall (s,t) \in \mathcal{OD} \label{eq:xttp4} \\
	&0 \leq x_a   \leq 1& \forall  a\in X(D) \label{eq:xttp5} \\
	& x_a \in \mathbb{Z}  & \forall a \in X(D) \label{eq:xttp6} \\
	& y_p \geq 0 & \forall p \in P_{st}, \forall (s, t) \in \mathcal{OD}. \label{eq:xttp7} 
\end{align}
In contrast to the formulation in \cite{BorndoerferHoppmannKarbstein2017}, we omit the transfer arcs, since the relevant transfers are already included in the passenger paths. 
Furthermore, we want to optimize not only the passenger flow variables but additionally the duration on the line cycles. 
We therefore add the arc variables to the objective.
Moreover, we add a coupling Constraint~\eqref{eq:xttp3}, such that there can be passenger flow on an arc only if the arc belongs to the subgraph of the solution.

The approach of integrating passenger routing can be transferred to cXPESP.
Then we call cXPESP with integrated passenger routing cXTTP defined by the following mixed integer programming formulation:
\begin{align}
	\min & \sum_{st \in D} \sum_{p \in P_{st}} d_{st} \tau_p y_p + \sum_{c \in C} \vartheta_c x_c & \text{cXTTP} \notag\\
	& \sum_{c \in C_L} x_c = 1 & \forall L \in \mathcal{L} \label{eq:pr1}\\
	& \sum_{p \in P_{st}} y_p = 1 & \forall (s,t) \in \mathcal{OD} \label{eq:pr2}\\
	&\sum_{c \in C:a \in c} x_c - \sum_{p\in P_{st}:a \in p} y_p \geq 0 & \forall a \in X(D),  \forall (s,t) \in \mathcal{OD} \label{eq:pr3}\\
	& x_c  \geq 0 & \forall c \in C \label{eq:pr4}\\
	& y_p \geq 0 & \forall p \in P_{st}, \forall (s,t) \in \mathcal{OD} \label{eq:pr5}\\
	& x_c \in \mathbb{Z}  & \forall c \in C. \label{eq:pr6}
\end{align}

Constraints~\eqref{eq:pr1} is inherited from cXPESP.  Constraint~\eqref{eq:pr2} is a partitioning constraint that defines the total flow between the nodes of an OD-pair to be one.
The aim of Constraint~\eqref{eq:pr3} is again a coupling between variables. 
There should only be a positive passenger flow on an arc if the arc is part of the resulting subgraph.
We want to minimize the duration both of passenger paths and of cycles.
Notice that we omit again the transfer variables, but a minimization of path durations automatically is a minimization of transfer durations.
Thus, we obtain a comparability of the objectives of XTTP and cXTTP.

\begin{remark}
    Relaxing the integrality constraints in the $\ttp$ model yields a quadratic program, whose optimal objective value is always the weighted sum taken over all OD-pairs $(s, t) \in \mathcal{OD}$, where each summand is obtained as demand $d_{st}$ times the cost of a shortest $s$-$t$-path w.r.t.\ the lower bounds $l$. Since cXTTP inherits the cycle variables from cXPESP, we expect by Example~\ref{ex:notzero} better LP relaxations as well.
\end{remark}

\subsection{Comparison of Solution Polytopes}

Since the objective of TTP is not linear,
we restrict the comparison between the integrated models to XTTP and cXTTP.

\begin{definition}
Denote by
\begin{align*}
	& P_{MIP}(\text{XTTP}) = \conv \{ (x, y) \in \mathbb{Z}^{X(D)} \times \mathbb{Q}^{P} | (x, y) \text{ satisfies } (\ref{eq:xttp1}) - (\ref{eq:xttp7})\},\\
	& P_{LP}(\text{XTTP}) = \{ (x, y) \in \mathbb{Q}^{X(D)} \times \mathbb{Q}^{P} | (x, y) \text{ satisfies } (\ref{eq:xttp1}) - (\ref{eq:xttp5}), (\ref{eq:xttp7})\}, \\
	& P_{MIP}(\text{cXTTP}) = \conv \{ (x, y)\in \mathbb{Z}^C \times \mathbb{Q}^{P} | (x, y)\text{ satisfies } (\ref{eq:pr1}) - (\ref{eq:pr6})\},\\
	& P_{LP}(\text{cXTTP}) =  \{ (x, y) \in \mathbb{Q}^C \times \mathbb{Q}^{P} | (x, y)\text{ satisfies } (\ref{eq:pr1}) - (\ref{eq:pr5})\}
\end{align*}
the solution spaces of the integer program and linear program relaxation for TTP, XTTP and cXTTP. 
Define the linear transformation
		\begin{align*}
			\phi : \mathbb{Q}^C \times \mathbb{Q}^P & \rightarrow \mathbb{Q}^{X(D)}\times \mathbb{Q}^{P}\\
			 \begin{pmatrix}{c} x \\ y \end{pmatrix} & \mapsto \begin{pmatrix}{cc} M_C & 0\\  0 &  I \end{pmatrix} \cdot \begin{pmatrix}{c} x \\ y \end{pmatrix}, 
		\end{align*}
		where 
		\begin{align*}
			&M_C = (m_{ac})_{a \in X(D), c \in C},
			& m_{ac} = \begin{cases} 1 &\text{if $a \in c$,}\\ 0 &\text{otherwise.} \end{cases} 
		\end{align*}
\end{definition}

As for cXPESP, the following theorem shows that cXTTP has a tighter linear programming relaxation as XTTP.

\begin{theorem} \label{thm:polycxttp}
The linear transformation $\phi$
has the property:
\begin{align*}
    \phi(P_{MIP}(\cxttp )) &= P_{MIP}( \xttp), \\
    \phi(P_{LP}(\cxttp )) &\subseteq P_{LP}( \xttp).
\end{align*}
\end{theorem}

\proof{Proof}
	 $\phi(P_{LP}(\cxttp )) \subseteq P_{LP}( \xttp)$:\\
		Let $( x^{\circ} , y^{\circ}) \in P_{LP}(\cxttp )$. We show that $(\bar{x}, \bar{y}) = \phi (x^{\circ} , y^{\circ}) \in P_{LP}(\xttp )$.
		Note that this transformation equals the transformation in Definition~\ref{def:transformation} when both are restricted to $\mathbb{Q}^C$. 
		Thus, Constraints~\eqref{eq:xttp1},  \eqref{eq:xttp2}, and  \eqref{eq:xttp5} are already proven.
		Furthermore, $\phi$ restricted to $\mathbb{Q}^P$ is the identity and, hence, Constraints~\eqref{eq:xttp4} and \eqref{eq:xttp7} hold. 
		It remains to show Constraint~\eqref{eq:xttp3}: Let $a \in X(D)$ and $(s,t) \in \mathcal{OD}$. Then 
		$$\bar{x}_a = \sum_{c \in C:a \in c} x^{\circ}_c \overset{(\ref{eq:pr3})}{\geq} \sum_{p \in P_{st}: a \in p} y^{\circ}_p = \sum_{p \in P_{st}: a \in p} \bar{y}_p. $$
	$\phi(P_{MIP}(\cxttp )) = P_{MIP}( \xttp)$ follows from the proof of Theorem~\ref{thm_transformation_IP}.	
 
\endproof

\begin{remark}
	In general, the inclusion in Theorem~\ref{thm:polycxttp} is not an equality. Consider again the example in Remark~\ref{rem:incl}, 
	where the set of OD-pairs is empty.
\end{remark}

\subsection{Column Generation}

The advantage of cXTTP in comparison to XTTP lies in its possibly tighter linear programming relaxation.
As for cXPESP, this is again rooted in the richer structure of the formulation, which includes information about the lines.
However, the same disadvantage appears in the huge amount of cycle variables.
Hence, here too, we consider column generation to deal with the large number of variables.

Consider the \textit{primal program} to be cXTTP with relaxed integrality constraint, that is, omit Constraints~\eqref{eq:pr6}.
For the \textit{dual program}, we introduce for each $L \in \mathcal{L}$ a dual variable $\mu_L$, for $(s,t) \in \mathcal{OD}$ a dual variable $\nu_{st}$ and for each $(s,t) \in \mathcal{OD}$ and $a \in X $ a dual variable $\rho_a^{st}$. 
Then the dual linear program reads
\begin{align}
	\max \quad &\sum_{L \in \mathcal{L}} \mu_L + \sum_{(s,t) \in \mathcal{OD}} \nu_{st} & dual \notag\\
	& \mu_L + \sum_{(s,t) \in \mathcal{OD}} \sum_{a \in c} \rho_a^{st} \leq \vartheta_c & \forall c \in C_L, \forall L \in \mathcal{L} \label{eq_dual:cycle}\\
	&  \nu_{st}  -  \sum_{a \in p} \rho_a^{st}  \leq d_{st} \tau_p & \forall p \in P_{st}, \forall (s,t) \in \mathcal{OD} \label{eq_dual:flow}\\
	& \rho_a^{st}  \geq 0 & \forall a \in X,  \forall (s,t) \in \mathcal{OD} \notag.
\end{align}

\subsubsection{Pricing Cycle Variables}

For pricing cycle variables, find a cycle that violates Constraint~\eqref{eq_dual:cycle}, that is, find $c \in C$ such that  
$$\mu_L > \vartheta_c - \sum_{(s,t) \in \mathcal{OD}} \sum_{a \in c} \rho_a^{st}.$$
This can be solved for each $L$ individually by
$$\min_{c \in C_L} \quad \vartheta_c - \sum_{(s,t) \in \mathcal{OD}} \sum_{a \in c} \rho_a^{st}$$
and checking if the minimal value is smaller than $\mu_L$.
Define for each arc $a$ in the subgraph $C_L$
the reduced costs by
$$\overline{c}_a \coloneqq \omega_a \tau_a - \sum_{(s,t) \in \mathcal{OD}} \rho_a^{st}.$$ 

Then the pricing problem can be solved as a sequence of $T$ shortest path problems on a directed acyclic graph as in Lemma~\ref{cycle_pricing}.

\subsubsection{Pricing Passenger Flow Variables}
\label{sec:pricing_passenger_flow}
For pricing passenger flow variables, find a passenger flow path that violates Constraint~\eqref{eq_dual:flow}, that is, find a path $p$ such that
 $$\nu_{st}   > d_{s,t} \tau_p  + \sum_{a \in p} \rho_a^{st},$$
which could be solved for each $(s,t) \in D$ individually by 
 $$\min_{p \in P_{s,t}} d_{s,t} \tau_p  + \sum_{a \in p} \rho_a^{st}.$$
Note that $s$ and $t$ are determined by $p \in P$. 
Then define for each arc $a \in p$ the reduced costs by
$$\overline{c}_a \coloneqq d_{st}\tau_a + \rho_a^{st}.$$
This is again a shortest path problem, but with the disadvantage that solving the problem for each $(s,t) \in \mathcal{OD}$ still involves the whole expanded event-activity network, in contrast to the pricing of cycle variables, where the problem only makes use of an acyclic subgraph. 
As the reduced costs are non-negative, 
the shortest path problem can be solved by the Dijkstra algorithm with a time complexity of $\mathcal{O}(|V(D)| \log |V(D)| + |A(D)|)$.
Observe that, different from the cycle variables, it is not necessary to solve the pricing problem for each $t \in [T]$.

\section{Computational Experiments} \label{section:computation}

In this section, we will assess the computational power of the optimization models presented in Section~\ref{section:cxpesp} and Section~\ref{section:integration}. We will first describe our instances in Section~\ref{ss:instances}, then describe our experimental setup in Section~\ref{ss:setup}, and finally evaluate the cXPESP and cXTTP models in Section~\ref{ss:results}.

\subsection{Instances}
\label{ss:instances}

\begin{figure}
\centering
\hfill
\subfigure[2linecross]{\raisebox{0.8cm}{{\resizebox{3cm}{!}{
  \input{graphics/figure20}}
}}} 
\hfill
\subfigure[3berlin]{{\resizebox{0.5\linewidth}{!}{
  \input{graphics/figure21}}
}}
\hfill
\subfigure[berlin]{{\resizebox{0.48\linewidth}{!}{
  \input{graphics/figure22}}
}}
\hfill
\subfigure[R1L1-10]{{\resizebox{0.48\linewidth}{!}{
  \includegraphics{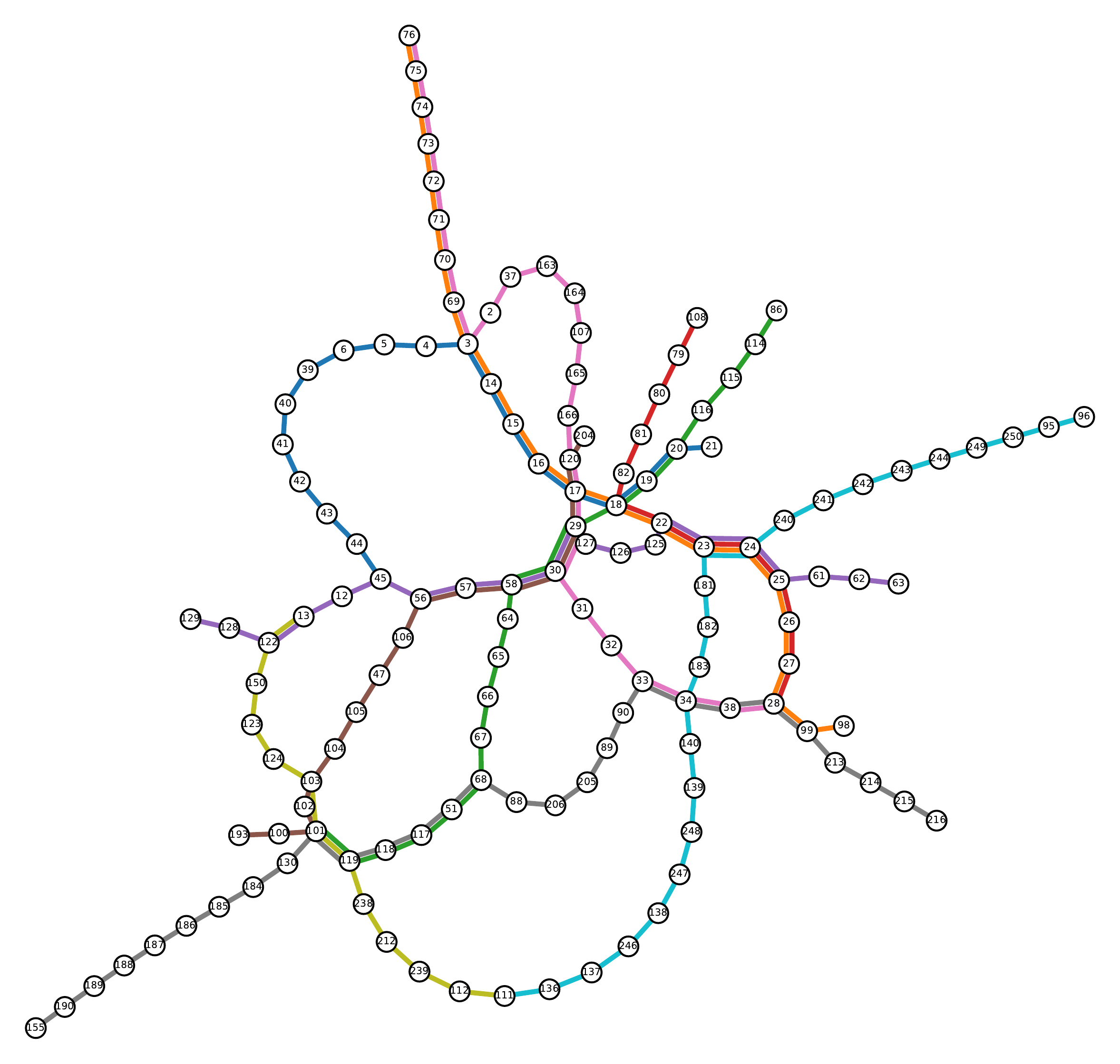}}
}}
\hfill
\caption{Line networks of the instances.} 
\label{fig:instance_networks}
\end{figure}

We consider four sets of instances: \emph{2linecross}, \emph{3berlin}, \emph{berlin}, and \emph{R1L1}. The instance 2linecross is a toy instance with 2 lines. The instances 3berlin and berlin are derived from the Berlin subway network, where 3berlin is a restriction to 3 lines, and berlin is the full network. For 2linecross and 3berlin, we consider varying period times from 5 to 60, and for berlin, there is a version with $T=5$ and one with $T = 10$. Finally, the R1L1 instances are subinstances of the smallest PESPlib \cite{} instance, comprising 1, 2, 5, and 10 lines according to the sorting procedure described in \cite{lindner_incremental_2023}. We consider the R1L1 instances with their original period time $T = 60$.
The line networks of the instances are depicted in Figure~\ref{fig:instance_networks}, and some characteristics are collected in Table~\ref{tab:instances}.

\begin{table}
\centering
    \caption{Instance characteristics}
    \label{tab:instances}
    \begin{tabular}{lr|rr|rr|rrr}
        & & \multicolumn{2}{c|}{Line Network $G$} & \multicolumn{2}{c|}{Event-Activity Network $N$} & \multicolumn{3}{
        c}{Expanded Network $D$} \\
         Instance & Period & Lines & Stations & Events & Activities & Nodes & Arcs & $\text{Arcs}^w$ \\ 
         & $T$ & $|\mathcal L|$ & $|V(G)|$ & $|V(N)|$ & $|A(N)|$ & $|V(D)|$ & $|A(D)|$ & $|A(D)|$    \\ 
         \hline  
        2linecross & 5 & \multirow{8}{*}{2} & \multirow{8}{*}{5} & \multirow{8}{*}{8} & \multirow{8}{*}{16} & 40 & 360 & 200   \\
        2linecross & 10 & & & & & 80 & 1320 & 600 \\
        2linecross & 15 & & & & & 120 & 2880 & 1200 \\
        2linecross & 20 & & & & & 160 & 5040 & 2000 \\
        2linecross & 30 & & & & & 240 & 11160 & 4200 \\
        2linecross & 40 & & & & & 320 & 19680 & 7200 \\
        2linecross & 50 & & & & & 400 & 30600 & 11000 \\
        2linecross & 60 & & & & & 480 & 43920 & 15600 \\ 
        \hline
        3berlin & 5 & \multirow{8}{*}{3} & \multirow{8}{*}{9} & \multirow{8}{*}{36} & \multirow{8}{*}{60} & 120 & 960 & 480 \\
        3berlin & 10 & & & & & 240 & 3420 & 1260 \\
        3berlin & 15 & & & & & 360 & 7380 & 2340\\
        3berlin & 20 & & & & & 480 & 12840 & 3720 \\
        3berlin & 30 & & & & & 720 & 28260 & 7380\\
        3berlin & 40 & & && & 960 & 49680 & 12240\\
        3berlin & 50 & & & & & 1200 & 77100 & 18300 \\
        3berlin & 60 & & & & & 1440 & 110520 & 25560\\ 
        \hline
        berlin & 5 & \multirow{2}{*}{9} & \multirow{2}{*}{35} & \multirow{2}{*}{220} & \multirow{2}{*}{502} & 960 & 9290 & 3650 \\
        berlin & 10 & & & & & 1920 & 33580 & 8200 \\ 
        \hline
        R1L1-1 & 60 & 1 & 26 & 100 & 100 & 120 & 3600 & 3600 \\
        R1L1-2 & 60 & 2 & 43 & 172 & 183 & 900 & 63180 & 24240 \\
        R1L1-5 & 60 & 5 & 79 & 360 & 424 & 4920 & 302400 & 75840 \\
        R1L1-10 & 60 & 10 & 139 & 718 & 1013 & 13080 & 1014660 & 189840
    \end{tabular}
\end{table}

We finally remark that all considered instances have exclusively free transfer activities, so that the LP relaxations of $\cxpesp$ and $\cxpespw$ will have the same optimal values by Theorem~\ref{thm:cxpespw-relaxation}.

\subsection{Experimental Setup}
\label{ss:setup}

We have implemented our new models $\cxpesp$, $\cxpespw$, and $\cxttp$ inside the ConcurrentPESP framework \citep{BorndoerferLindnerRoth2020}, which already features various $\pesp$ models, heuristics, and preprocessing techniques. To our time-expanded models, we apply the following two preprocessing steps:
\begin{itemize}
    \item The inherent symmetry of periodic timetables allows to fix a single event $v \in V(N)$ to time $\pi_v = 0$ (see, e.g., \citealt{Liebchen2006PeriodicTO}). We choose an event $v$ of maximum degree in $N$, create only the expanded node $v[0]$, and delete $v[t]$ for all $t \in \{1, \dots, T-1\}$.
    \item Before building the expanded network, we contract all events of degree $2$ in $N$. For baseline PESP, this leads to an objective function which is only piecewise linear (see, e.g., \citealt{goerigk_improved_2017}). However, the time expansion provides a linearization, cf.\ Remark~\ref{rem:non-linear}, so that the contraction is exact.
\end{itemize}

Due to the size of the time-expanded models (cf.\ Table~\ref{tab:instances}), we resort to column generation. As we want to have a full and flexible control of the column generation process, we choose SCIP 8.0.3 \citep{BestuzhevaEtal2021OO} as branch-cut-and-price engine, with Gurobi 10 \citep{gurobi} as underlying LP solver. The pricing problem for cycle variables is delegated to a custom implementation that uses topological search on directed acyclic graphs as described in Section~\ref{sss:pricing_cycle}, while for the transfer arc variables, we resort to SCIP's built-in variable pricer. We enhance the column generation process by a dual smoothing technique using stabilization centers along the lines of \cite{pessoa_exact_2010}.

Concerning passenger routing within the $\cxttp$ model, we restrict ourselves to the two smallest instance sets \emph{2linecross} and \emph{3berlin}. For each set, we generate a dense random OD-matrix, resulting in 20 and 72 OD-pairs, respectively. We again fix a single event $v$ to time $\pi_v = 0$. For the passenger paths, for each OD-pair $(s, t)$, we first enumerate all paths on the original event-activity network $N$ that can potentially be a shortest path (\citealp{karasan_robust_2012}, Proposition~2.3; \citealp{MASING2025100163}, Theorem~5). We then expand these paths to $D$, and use the union of all these paths as a routing graph for passengers from $s$ to $t$, implicitly defining the set $P_{st}$. The pricing routine for the passenger path variables $y_p$ (cf.\ Section~\ref{sec:pricing_passenger_flow}) is implemented using Dijkstra's algorithm on the corresponding routing graph. In each pricing round, we first price cycles as for $\cxpesp$, and then paths for each OD-pair.

All experiments are run on an Intel Xeon E3-1270 v6 CPU running at 3.8 GHz with 32 GB RAM.

\subsection{Results}
\label{ss:results}

\paragraph{Quality of the LP relaxation.}
We first evaluate the integrality gap of cXPESP by comparing the LP relaxation of $\cxpesp$ to the optimal objective value of the integer program, the latter being the same as for the standard PESP model. Instead of the weighted periodic tensions $\omega^T \mathbf{x}$, we measure the gap in terms of weighted periodic slacks $\omega^T (\mathbf{x} - l)$ by subtracting the sum of weighted lower bounds of all activities in $N$, so that the LP relaxations of both the standard PESP model and the XPESP model have an optimal objective value of $0$ (cf.\ Lemma~\ref{lem:xpesplowerbounds}).
Table~\ref{tab:obj_comparison} shows that the cXPESP formulation is indeed able to close a significant amount of the integrality gap. The closed gap is stable across different periods, but decreases for larger instances.

\begin{table}
\centering
     \caption{Optimal objective values in terms of weighted slack for the LP relaxation of the cXPESP model, the IP, and the closed gap. The instances \emph{berlin} (for $T = 10$) and \emph{R1L1-10} could not be solved to optimality within 24 hours, and the best dual and primal bounds are given as an interval.}
     \label{tab:obj_comparison}
     \begin{tabular}{lrrrrr}
    Instance & $T$ & Weighted slack (cxPESP-LP) & Weighted slack (IP) & Closed gap [\%]\\ 
    \hline
    2linecross & 5 & 4.00 & 4 & 100.00 \\
    2linecross & 10 & 17.00 & 24 & 70.83 \\
    2linecross & 15 & 34.13 & 44 & 77.57 \\
    2linecross & 20 & 51.33 & 64 & 80.20 \\
    2linecross & 30 & 86.18 & 104 & 82.87 \\
    2linecross & 40 & 121.13 & 144 & 84.12 \\
    2linecross & 50 & 156.10 & 184 & 84.84 \\
    2linecross & 60 & 191.08 & 224 & 85.30 \\ \hline
    3berlin & 5 & 31.67 & 45 & 70.38 \\
    3berlin & 10 & 77.00 & 106 & 72.64 \\
    3berlin & 15 & 135.24 & 144 & 93.92 \\
    3berlin & 20 & 182.40 & 300 & 60.80 \\
    3berlin & 30 & 309.70 & 392 & 79.01 \\
    3berlin & 40 & 389.33 & 568 & 68.54 \\
    3berlin & 50 & 534.50 & 902 & 59.26 \\
    3berlin & 60 & 684.78 & 1064 & 64.36 \\ 
    \hline
    berlin & 5 & 68.80 & 459 & 14.99 \\
    berlin & 10 & 156.40 & [1262, 1369] & [11.42, 12.39] \\ 
    \hline
    R1L1-1 & 60 & 0.00 & 0 & -- \\
    R1L1-2 & 60 & 120057.00 & 159586 & 75.23 \\
    R1L1-5 & 60 & 627475.00 & 1601566 & 39.18 \\
    R1L1-10 & 60 & 1707740.90 & [6460165, 6963841] & [24.52, 26.43] \\
    \end{tabular}
\end{table}

\paragraph{Computation times and pricing statistics.}
Having obtained a convincing quality of the LP relaxations, we now turn to a quantitative evaluation of computation times and the column generation process, comparing the pure cXPESP model, the slimmer $\cxpespw$ reformulation, and $\cxpespw$ with additional stabilization. Table~\ref{tab:pricing} shows for each instance the number of pricing rounds, the total number of added cycle variables $x_c$, the total time spent in pricing, and the total computation time for solving the LP at the root node. We have experimented with several stabilization factors $\zeta \in (0, 1]$, using convex combinations of the current dual solution with a weight of $\zeta$ and a stability center with a weight of $1-\zeta$, so that $\zeta = 1$ corresponds to no stabilization. In Table~\ref{tab:pricing}, we indicate the best performing value for $\zeta$ in terms of total computation time.
Indeed, the $\cxpespw$ model requires much less variables, which has a very positive impact on computation times in comparison to $\cxpesp$, while the pricing effort is comparable. The number of priced variables first of all does not explode, and can be significantly reduced with stabilization techniques. For the larger instances, pricing is not the bottleneck, but solving the LP is, and stabilization is indispensable to obtain a solution within a reasonable amount of time at all. However, a good stabilization factor seems hard to predict.

\begin{landscape}
\centering
     \begin{table}
     \caption{Pricing statistics. Depicted are the number of pricing rounds, number of priced cycle variables, total pricing time and total computation time for the LP relaxations of $\cxpesp$, $\cxpespw$ and $\cxpespw$ with best performing stabilization factor $\zeta$.}
     \label{tab:pricing}
     \begin{tabular}{lr|rrrr|rrrr|rrrrr}
     & & \multicolumn{4}{c}{$\cxpesp$-LP} & \multicolumn{4}{c}{$\cxpespw$-LP} & \multicolumn{5}{c}{$\cxpespw$-LP with stabilization} \\ 
     Instance & $T$ &
     Rounds & Cycles & Pricing [s] & Total [s] &
     Rounds & Cycles & Pricing [s] & Total [s] &
     Rounds & Cycles & Pricing [s] & Total [s] & Factor \\ 
     \hline
     2linecross & 5 & 6 & 27 & 0.014 & 0.113 & 5 & 32 & 0.015 & 0.063 & 5 & 32 & 0.015 & 0.063 & 1.0 \\
    2linecross & 10 & 16 & 63 & 0.004 & 0.046 & 10 & 52 & 0.028 & 0.080 & 12 & 39 & 0.025 & 0.073 & 0.5\\
    2linecross & 15 & 20 & 112 & 0.027 & 0.079 & 17 & 152 & 0.050 & 0.123 & 17 & 152 & 0.050 & 0.123 & 1.0  \\
    2linecross & 20 & 22 & 255 & 0.054 & 0.118 & 21 & 265 & 0.090 & 0.173 & 23 & 247 & 0.097 & 0.171 & 0.7  \\
    2linecross & 30 & 38 & 592 & 0.270 & 0.473 & 37 & 603 & 0.317 & 0.440 & 32 & 412 & 0.288 & 0.387 & 0.9 \\
    2linecross & 40 & 49 & 460 & 0.710 & 1.093 & 47 & 860 & 0.742 & 0.885 & 47 & 473 & 0.725 & 0.845 & 0.4 \\
    2linecross & 50 & 57 & 697 & 1.536 & 2.238 & 59 & 1235 & 1.608 & 1.875 & 55 & 784 & 1.473 & 1.685 & 0.9 \\
    2linecross & 60 & 74 & 953 & 3.330 & 4.550 & 67 & 1382 & 2.927 & 3.167 & 62 & 1285 & 2.774 & 3.016 & 0.8 \\ 
    \hline
    3berlin & 5 & 31 & 200 & 0.020 & 0.107 & 34 & 242 & 0.078 & 0.228 & 24 & 195 & 0.061 & 0.176 & 0.9 \\
    3berlin & 10 & 44 & 547 & 0.088 & 0.474 & 48 & 577 & 0.158 & 0.475 & 37 & 463 & 0.145 & 0.386 & 0.6 \\
    3berlin & 15 & 51 & 702 & 0.234 & 1.263 & 46 & 729 & 0.281 & 0.652 & 40 & 696 & 0.262 & 0.600 & 0.9 \\
    3berlin & 20 & 76 & 1564 & 0.661 & 4.088 & 108 & 1747 & 1.045 & 2.590 & 65 & 1050 & 0.675 & 1.595 & 0.3 \\
    3berlin & 30 & 75 & 2682 & 1.705 & 11.406 & 67 & 2311 & 1.586 & 3.426 & 64 & 1767 & 1.553 & 2.956 & 0.3 \\
    3berlin & 40 & 95 & 3736 & 4.214 & 30.458 & 82 & 3541 & 3.549 & 7.313 & 71 & 3143 & 3.210 & 6.441 & 0.8 \\
    3berlin & 50 & 118 & 7034 & 9.288 & 65.237 & 132 & 7813 & 9.900 & 26.213 & 87 & 5877 & 6.784 & 15.281 & 0.3\\
    3berlin & 60 & 132 & 8131 & 16.598 & 122.47 & 135 & 8015 & 16.523 & 32.015 & 105 &  7210 & 13.054 & 24.421 & 0.8 \\ 
    \hline
    berlin & 5 & 448 & 7430 & 3.231 & 212.201 & 413 & 7044 & 3.004 & 123.391 & 239 & 1360 & 1.835 & 19.232 & 0.01\\
    berlin & 10 & 192 & 8575 & 4.424 & 2200.058 & 355 & 10924 & 7.947 & 1122.544 & 242 & 3962 & 5.999 & 253.817 & 0.05 \\ 
    \hline
    R1L1-1 & 60 & 0 & 0 & 0.000 & 0.075 & 0 & 0 & 0.000 & 0.055 & 0 & 0 & 0.000 & 0.055 & 1.0 \\
    R1L1-2 & 60 & 25 & 520 & 3.130 & 3.805 & 55 & 1480 & 7.544 & 9.572 & 55 & 1480 & 7.544 & 9.572 & 1.0 \\
    R1L1-5 & 60 & 125 & 27253 & 85.557 & 2189.128 & 143 & 24821 & 98.305 & 205.493 & 120 & 20614 & 84.523 & 144.486 & 0.2 \\
    R1L1-10 & 60 & -- & -- & -- & $\geq$ 24h & -- & -- & -- & $\geq$ 24h & 321 & 117934 & 696.124 & 17000.746 & 0.5 \\
    \end{tabular}
    \end{table}
\end{landscape}

\paragraph{Branch-cut-price.}
Beyond LP relaxations, we attempt to solve $\cxpespw$ as an integer program using column generation. We compare $\cxpespw$ (with stabilization) to the standard incidence-based integer programming formulation of PESP as presented in Section~\ref{section:PESP}, and to the common cycle-based formulation using a fundamental cycle basis from a minimum spanning tree (see, e.g., \citealt{nachtigall_periodic_1998} and \citealt{Liebchen2006PeriodicTO}). The PESP models are solved with SCIP using default settings, still with Gurobi as LP solver. Table~\ref{tab:bb} collects the total running time and the number of nodes of the branch-and-bound tree for $\cxpespw$ and the two compact PESP formulations. For $\cxpespw$, we also list the total number of pricing rounds, the total number of added cycle variables, and the total pricing time, summed over all branch-and-bound nodes.
The upshot is that $\cxpespw$ is impractical to solve PESP instances to optimality, although the method works in principle. The pricing time becomes much more dominant in comparison to the LP time in a branch-and-bound context. Although there is an advantage in terms of the number of required nodes for quite some instances, the overall process is too slow to be competitive. Moreover, it is striking that the cycle-based formulation for PESP performs much better than the incidence-based formulation, e.g., by one order of magnitude in terms of computation time and nodes for the 3berlin instances.

\begin{table}
\centering
    \caption{Branch-cut-price statistics. Depicted are the total number of pricing rounds, the number of priced cycle variables, and the pricing time for solving $\cxpespw$ as an IP. Furthermore, the total computation time and the number of branch-and-bound nodes for $\cxpespw$-IP, for a cycle-based, and for the incidence-based IP formulation for PESP.}
    \label{tab:bb}
    \begin{tabular}{lr|rrrrr|rr|rr}
     & & \multicolumn{5}{c|}{$\cxpespw$-IP} & \multicolumn{2}{c|}{PESP-IP (cycle)} & \multicolumn{2}{c}{PESP-IP (incidence)}  \\ 
     Instance & $T$ & Rounds & Cycles & Pricing [s] & Total [s] & Nodes & Total [s] & Nodes & Total[s] & Nodes \\ 
     \hline
    2linecross & 5 & 5 & 32 & 0.015 & 0.095 & 1 & 0.115 & 1 & 0.11 & 1\\
    2linecross & 10 & 60 & 109 & 0.064 & 0.293 & 19 & 0.102 & 14 & 0.18 & 179 \\
    2linecross & 15 & 191 & 510 & 0.362 & 0.991 & 70 & 0.103 & 1 & 0.12 & 34\\
    2linecross & 20 & 315 & 918 & 0.994 & 2.083 & 119 & 0.096 & 1 & 0.14 & 160\\
    2linecross & 30 & 678 & 2096 & 5.141 & 7.863 & 219 & 0.090 & 1 & 0.13 & 100\\
    2linecross & 40 & 886 & 3971 & 13.135 & 17.544 & 319 & 0.105 & 1 & 0.15 & 115\\
    2linecross & 50 & 1582 & 6333 & 41.400 & 51.265 & 419 & 0.096 &  1 & 0.16 & 57\\
    2linecross & 60 & 1683 & 8461 & 73.452 & 85.963 & 519 & 0.097 & 1 & 0.16 & 188\\ 
    \hline
    3berlin & 5 & 151 & 570 & 0.226 & 1.001 & 18 & 0.601 & 287 & 1.26 & 1051 \\
    3berlin & 10 & 297 & 1525 & 0.845 & 3.084 & 18 & 0.301 & 119 & 3.77 & 1492\\
    3berlin & 15 & 263 & 1696 & 1.471 & 4.169 & 13 & 0.301 & 22 & 2.15 & 1339\\
    3berlin & 20 & 881 & 5646 & 8.272 & 24.267 & 95 & 0.599 & 295 & 7.58 & 5056\\
    3berlin & 30 & 569 & 6442 & 13.025 & 36.941 & 35 & 0.405 & 77 & 2.13 & 1127\\
    3berlin & 40 & 894 & 11520 & 37.415 & 93.951 & 62 & 0.604 & 181 & 4.08 & 1398 \\
    3berlin & 50 & 4012 & 37089 & 271.431 & 811.713 & 322 & 0.706 & 352 & 9.93 & 5154\\
    3berlin & 60 & 2007 & 31472 & 234.081 & 511.205 & 126 & 0.489 & 273 & 8.86 & 3242\\ 
    \hline
    berlin & 5 & -- & -- & -- & -- & -- & -- & -- & -- & --  \\
    berlin & 10 & -- & -- & -- & -- & -- & -- & -- & -- & --\\ 
    \hline
    R1L1-1 & 60 & 0 & 0 & 0.000 & 0.055 & 1 & 0.125 & 1 & 0.01 & 1\\
    R1L1-2 & 60 & 145 & 1770 & 18.423 & 24.046 & 33 & 0.142 & 7 & 3.38 & 3265\\
    R1L1-5 & 60 & -- & -- & -- & -- & -- & 8.026 & 1792 & -- & --\\
    R1L1-10 & 60 & -- & -- & -- & -- & -- & -- & -- & -- & --\\
    \end{tabular}
\end{table}

\paragraph{Integrating passenger routing.} The quality of the LP relaxation given by $\cxttp$ is convincing: For \emph{2linecross} and all considered period times, the gap is always closed at the root node, as is for \emph{3berlin} and $T = 40$. Otherwise, higher $T$ seem to imply a tighter gap, see Table~\ref{tab:obj_comparison_cxttp}. Examinating the pricing statistics in Table~\ref{tab:pricing_cxttp}, we note that $\cxttp$ requires more cycles than $\cxpesp$, the number of cycles is roughly comparable to the number of paths, and that the pricing procedure for paths is faster. What is however striking is the large amount of time required to solve the arising linear programs: For example, \emph{3berlin-40} spends almost 12 hours in LP, while all pricing steps together take in total less than 90 seconds. The huge computational demand for LP solving makes it practically impossible to solve larger instances, and this is why we restrict to only two instances, and omit a detailed analysis of branch-and-cut experiments for the IP. While the number of nodes is smaller -- the IP is solved at the root node for \emph{2linecross} and \emph{3berlin-40} -- computation times explode even further, while SCIP with Gurobi as LP solver always manages to solve the bilinear $\ttp$ integer program for \emph{2linecross} and \emph{3berlin} in less than 60 seconds per instance.

\begin{table}
\centering
     \caption{Optimal objective values for the trivial QP relaxation of the TTP model, for the LP relaxation of the cXTTP model, the IP, and the closed gap. The closed gap is computed as $(\text{IP obj.} - \text{TTP-QP obj.})/(\text{cXTTP-LP obj.} - \text{TTP-QP obj.})$. 
     The cXTTP root LP computation did not terminate within 24 hours for \emph{3berlin-50} and stopped early after approximately 19 hours due to numerical troubles for \emph{3berlin-60}.}
     \label{tab:obj_comparison_cxttp}
     \begin{tabular}{lrrrrrr}
    Instance & $T$ & Objective ($\ttp$-QP) & Objective ($\cxttp$-LP) & Objective (IP) & Closed gap [\%]\\ 
    \hline
    2linecross &  5 & \multirow{8}{*}{132} & 143 & 143 & 100.00\\
    2linecross & 10 & & 162 & 162 & 100.00\\
    2linecross & 15 & & 182 & 182 & 100.00\\
    2linecross & 20 & & 202 & 202 & 100.00\\
    2linecross & 30 & & 242 & 242 & 100.00\\
    2linecross & 40 & & 282 & 282 & 100.00\\
    2linecross & 50 & & 322 & 322 & 100.00\\
    2linecross & 60 & & 362 & 362 & 100.00\\ \hline
    3berlin &  5 & \multirow{8}{*}{5632} & 5694.60 & 5765 & 47.07\\
    3berlin & 10 & & 5764.80 & 5820 & 70.64 \\
    3berlin & 15 & & 5863.93 & 5915 & 81.96 \\
    3berlin & 20 & & 5945.33 & 5990 & 87.52 \\
    3berlin & 30 & & 6036.00 & 6040 & 99.02 \\
    3berlin & 40 & & 6050.00 & 6050 & 100.00 \\
    3berlin & 50 & & -- & 6080 & --\\
    3berlin & 60 & & -- & 6090 & --\\
    \end{tabular}
\end{table}

 \begin{table}
 \centering
     \caption{Pricing statistics. The table shows the number of pricing rounds, priced cycle and passenger path variables, total pricing time for cycles and paths, and total computation time for the LP relaxation of $\cxttp$ using column generation.}
     \label{tab:pricing_cxttp}
     \begin{tabular}{lr|rrrrrr}
     & & \multicolumn{6}{c}{$\cxttp$-LP} \\ 
     Instance & $T$ & Rounds & Cycles & Paths & Cycle pricing [s] & Path pricing [s] & Total [s]  \\ 
     \hline
     2linecross & 5 & 20 & 80 & 142 & 0.017 & 0.007 & 0.074 \\
    2linecross & 10 & 34 & 232 & 277 & 0.090 & 0.034 & 0.390 \\
    2linecross & 15 & 48 & 403 & 392 & 0.246 & 0.086 & 1.110 \\
    2linecross & 20 & 54 & 317 & 406 & 0.459 & 0.139 & 1.651 \\
    2linecross & 30 & 64 & 724 & 636 & 1.272 & 0.319 & 4.925 \\
    2linecross & 40 & 107 & 1822 & 1060 & 4.026 & 0.849 & 18.226 \\
    2linecross & 50 & 126 & 2523 & 1279 & 7.525 & 1.365 & 31.006 \\
    2linecross & 60 & 409 & 3400 & 1798 & 36.881 & 6.431 & 107.381 \\
    \hline
     3berlin & 5 & 42 & 387 & 846 & 0.235 & 0.108 & 3.660 \\
    3berlin & 10 & 75 & 1305 & 1975 & 1.153 & 0.550 & 60.119  \\
    3berlin & 15 & 120 & 3021 & 3055 & 3.466 & 2.289 & 388.485  \\
    3berlin & 20 & 137 & 4320 & 4316 & 6.23 & 3.968 & 1750.609  \\
    3berlin & 30 & 204 & 8106 & 6880 & 18.976 & 11.189 & 12310.673  \\
    3berlin & 40 & 360 & 13742 & 9337 & 56.958 & 32.002 & 41456.027  \\
    3berlin & 50 & -- & -- & -- & -- & -- & $\geq$ 24h  \\
    3berlin & 60 & -- & -- & -- & -- & -- & -- 
    \end{tabular}
    \end{table}

\section{Conclusions}
\label{section:conclusions}

We presented a new model for periodic timetabling based on a graph expansion of the event-activity network used in PESP. For this new model, we introduced a novel path-based ($\cxpesp$) integer programming formulation. We demonstrated that the solution space of the integer program of $\cxpesp$ is identical to the corresponding solution space for PESP, while providing a tighter linear programming relaxation. 
The resulting lower bound on the linear programming relaxation to $\cxpesp$ is, to our knowledge, the best known to date, which is supported by the validity of cycle, change cycle, and flip-cycle inequalities on the line cycles of the underlying network.
This effect is a result of including more of the problem's inherent structure into the programming formulation for the operated lines.
The enhanced structure comes with an increased number of variables. We handled the increased size with the use of column generation and, therefore, introduced the pricing problems for different variable types in $\cxpesp$ and gave a suggestion on how to solve them. The pricing of cycle variables results in a set of shortest path problems on an acyclic graph. We further described an alternative linearization of the transfer arcs to deal with their number, which comes with a slightly weaker LP relaxation.

Computational experiments confirmed that solving the linear programming of $\cxpesp$ closes a large part of the integrality gap. The column generation procedure effectively controls the number of generated variables. The linearization of the transfer arcs further reduces the number of variables. However, even for small instances, the bottleneck lies in solving the linear program rather than the pricing. While solving the full model to integrality reduces the number of nodes in the branch-and-bound tree, it remains computationally intensive.

Finally, we extended the path-based timetabling model to integrate passenger routing ($\cxttp$), which inherits the advantages of $\cxpesp$. Again, we tackled the increased problem size by the use of column generation. The pricing problem for the passenger flow variables results in a shortest path problem on the expanded event-activity network. While the pricing itself scales reasonably, solving the linear programs is even more tedious than in the $\cxpesp$ case, prohibiting successful computations for meaningful $\cxttp$ instances. However, the theoretical strengths are worth mentioning: Once the LP relaxation has been computed, the remaining integrality gap is small.

In summary, the proposed path-based timetabling model and its passenger-flow-integrated variant demonstrate theoretical advantages and also improvements in closing the integrality gap. While column generation effectively mitigates the growth in model size, solving the integer programming formulation remains computationally challenging and not competitive in practice. The contribution is therefore mostly on theory. A potential direction for future is to further exploit the inherent symmetries of the problem to get a better control of the number of generated columns and hence to accelerate solving times.

\section*{Acknowledgement}
The work for this article has been conducted in the Research Campus MODAL funded by the Federal Ministry of Research, Technology and Space (BMFTR) (fund numbers 05M14ZAM, 05M20ZBM, 05M2025).

\bibliographystyle{abbrvnat}
\bibliography{references} 

\end{document}

%% file: graphics/figure1.tex
    \begin{tikzpicture}[node distance={20mm}, minimum size=2mm, main/.style = {draw, circle}, time/.style = {}]
	     
     \node[main] (1c) at (-9,1) { };
    \node[main] (2c) [left of=1c] { };
    \node[main] (3c) [right of=1c] { };
    \node[main] (4c) [above of=1c] {};
    \node[main] (5c) [below of=1c] { };
    
    \draw[blue] (2c) to node[midway, above]  {}  (1c) to node[midway, above]  {} (3c);
    \draw[red] (4c) to node[midway, right]  {}  (1c) to node[midway, right]  {} (5c);

    \end{tikzpicture}

    

%% file: graphics/figure2.tex
    \begin{tikzpicture}[node distance={20mm}, minimum size=2mm, main/.style = {draw, circle}, time/.style = {}]
    
    \node[main] (1) at (-5,4) {};
    \node[main ](2) [right of=1] {};
    \node[main ](3) [right of=2]  {};
    \node[main ](4) [right of=3] {};
    \node[main] (5) [below of=4]{};
    \node[main ](6) [below of=3] {};
    \node[main ](7) [below of=2]{};
    \node[main ](8) [below of=1]{};
    
   \node[main] (1b) [below of=8] {};
    \node[main ](2b) [right of=1b] {};
    \node[main ](3b) [right of=2b]  {};
    \node[main ](4b) [right of=3b] {};
    \node[main] (5b) [below of=4b]{};
    \node[main ](6b) [below of=3b] {};
    \node[main ](7b) [below of=2b]{};
    \node[main ](8b) [below of=1b]{};

    \draw[blue, ->] (1) to node[midway, above]  {$[1,1]$}  (2);
    \draw[blue, ->] (2) to node[midway, above]  {$[0,1]$} (3);
    \draw[blue, ->] (3) to node[midway, above]  {$[1,1]$} (4);
    \draw[blue, ->] (4) to node[midway, right]  {$[2,3]$} (5); 
    \draw[blue, ->] (5) to node[midway, above]  {$[1,1]$} (6);
    \draw[blue, ->] (6) to node[midway, above]  {$[0,1]$} (7);
    \draw[blue, ->] (7) to node[midway, above]  {$[1,1]$} (8);
    \draw[blue, ->] (8) to node[midway, left]  {$[2,3]$} (1);
    
    \draw[red, ->] (1b) to node[midway, below]  {$[1,1]$}  (2b);
    \draw[red, ->] (2b) to node[midway, below]  {$[0,1]$} (3b);
    \draw[red, ->] (3b) to node[midway, below]  {$[1,1]$} (4b);
    \draw[red, ->] (4b) to node[midway, right]  {$[2,3]$} (5b); 
    \draw[red, ->] (5b) to node[midway, below]  {$[1,1]$} (6b);
    \draw[red, ->] (6b) to node[midway, below]  {$[0,1]$} (7b);
    \draw[red, ->] (7b) to node[midway, below]  {$[1,1]$} (8b);
    \draw[red, ->] (8b) to node[midway, left]  {$[2,3]$} (1b);
    
    \draw[gray, ->] (2b) to node[midway, left] {$[0,1]$} (7);
    \draw[gray, ->] (6) to node[midway, right]  {$[0,1]$}  (3b);
    \draw[gray, ->] (6b) to node[midway, left]  {$[0,1]$} (7);
    \draw[gray, ->] (2) to node[midway, left]  {$[0,1]$}  (3b);
    \draw[gray, ->] (2b) to node[midway, left]  {$[0,1]$}  (3);
    \draw[gray, ->] (6) to node[midway, left]  {$[0,1]$}  (7b);
    \draw[gray, ->] (2) to [out=220, in=140, looseness=0.5] node[midway, left]  {$[0,1]$} (7b);
    \draw[gray, ->] (6b) to [out=40, in=320, looseness=0.5] node[midway, right]  {$[0,1]$}  (3);
    
     \end{tikzpicture} 

%% file: graphics/figure3.tex
   \begin{tikzpicture}[node distance={7mm}, inner sep=1pt, font=\small, minimum size=5mm, main/.style = {draw, circle}, time/.style = {}]
    
    \node[main] (11d1p) at (4,5.5) {$1$};
    \node[main] (11d0p) [below of=11d1p] {$0$};    
    \node[main] (12a1p) at (7,5.5) {$1$};
    \node[main] (12a0p) [below of=12a1p] {$0$};    
    \node[main] (12d1p) at (9,5.5) {$1$};
    \node[main] (12d0p) [below of=12d1p] {$0$};    
    \node[main] (13d1p) at (12,5.5) {$1$};
    \node[main] (13d0p) [below of=13d1p] {$0$};   
    \node[main] (13d1m) at (12,3) {$1$};
    \node[main] (13d0m) [below of=13d1m]  {$0$};        
    \node[main] (12a1m) at (9,3) {$1$};
    \node[main] (12a0m) [below of=12a1m] {$0$};    
    \node[main] (12d1m) at (7,3) {$1$};
    \node[main] (12d0m) [below of=12d1m] {$0$};    
    \node[main] (11a1m) at (4,3) {$1$};
    \node[main] (11a0m) [below of=11a1m]{$0$};
       
    \node[main] (24d1p) at (4,0) {$1$};
    \node[main] (24d0p) [below of=24d1p] {$0$};    
    \node[main] (22a1p) at (7,0) {$1$};
    \node[main] (22a0p) [below of=22a1p] {$0$};    
    \node[main] (22d1p) at (9,0) {$1$};
    \node[main] (22d0p) [below of=22d1p] {$0$};
    \node[main] (25a1p) at (12,0) {$1$};
    \node[main] (25a0p) [below of=25a1p]{$0$};    
    \node[main] (25d1m) at (12,-2.5) {$1$};
    \node[main] (25d0m) [below of=25d1m] {$0$};    
    \node[main] (22a1m) at (9,-2.5) {$1$};
    \node[main] (22a0m) [below of=22a1m] {$0$};    
    \node[main] (22d1m) at (7,-2.5) {$1$};
    \node[main] (22d0m) [below of=22d1m] {$0$};    
    \node[main] (24a1m) at (4,-2.5) {$1$};
    \node[main] (24a0m) [below of=24a1m] {$0$};

    \draw[blue, ->] (11d1p) --  (12a0p);
    \draw[blue, ->] (11d0p) --  (12a1p);
    \draw[blue, ->] (12a0p) --  (12d0p);
    \draw[blue, ->] (12a1p) --  (12d0p);
    \draw[blue, ->] (12a0p) --  (12d1p);
    \draw[blue, ->] (12a1p) --  (12d1p);
    \draw[blue, ->] (12d0p) --  (13d1p);
    \draw[blue, ->] (12d1p) --  (13d0p);
    \draw[blue, ->] (13d0p) to [out=330, in=30, looseness=1] (13d1m);
    \draw[blue, ->] (13d1p) to [out=330, in=30, looseness=1] (13d0m);
    \draw[blue, ->] (13d1p) to [out=330, in=30, looseness=1] (13d1m);
    \draw[blue, ->] (13d0p) to [out=330, in=30, looseness=1] (13d0m);
    \draw[blue, ->] (13d0m) --  (12a1m);
    \draw[blue, ->] (13d1m) --  (12a0m);
    \draw[blue, ->] (12a0m) --  (12d0m);
    \draw[blue, ->] (12a1m) --  (12d0m);
    \draw[blue, ->] (12a0m) --  (12d1m);
    \draw[blue, ->] (12a1m) --  (12d1m);
    \draw[blue, ->] (12d0m) --  (11a1m);
    \draw[blue, ->] (12d1m) --  (11a0m);
    \draw[blue, ->] (11a1m) to [out=150, in=210, looseness=1] (11d0p);
    \draw[blue, ->] (11a0m) to [out=150, in=210, looseness=1] (11d1p);
    \draw[blue, ->] (11a0m) to [out=150, in=210, looseness=1] (11d0p);
    \draw[blue, ->] (11a1m) to [out=150, in=210, looseness=1] (11d1p);
    
    \draw[red, ->] (24d0p) --  (22a1p);
    \draw[red, ->] (24d1p) --  (22a0p);
    \draw[red, ->] (22a1p) --  (22d1p);
    \draw[red, ->] (22a0p) --  (22d1p);
    \draw[red, ->] (22a1p) --  (22d0p);
    \draw[red, ->] (22a0p) --  (22d0p);
    \draw[red, ->] (22d0p) --  (25a1p);
    \draw[red, ->] (22d1p) --  (25a0p);
    \draw[red, ->] (25a0p) to [out=330, in=30, looseness=1] (25d0m);
    \draw[red, ->] (25a0p) to [out=330, in=30, looseness=1] (25d1m);
    \draw[red, ->] (25a1p) to [out=330, in=30, looseness=1] (25d0m);
    \draw[red, ->] (25a1p) to [out=330, in=30, looseness=1] (25d1m);
    \draw[red, ->] (25d0m) --  (22a1m);
    \draw[red, ->] (25d1m) --  (22a0m);
    \draw[red, ->] (22a1m) --  (22d1m);
    \draw[red, ->] (22a0m) --  (22d1m);
    \draw[red, ->] (22a1m) --  (22d0m);
    \draw[red, ->] (22a0m) --  (22d0m);
    \draw[red, ->] (22d0m) --  (24a1m);
    \draw[red, ->] (22d1m) --  (24a0m);
    \draw[red, ->] (24a0m) to [out=150, in=210, looseness=1] (24d0p);
    \draw[red, ->] (24a1m) to [out=150, in=210, looseness=1] (24d0p);
    \draw[red, ->] (24a1m) to [out=150, in=210, looseness=1] (24d1p);
    \draw[red, ->] (24a0m) to [out=150, in=210, looseness=1] (24d1p);
    
    \draw[gray, ->] (12a0p) -- (22d1p);
    \draw[gray, ->] (12a0p) -- (22d0p);
    \draw[gray, ->] (12a0p) to [out=220, in=140, looseness=0.5] (22d0m);
    \draw[gray, ->] (12a0p) to [out=220, in=140, looseness=0.5] (22d1m);
    \draw[gray, ->] (12a1p) -- (22d1p);
    \draw[gray, ->] (12a1p) -- (22d0p);
    \draw[gray, ->] (12a1p) to [out=220, in=140, looseness=0.5] (22d0m);
    \draw[gray, ->] (12a1p) to [out=220, in=140, looseness=0.5] (22d1m);
    
    \draw[gray, ->] (12a0m) -- (22d1m);
    \draw[gray, ->] (12a0m) -- (22d0m);
    \draw[gray, ->] (12a0m) to [out=320, in=40, looseness=0.5] (22d1p);
    \draw[gray, ->] (12a0m) to [out=320, in=40, looseness=0.5] (22d0p);
    \draw[gray, ->] (12a1m) -- (22d1m);
    \draw[gray, ->] (12a1m) -- (22d0m);
    \draw[gray, ->] (12a1m) to [out=320, in=40, looseness=0.5] (22d1p);
    \draw[gray, ->] (12a1m) to [out=320, in=40, looseness=0.5] (22d0p);
    
    \draw[gray, ->] (22a0p) -- (12d0p);
    \draw[gray, ->] (22a0p) -- (12d1p);
    \draw[gray, ->] (22a0p) to [out=140, in=220, looseness=0.5] (12d0m);
    \draw[gray, ->] (22a0p) to [out=140, in=220, looseness=0.5] (12d1m);
    \draw[gray, ->] (22a1p) -- (12d0p);
    \draw[gray, ->] (22a1p) -- (12d1p);
    \draw[gray, ->] (22a1p) to [out=140, in=220, looseness=0.5] (12d0m);
    \draw[gray, ->] (22a1p) to [out=140, in=220, looseness=0.5] (12d1m);
    
    \draw[gray, ->] (22a0m) -- (12d0m);
    \draw[gray, ->] (22a0m) -- (12d1m);
    \draw[gray, ->] (22a0m) to [out=40, in=320, looseness=0.5] (12d0p);
    \draw[gray, ->] (22a0m) to [out=40, in=320, looseness=0.5] (12d1p);
    \draw[gray, ->] (22a1m) -- (12d0m);
    \draw[gray, ->] (22a1m) -- (12d1m);
    \draw[gray, ->] (22a1m) to [out=40, in=320, looseness=0.5] (12d0p);
    \draw[gray, ->] (22a1m) to [out=40, in=320, looseness=0.5] (12d1p);    
    
    \end{tikzpicture}  

%% file: graphics/figure4.tex
  \begin{tikzpicture}[node distance={5mm}, inner sep=1pt, minimum size=3mm, font=\scriptsize, main/.style = {draw, circle}, time/.style = {}, invisi/.style = {}]

    \node[main, blue, thick] (11d2p) at (4,5.5) {$2$};
    \node[main, red, thick] (11d1p) [below of=11d2p]  {$1$};
    \node[main, red, thick] (11d0p) [below of=11d1p] {$0$};
    
    \node[main, red, thick] (12a2p) at (7,5.5) {$2$};
    \node[main, red, thick] (12a1p) [below of=12a2p] {$1$};
    \node[main, blue, thick] (12a0p) [below of=12a1p] {$0$};
    
    \node[main, red, thick] (12d2p) at (9,5.5) {$2$};
    \node[main, blue, thick] (12d1p) [below of=12d2p] {$1$};
    \node[main, red, thick] (12d0p) [below of=12d1p] {$0$};
    
    \node[main, blue, thick] (13d2p) at (12,5.5) {$2$};
    \node[main, red, thick] (13d1p) [below of=13d2p] {$1$};
    \node[main, red, thick] (13d0p) [below of=13d1p] {$0$};
    
    \node[main, red, thick] (13d2m) at (12,3) {$2$};
    \node[main, red, thick] (13d1m) [below of=13d2m] {$1$};
    \node[main, blue, thick] (13d0m) [below of=13d1m]  {$0$};
        
    \node[main, red, thick] (12a2m) at (9,3) {$2$};
    \node[main, blue, thick] (12a1m) [below of=12a2m] {$1$};
    \node[main, red, thick] (12a0m) [below of=12a1m] {$0$};
    
    \node[main, red, thick] (12d2m) at (7,3) {$2$};
    \node[main, red, thick] (12d1m) [below of=12d2m] {$1$};
    \node[main, blue, thick] (12d0m) [below of=12d1m] {$0$};
    
    \node[main, red, thick] (11a2m) at (4,3) {$2$};
    \node[main, blue, thick] (11a1m) [below of=11a2m] {$1$};
    \node[main, red, thick] (11a0m) [below of=11a1m]{$0$};

    \draw[blue, ->, thick] (11d2p) --  (12a0p);
    \draw[red, thick, ->] (11d1p) --  (12a2p);
    \draw[red, thick, ->] (11d0p) --  (12a1p);
    
    \draw[->] (12a0p) --  (12d2p);
    \draw[->] (12a1p) --  (12d2p);
    \draw[->] (12a2p) --  (12d0p);
    \draw[->] (12a2p) --  (12d1p);
    \draw[->] (12a1p) --  (12d1p);
    \draw[->] (12a0p) --  (12d0p);
    \draw[blue, ->, thick] (12a0p) --  (12d1p);
    \draw[red, thick, ->] (12a1p) --  (12d0p);
    \draw[red, thick, ->] (12a2p) --  (12d2p);
    
    \draw[red, thick, ->] (12d2p) --  (13d0p);
    \draw[blue, ->, thick] (12d1p) --  (13d2p);
    \draw[red, thick, ->] (12d0p) --  (13d1p);
    
    \draw[->] (13d0p) to [out=330, in=30, looseness=1]  (13d2m);
    \draw[->] (13d1p) to [out=330, in=30, looseness=1]  (13d0m);
    \draw[->] (13d2p) to [out=330, in=30, looseness=1]  (13d1m);
    \draw[red, thick, ->] (13d1p) to [out=330, in=30, looseness=1] (13d2m);
    \draw[blue, ->, thick] (13d2p) to [out=330, in=30, looseness=1]  (13d0m);
    \draw[red, thick, ->] (13d0p) to [out=330, in=30, looseness=1]  (13d1m);

    \draw[red, thick, ->] (13d2m) --  (12a0m);
    \draw[red, thick, ->] (13d1m) --  (12a2m);
    \draw[blue, ->, thick] (13d0m) --  (12a1m);
    
    \draw[->] (12a0m) --  (12d1m);
    \draw[->] (12a1m) --  (12d2m);
    \draw[->] (12a2m) --  (12d0m);
    \draw[->] (12a2m) --  (12d2m);
    \draw[->] (12a1m) --  (12d1m);
    \draw[->] (12a0m) --  (12d0m);
    \draw[red, thick, ->] (12a0m) --  (12d2m);
    \draw[blue, ->, thick] (12a1m) --  (12d0m);
    \draw[red, thick, ->] (12a2m) --  (12d1m);
    
    \draw[red, thick, ->] (12d2m) --  (11a0m);
    \draw[red, thick, ->] (12d1m) --  (11a2m);
    \draw[blue, ->, thick] (12d0m) --  (11a1m);
    
    \draw[->] (11a2m) to [out=150, in=210, looseness=1] (11d1p);
    \draw[->] (11a0m) to [out=150, in=210, looseness=1] (11d2p);
    \draw[->] (11a1m) to [out=150, in=210, looseness=1] (11d0p);
    \draw[blue, ->, thick] (11a1m) to [out=150, in=210, looseness=1] (11d2p);
    \draw[red, thick, ->] (11a2m) to [out=150, in=210, looseness=1] (11d0p);
    \draw[red, thick, ->] (11a0m) to [out=150, in=210, looseness=1] (11d1p);
    
    \end{tikzpicture}  

%% file: graphics/figure5.tex
        \begin{tikzpicture}[node distance={15mm}, minimum size=5mm, main/.style = {draw, circle}, time/.style = {}, invisi/.style = {}]
    \node[main] (4) at (9,2) {};
    \node[main] (5) at (9,1) {}; 
    \node[main] (6) at (9,0) {};
    \node[main] (10) at (12,2) {};
    \node[main] (11) at (12,1) {};
    \node[main] (12) at (12,0) {};
    
     \node[invisi] (i1) at (7,1.2) {};
     \node[invisi] (i2) at (7,0.2) {};
     \node[invisi] (i3) at (7,-0.8) {};
     \node[invisi] (i4) at (6,2.6) {};
     \node[invisi] (i5) at (6,1.6) {};
     \node[invisi] (i6) at (6,0.6) {};
     
     \node[invisi] (j1) at (14,1.2) {};
     \node[invisi] (j2) at (14,0.2) {};
     \node[invisi] (j3) at (14,-0.8) {};
     \node[invisi] (j4) at (15,2.6) {};
     \node[invisi] (j5) at (15,1.6) {};
     \node[invisi] (j6) at (15,0.6) {};
     
     \node[invisi] (b4) at (9,2.8) {$v$};
     \node[invisi] (b5) at (12,2.8) {$w$};
     
     \node[invisi] (o1) at (10.5, 2) {};
     \node[invisi] (o2) at (10.5, 1.5) {};
     \node[invisi] (o3) at (10.5, 1) {};
     \node[invisi] (o4) at (10.5, 0.5) {};
     \node[invisi] (o5) at (10.5, 0) {};
    
    \draw[orange, -] (4) -- (o1);
    \draw[orange, -] (4) -- (o2);
    \draw[orange, -] (4) -- (o3);
    \draw[red, -] (5) -- (o2);
    \draw[red, -] (5) -- (o3);
    \draw[red, -] (5) -- (o4);
    \draw[violet, -] (6) -- (o3);
    \draw[violet, -] (6) -- (o4);
    \draw[violet, -] (6) -- (o5);
    \draw[blue, ->] (o1) -- (10);
    \draw[blue, ->] (o2) -- (10);
    \draw[blue, ->] (o3) -- (10);
    \draw[cyan, ->] (o2) -- (11);
    \draw[cyan, ->] (o3) -- (11);
    \draw[cyan, ->] (o4) -- (11);
    \draw[green, ->] (o3) -- (12);
    \draw[green, ->] (o4) -- (12);
    \draw[green, ->] (o5) -- (12);
    
    \draw[orange, ->] (i4) to [out=0, in=10, looseness=4.5] (i1);
    \draw[red, ->] (i5) to [out=0, in=10, looseness=4.5] (i2);
    \draw[violet, ->] (i6) to [out=0, in=10, looseness=4.5] (i3);
    
    \draw[blue, ->] (j4) to [out=180, in=170, looseness=4.5] (j1);
    \draw[cyan, ->] (j5) to [out=180, in=170, looseness=4.5] (j2);
    \draw[green, ->] (j6) to [out=180, in=170, looseness=4.5] (j3);

    \end{tikzpicture} 

%% file: graphics/figure7.tex
    \begin{tikzpicture}[node distance={15mm}, minimum size=5mm, main/.style = {draw, circle}, time/.style = {}]
    \node[main] (1) at (0,0) {$v_0$};
    \node[main ](2) at (3,0) {$v_1$};
    \node[main ](3) at (3,-2) {$v_2$};
    \node[main ](4) at (0,-2) {$v_3$};
    \draw[->] (1) to node[midway, above]  {$[1,1]$}  (2);
    \draw[->] (3) to node[midway, below]  {$[1,1]$} (4);
    \draw[->] (2) to [out=330,in=30,looseness=1] node[midway, right]  {$[1,2]$} (3);
    \draw[->] (4) to [out=150,in=210,looseness=1] node[midway, left]  {$[1,2]$} (1);
     \end{tikzpicture}

%% file: graphics/figure8.tex
    \begin{tikzpicture}[node distance={7mm}, inner sep=1pt, minimum size=4mm, main/.style = {draw, circle}, time/.style = {}] 
    
    \node[main] (4) at (2,2) {$2$};
    \node[main] (5) [below of =4] {$1$};
    \node[main] (6) [below of =5]  {$0$};
    \node[main] (7) at (2,-1) {$2$};
    \node[main] (8) [below of =7]  {$1$};
    \node[main] (9) [below of =8]  {$0$};
    \node[main] (10) at (5,2) {$2$};
    \node[main] (11) [below of =10]  {$1$};
    \node[main] (12) [below of =11]  {$0$};
    \node[main] (13) at (5,-1) {$2$};
    \node[main] (14) [below of =13]  {$1$};
    \node[main] (15) [below of =14]  {$0$};
    \draw[orange, ->] (6) --  (11);
    \draw[orange, ->] (5) to node[midway, above] {$\frac{1}{3}$}  (10);
    \draw[orange, ->] (4) --  (12);
    \draw[orange, ->] (13) --  (9);
    \draw[orange, ->] (14) --  (7);
    \draw[orange, ->] (15) --  (8);
    \draw[violet, ->] (10) to [out=330,in=30,looseness=1] (14);
    \draw[violet, ->] (10) to [out=330,in=30,looseness=1] (15);
    \draw[violet, ->] (11) to [out=330,in=30,looseness=1] (13);
    \draw[violet, ->] (11) to [out=330,in=30,looseness=1] (15);
    \draw[violet, ->] (12) to [out=330,in=30,looseness=1] (13);
    \draw[violet, ->] (12) to [out=330,in=30,looseness=1] (14);
    \draw[violet, ->] (9) to [out=150,in=210,looseness=1] node[midway, left] {$\frac{1}{6}$} (4);
    \draw[violet, ->] (9) to [out=150,in=210,looseness=1] (5);
    \draw[violet, ->] (8) to [out=150,in=210,looseness=1] (4);
    \draw[violet, ->] (8) to [out=150,in=210,looseness=1] (6);
    \draw[violet, ->] (7) to [out=150,in=210,looseness=1] (5);
    \draw[violet, ->] (7) to [out=150,in=210,looseness=1] (6);

    \end{tikzpicture}

%% file: graphics/figure9.tex
    \begin{tikzpicture}[node distance={7mm}, inner sep=1pt, minimum size=4mm, main/.style = {draw, circle}, time/.style = {}]
    
    \node[main] (4) at (9,2) {$2$};
    \node[main] (5) [below of =4] {$1$};
    \node[main] (6) [below of =5]  {$0$};
    \node[main] (7) at (9,-1) {$2$};
    \node[main] (8) [below of =7]  {$1$};
    \node[main] (9) [below of =8]  {$0$};
    \node[main] (10) at (12,2) {$2$};
    \node[main] (11) [below of =10]  {$1$};
    \node[main] (12) [below of =11]  {$0$};
    \node[main] (13) at (12,-1) {$2$};
    \node[main] (14) [below of =13]  {$1$};
    \node[main] (15) [below of =14]  {$0$};
    \draw[black!30!green, ->] (6) --  (11);
    \draw[red, ->] (5) --  (10);
    \draw[blue, ->] (4) --  (12);
    \draw[blue, ->] (13) --  (9);
    \draw[red, ->] (14) --  (7);
    \draw[black!30!green, ->] (15) --  (8);
    \draw[dotted, ->] (12) to [out=330,in=30,looseness=1] (14);
    \draw[red, ->] (10) to [out=330,in=30,looseness=1] node[pos=0.2, above right]  {$c_2$} (14);
    \draw[dotted, ->] (10) to [out=330,in=30,looseness=1] (15);
    \draw[dotted, ->] (11) to [out=330,in=30,looseness=1] (13);
    \draw[black!30!green, ->] (11) to [out=330,in=30,looseness=1] node[pos=0.8, below right]  {$c_3$} (15);
    \draw[dotted, ->] (8) to [out=150,in=210,looseness=1] (4);
    \draw[blue, ->] (12) to [out=330,in=30,looseness=1]node[midway, left]  {$c_1$} (13);
    \draw[blue, ->] (9) to [out=150,in=210,looseness=1] (4);
    \draw[dotted, ->] (9) to [out=150,in=210,looseness=1] (5);
    \draw[dotted, ->] (7) to [out=150,in=210,looseness=1] (6);
    \draw[black!30!green, ->] (8) to [out=150,in=210,looseness=1] (6);
    \draw[red, ->] (7) to [out=150,in=210,looseness=1] (5);
    
    \end{tikzpicture}

%% file: graphics/figure10.tex
    \begin{tikzpicture}[node distance={5mm}, inner sep=1pt, minimum size=2mm, main/.style = {draw, circle}, time/.style = {}]
    
    \node[main] (11d3p) at (0,5.5) {$3$};
    \node[main] (11d2p) [below of=11d3p]  {$2$};
    \node[main] (11d1p) [below of=11d2p] {$1$};
    \node[main] (11d0p) [below of=11d1p] {$0$};
    
    \node[main] (12a3p) at (2,5.5) {$3$};
    \node[main] (12a2p) [below of=12a3p] {$2$};
    \node[main] (12a1p) [below of=12a2p] {$1$};
    \node[main] (12a0p) [below of=12a1p] {$0$};
    
    \node[main] (12d3p) at (3,5.5) {$3$};
    \node[main] (12d2p) [below of=12d3p] {$2$};
    \node[main] (12d1p) [below of=12d2p] {$1$};
    \node[main] (12d0p) [below of=12d1p] {$0$};
    
    \node[main] (13d3p) at (5,5.5) {$3$};
    \node[main] (13d2p) [below of=13d3p] {$2$};
    \node[main] (13d1p) [below of=13d2p] {$1$};
    \node[main] (13d0p) [below of=13d1p] {$0$};
    
    \node[main] (13d3m) at (5,3) {$3$};
    \node[main] (13d2m) [below of=13d3m] {$2$};
    \node[main] (13d1m) [below of=13d2m] {$1$};
    \node[main] (13d0m) [below of=13d1m]  {$0$};
        
    \node[main] (12a3m) at (3,3) {$3$};
    \node[main] (12a2m) [below of=12a3m] {$2$};
    \node[main] (12a1m) [below of=12a2m] {$1$};
    \node[main] (12a0m) [below of=12a1m] {$0$};
    
    \node[main] (12d3m) at (2,3) {$3$};
    \node[main] (12d2m) [below of=12d3m] {$2$};
    \node[main] (12d1m) [below of=12d2m] {$1$};
    \node[main] (12d0m) [below of=12d1m] {$0$};
    
    \node[main] (11a3m) at (0,3) {$3$};
    \node[main] (11a2m) [below of=11a3m] {$2$};
    \node[main] (11a1m) [below of=11a2m] {$1$};
    \node[main] (11a0m) [below of=11a1m]{$0$};

    \node[main] (24d3p) at (0,0) {$3$};
    \node[main] (24d2p) [below of=24d3p] {$2$};
    \node[main] (24d1p) [below of=24d2p] {$1$};
    \node[main] (24d0p) [below of=24d1p] {$0$};
    
    \node[main] (22a3p) at (2,0) {$3$};
    \node[main] (22a2p) [below of=22a3p] {$2$};
    \node[main] (22a1p) [below of=22a2p] {$1$};
    \node[main] (22a0p) [below of=22a1p] {$0$};
    
    \node[main] (22d3p) at (3,0) {$3$};
    \node[main] (22d2p) [below of=22d3p] {$2$};
    \node[main] (22d1p) [below of=22d2p] {$1$};
    \node[main] (22d0p) [below of=22d1p] {$0$};

    \node[main] (25a3p) at (5,0) {$3$};
    \node[main] (25a2p) [below of=25a3p] {$2$};
    \node[main] (25a1p) [below of=25a2p] {$1$};
    \node[main] (25a0p) [below of=25a1p]{$0$};
    
    \node[main] (25d3m) at (5,-2.5) {$3$};
    \node[main] (25d2m) [below of=25d3m] {$2$};
    \node[main] (25d1m) [below of=25d2m] {$1$};
    \node[main] (25d0m) [below of=25d1m] {$0$};
    
    \node[main] (22a3m) at (3,-2.5) {$3$};
    \node[main] (22a2m) [below of=22a3m] {$2$};
    \node[main] (22a1m) [below of=22a2m] {$1$};
    \node[main] (22a0m) [below of=22a1m] {$0$};
    
    \node[main] (22d3m) at (2,-2.5) {$3$};
    \node[main] (22d2m) [below of=22d3m] {$2$};
    \node[main] (22d1m) [below of=22d2m] {$1$};
    \node[main] (22d0m) [below of=22d1m] {$0$};
    
    \node[main] (24a3m) at (0,-2.5) {$3$};
    \node[main] (24a2m) [below of=24a3m] {$2$};
    \node[main] (24a1m) [below of=24a2m] {$1$};
    \node[main] (24a0m) [below of=24a1m] {$0$};

    \draw[black!30!blue, ->] (11d3p) --  (12a0p);
    \draw[black!30!blue, ->] (12a0p) --  (12d0p);
    \draw[black!30!blue, ->] (12d0p) --  (13d1p);
    \draw[black!30!blue, ->] (13d1p) to [out=330, in=30, looseness=1] node[midway, right] {$\frac{1}{2}$} (13d3m);
    \draw[black!30!blue, ->] (13d3m) --  (12a0m);
    \draw[black!30!blue, ->] (12a0m) --  (12d0m);
    \draw[black!30!blue, ->] (12d0m) --  (11a1m);
    \draw[black!30!blue, ->] (11a1m) to [out=150, in=210, looseness=1] (11d3p);
    
    \draw[cyan, ->] (11d1p) --  (12a2p);
    \draw[cyan, ->] (12a2p) --  (12d2p);
    \draw[cyan, ->] (12d2p) --  (13d3p);
    \draw[cyan, ->] (13d3p) to [out=330, in=30, looseness=1] node[midway, right] {$\frac{1}{2}$} (13d1m);
    \draw[cyan, ->] (13d1m) --  (12a2m);
    \draw[cyan, ->] (12a2m) --  (12d2m);
    \draw[cyan, ->] (12d2m) --  (11a3m);
    \draw[cyan, ->] (11a3m) to [out=150, in=210, looseness=1] (11d1p);
    
    \draw[black!30!red, ->] (24d0p) --  (22a1p);
    \draw[black!30!red, ->] (22a1p) --  (22d1p);
    \draw[black!30!red, ->] (22d1p) --  (25a2p);
    \draw[black!30!red, ->] (25a2p) to [out=330, in=30, looseness=1]  (25d0m);
    \draw[black!30!red, ->] (25d0m) --  (22a1m);
    \draw[black!30!red, ->] (22a1m) --  (22d1m);
    \draw[black!30!red, ->] (22d1m) --  (24a2m);
    \draw[black!30!red, ->] (24a2m) to [out=150, in=210, looseness=1] node[midway, left] {$\frac{1}{2}$}(24d0p);
    
    \draw[orange, ->] (24d2p) --  (22a3p);
    \draw[orange, ->] (22a3p) --  (22d3p);
    \draw[orange, ->] (22d3p) --  (25a0p);
    \draw[orange, ->] (25a0p) to [out=330, in=30, looseness=1]  (25d2m);
    \draw[orange, ->] (25d2m) --  (22a3m);
    \draw[orange, ->] (22a3m) --  (22d3m);
    \draw[orange, ->] (22d3m) --  (24a0m);
    \draw[orange, ->] (24a0m) to [out=150, in=210, looseness=1] node[midway, left] {$\frac{1}{2}$} (24d2p);
    
    \draw[gray, ->] (22a3p) -- (12d0p);
    \draw[gray, ->] (22a1p) -- (12d2p);
    \draw[gray, ->] (12a2m) -- (22d3m);
    \draw[gray, ->] (12a0m) -- (22d1m);
    \draw[gray, ->] (22a1m) -- (12d2m);
    \draw[gray, ->] (22a3m) -- (12d0m);
    \draw[gray, ->] (12a2p) -- (22d3p);
    \draw[gray, ->] (12a0p) -- (22d1p);
    \draw[gray, ->] (12a2p) to [out=220, in=140, looseness=0.5] (22d3m);
    \draw[gray, ->] (12a0p) to [out=220, in=140, looseness=0.5] (22d1m);
    \draw[gray, ->] (22a3p) to [out=140, in=220, looseness=0.5] (12d0m);
    \draw[gray, ->] (22a1p) to [out=140, in=220, looseness=0.5] (12d2m);
    \draw[gray, ->] (22a3m) to [out=40, in=320, looseness=0.5] (12d0p);
    \draw[gray, ->] (22a1m) to [out=40, in=320, looseness=0.5] node[midway, right] {$\frac{1}{2}$} (12d2p);
    \draw[gray, ->] (12a0m) to [out=320, in=40, looseness=0.5] (22d1p);
    \draw[gray, ->] (12a2m) to [out=320, in=40, looseness=0.5] (22d3p);
    
    \end{tikzpicture}

%% file: graphics/figure12.tex
    \begin{tikzpicture}[node distance={12mm}, minimum size=4mm, main/.style = {draw, circle}, time/.style = {}]

    \node[main] (1) at (0,0) {};
    \node[main ](2) [right of=1] {};
    \node[main ](3) [right of=2] {};
    \node[main ](4) [right of=3] {};
    \node[main] (5) [below of=4] {};
    \node[main ](6) [left of=5] {};
    \node[main ](7) [left of=6] {};
    \node[main ](8) [left of=7] {};
    
    \node[main] (1b) [below of=8] {};
    \node[main ](2b) [right of=1b] {};
    \node[main ](3b) [right of=2b] {};
    \node[main ](4b) [right of=3b] {};
    \node[main] (5b) [below of=4b] {};
    \node[main ](6b) [left of=5b] {};
    \node[main ](7b) [left of=6b] {};
    \node[main ](8b) [left of=7b] {};

    \draw[blue, ->] (1) to node[midway, above]  {$1$}  (2);
    \draw[blue, ->] (2) to node[midway, above]  {$0$} (3);
    \draw[blue, ->] (3) to node[midway, above]  {$1$} (4);
    \draw[blue, ->] (4) to node[midway, right]  {$2$} (5); 
    \draw[blue, ->] (5) to node[midway, above]  {$1$} (6);
    \draw[blue, ->] (6) to node[midway, above]  {$0$} (7);
    \draw[blue, ->] (7) to node[midway, above]  {$1$} (8);
    \draw[blue, ->] (8) to node[midway, left]  {$2$} (1);
    
    \draw[red, ->] (1b) to node[midway, above]  {$1$}  (2b);
    \draw[red, ->] (2b) to node[midway, above]  {$0$} (3b);
    \draw[red, ->] (3b) to node[midway, above]  {$1$} (4b);
    \draw[red, ->] (4b) to node[midway, right]  {$2$} (5b); 
    \draw[red, ->] (5b) to node[midway, above]  {$1$} (6b);
    \draw[red, ->] (6b) to node[midway, above]  {$0$} (7b);
    \draw[red, ->] (7b) to node[midway, above]  {$1$} (8b);
    \draw[red, ->] (8b) to node[midway, left]  {$2$} (1b);
    
    \draw[gray, ->] (2b) to node[midway, left] {$1$} (7);
    \draw[gray, ->] (6) -- (3b);
    \draw[gray, ->] (6b) -- (7);
    \draw[gray, ->] (2) -- (3b);
    \draw[gray, ->] (2b) -- (3);
    \draw[gray, ->] (6) -- (7b);
    \draw[gray, ->] (2) to [out=220, in=140, looseness=0.5] (7b);
    \draw[gray, ->] (6b) to [out=40, in=320, looseness=0.5] (3);
    
    \end{tikzpicture}

%% file: graphics/figure13.tex
    \begin{tikzpicture}[node distance={12mm}, minimum size=4mm, main/.style = {draw, circle}, time/.style = {}]
    
    \node[time] (1c) at (8,0) {};
    \node[time ](2c) [right of=1c] {};
    \node[main ](3c) [right of=2c] {};
    \node[main ](4c) [right of=3c] {};
    \node[main] (5c) [below of=4c] {};
    \node[main ](6c) [left of=5c] {};
    \node[time ](7c) [left of=6c] {};
    \node[time ](8c) [left of=7c] {};
    
    \node[main] (1d) [below of=8c] {};
    \node[main ](2d) [right of=1d] {};
    \node[time ](3d) [right of=2d] {};
    \node[time ](4d) [right of=3d] {};
    \node[time] (5d) [below of=4d] {};
    \node[time ](6d) [left of=5d] {};
    \node[main ](7d) [left of=6d] {};
    \node[main ](8d) [left of=7d] {};
    
    \draw[blue, ->] (3c) to node[midway, above]  {$1$} (4c);
    \draw[blue, ->] (4c) to node[midway, right]  {$2$} (5c); 
    \draw[blue, ->] (5c) to node[midway, above]  {$1$} (6c);    
    \draw[red, ->] (1d) to node[midway, above]  {$1$}  (2d);
    \draw[red, ->] (7d) to node[midway, above]  {$1$} (8d);
    \draw[red, ->] (8d) to node[midway, left]  {$2$} (1d);
    \draw[gray, ->] (2d) to node[midway, left] {$1$} (3c);
    \draw[gray, ->] (6c) to node[midway, right] {$1$} (7d);

\end{tikzpicture}

%% file: graphics/figure11.tex
    \begin{tikzpicture}[node distance={5mm}, inner sep=1pt, minimum size=4mm, main/.style = {draw, circle}, time/.style = {}]

    \node[main] (12d3p) at (11,5.5) {$3$};
    \node[main] (12d2p) [below of=12d3p] {$2$};
    \node[main] (12d1p) [below of=12d2p] {$1$};
    \node[main] (12d0p) [below of=12d1p] {$0$};
    
    \node[main] (13d3p) at (13,5.5) {$3$};
    \node[main] (13d2p) [below of=13d3p] {$2$};
    \node[main] (13d1p) [below of=13d2p] {$1$};
    \node[main] (13d0p) [below of=13d1p] {$0$};
    
    \node[main] (13d3m) at (13,3) {$3$};
    \node[main] (13d2m) [below of=13d3m] {$2$};
    \node[main] (13d1m) [below of=13d2m] {$1$};
    \node[main] (13d0m) [below of=13d1m]  {$0$};
        
    \node[main] (12a3m) at (11,3) {$3$};
    \node[main] (12a2m) [below of=12a3m] {$2$};
    \node[main] (12a1m) [below of=12a2m] {$1$};
    \node[main] (12a0m) [below of=12a1m] {$0$};

   \node[main] (24d3p) at (8,0) {$3$};
    \node[main] (24d2p) [below of=24d3p] {$2$};
    \node[main] (24d1p) [below of=24d2p] {$1$};
    \node[main] (24d0p) [below of=24d1p] {$0$};
    
    \node[main] (22a3p) at (10,0) {$3$};
    \node[main] (22a2p) [below of=22a3p] {$2$};
    \node[main] (22a1p) [below of=22a2p] {$1$};
    \node[main] (22a0p) [below of=22a1p] {$0$};

    \node[main] (22d3m) at (10,-2.5) {$3$};
    \node[main] (22d2m) [below of=22d3m] {$2$};
    \node[main] (22d1m) [below of=22d2m] {$1$};
    \node[main] (22d0m) [below of=22d1m] {$0$};
    
    \node[main] (24a3m) at (8,-2.5) {$3$};
    \node[main] (24a2m) [below of=24a3m] {$2$};
    \node[main] (24a1m) [below of=24a2m] {$1$};
    \node[main] (24a0m) [below of=24a1m] {$0$};
    
    \draw[black!30!blue, ->] (12d0p) --  (13d1p);
    \draw[black!30!blue, ->] (13d1p) to [out=330, in=30, looseness=1] node[midway, right] {$\frac{1}{2}$} (13d3m);
    \draw[black!30!blue, ->] (13d3m) --  (12a0m);
    \draw[cyan, ->] (12d2p) --  (13d3p);
    \draw[cyan, ->] (13d3p) to [out=330, in=30, looseness=1] node[midway, right] {$\frac{1}{2}$} (13d1m);
    \draw[cyan, ->] (13d1m) --  (12a2m);        
    \draw[black!30!red, ->] (24d0p) --  (22a1p);
    \draw[black!30!red, ->] (22d1m) --  (24a2m);
    \draw[black!30!red, ->] (24a2m) to [out=150, in=210, looseness=1] node[midway, left] {$\frac{1}{2}$}(24d0p);   
    \draw[orange, ->] (24d2p) --  (22a3p);
    \draw[orange, ->] (22d3m) --  (24a0m);
    \draw[orange, ->] (24a0m) to [out=150, in=210, looseness=1] node[midway, left] {$\frac{1}{2}$} (24d2p);
    
    \draw[gray, ->] (22a3p) -- (12d0p);
    \draw[gray, ->] (22a1p) -- (12d2p);
    \draw[gray, ->] (12a2m) -- (22d3m);
    \draw[gray, ->] (12a0m) to node[midway, right] {$\frac{1}{2}$} (22d1m);

    \end{tikzpicture}

%% file: graphics/figure14.tex
    \begin{tikzpicture}[node distance={15mm}, minimum size=5mm, main/.style = {draw, circle}, time/.style = {}]
    \node[main] (1) at (0,0) {};
    \node[main ](2) at (3,0) {};
    \node[main ](3) at (3,-2) {};
    \node[main ](4) at (0,-2) {};
    \draw[->] (1) to node[midway, above]  {\textcolor{red}{$1$}}  (2);
    \draw[->] (3) to node[midway, below]  {\textcolor{red}{$1$}} (4);
    \draw[->] (2) to [out=330,in=30,looseness=1] node[midway, right]  {\textcolor{red}{$3$}} (3);
    \draw[->] (4) to [out=150,in=210,looseness=1] node[midway, left]  {\textcolor{red}{$3$}} (1);
     \end{tikzpicture}

%% file: graphics/figure15.tex
    \begin{tikzpicture}[node distance={15mm}, minimum size=5mm, main/.style = {draw, circle}, time/.style = {}]
    \node[main] (1) at (0,0) {};
    \node[main ](2) at (3,0) {};
    \node[main ](3) at (3,-2) {};
    \node[main ](4) at (0,-2) {};
    \draw[->] (1) to node[midway, above]  {\textcolor{red}{$1$}}  (2);
    \draw[->] (3) to node[midway, below]  {\textcolor{red}{$1$}} (4);
    \draw[->] (2) to [out=330,in=30,looseness=1] node[midway, right]  {\textcolor{red}{$2$}} (3);
    \draw[->] (4) to [out=150,in=210,looseness=1] node[midway, left]  {\textcolor{red}{$2$}} (1);
     \end{tikzpicture}

%% file: graphics/figure16.tex
    \begin{tikzpicture}[node distance={6mm}, inner sep=1pt, minimum size=4mm, main/.style = {draw, circle}, time/.style = {}]
    
    \node[main] (4) at (9,2) {$3$};
    \node[main] (45) [below of =4] {$2$};
    \node[main] (5) [below of =45] {$1$};
    \node[main] (6) [below of =5]  {$0$};
    \node[main] (7) at (9,-1) {$3$};
    \node[main] (78) [below of =7]  {$2$};
    \node[main] (8) [below of =78]  {$1$};
    \node[main] (9) [below of =8]  {$0$};
    \node[main] (10) at (12,2) {$3$};
    \node[main] (1011) [below of =10]  {$2$};
    \node[main] (11) [below of =1011]  {$1$};
    \node[main] (12) [below of =11]  {$0$};
    \node[main] (13) at (12,-1) {$3$};
    \node[main] (1314) [below of =13]  {$2$};
    \node[main] (14) [below of =1314]  {$1$};
    \node[main] (15) [below of =14]  {$0$};
    
    \draw[gray, ->] (45) --  (10);
    \draw[gray, ->] (5) --  (1011);
    \draw[gray, ->] (4) --  (12);
    \draw[gray, ->] (13) --  (9);
    \draw[gray, ->] (1314) --  (7);
    \draw[gray, ->] (14) --  (78);
    \draw[thick, red, ->] (15) --  (8);
    \draw[thick, red, ->] (6) --  (11);
    \draw[ gray, ->] (10) to [out=330,in=30,looseness=1] (1314);
    \draw[gray,  ->] (10) to [out=330,in=30,looseness=1] (14);
    \draw[gray, ->] (1011) to [out=330,in=30,looseness=1] (15);
    \draw[gray,  ->] (1011) to [out=330,in=30,looseness=1] (14);
    \draw[gray, ->] (11) to [out=330,in=30,looseness=1] (13);
    \draw[gray, ->] (12) to [out=330,in=30,looseness=1] (13);
    \draw[gray, ->] (12) to [out=330,in=30,looseness=1] (1314);
    \draw[gray, ->] (9) to [out=150,in=210,looseness=1] (45);
    \draw[gray, ->] (9) to [out=150,in=210,looseness=1] (4);
    \draw[gray, ->] (8) to [out=150,in=210,looseness=1] (4);
    \draw[gray, ->] (7) to [out=150,in=210,looseness=1] (45);
    \draw[gray, ->] (7) to [out=150,in=210,looseness=1] (5);
    \draw[gray, ->] (78) to [out=150,in=210,looseness=1] (5);
    \draw[gray, ->] (78) to [out=150,in=210,looseness=1] (6);
    \draw[thick, red, ->] (8) to [out=150,in=210,looseness=1] (6);
    \draw[thick, red, ->] (11) to [out=330,in=30,looseness=1] (15);
    
    \end{tikzpicture}

%% file: graphics/figure17.tex
    \begin{tikzpicture}[node distance={5mm}, inner sep=1pt, font = \scriptsize, minimum size=4mm, main/.style = {draw, circle}, time/.style = {}, invisi/.style = {}]
    
    \node[invisi] (A) at (4.1,6) {};
    \node[invisi] (B) at (4.1,3.8) {\ding{34}};
    \draw[gray, thick, dashed] (A) -- (B);
    
    \node[invisi, font=\small, red] (11d2p) at (4,6) {$\hat{v}$};
    \node[main, red] (11d2p) at (4,5.5) {$2$};
    \node[main, red] (11d1p) [below of=11d2p]  {$1$};
    \node[main, red] (11d0p) [below of=11d1p] {$0$};
    
    \node[main] (12a2p) at (7,5.5) {$2$};
    \node[main] (12a1p) [below of=12a2p] {$1$};
    \node[main] (12a0p) [below of=12a1p] {$0$};
    
    \node[main] (12d2p) at (9,5.5) {$2$};
    \node[main] (12d1p) [below of=12d2p] {$1$};
    \node[main] (12d0p) [below of=12d1p] {$0$};
    
    \node[main] (13d2p) at (12,5.5) {$2$};
    \node[main] (13d1p) [below of=13d2p] {$1$};
    \node[main] (13d0p) [below of=13d1p] {$0$};
    
    \node[main] (13d2m) at (12,3) {$2$};
    \node[main] (13d1m) [below of=13d2m] {$1$};
    \node[main] (13d0m) [below of=13d1m]  {$0$};
        
    \node[main] (12a2m) at (9,3) {$2$};
    \node[main] (12a1m) [below of=12a2m] {$1$};
    \node[main] (12a0m) [below of=12a1m] {$0$};
    
    \node[main] (12d2m) at (7,3) {$2$};
    \node[main] (12d1m) [below of=12d2m] {$1$};
    \node[main] (12d0m) [below of=12d1m] {$0$};
    
    \node[main] (11a2m) at (4,3) {$2$};
    \node[main] (11a1m) [below of=11a2m] {$1$};
    \node[main] (11a0m) [below of=11a1m]{$0$};

    \draw[red, ->] (11d2p) --  (12a0p);
    \draw[red, ->] (11d1p) --  (12a2p);
    \draw[red, ->] (11d0p) --  (12a1p);
    
    \draw[->] (12a0p) --  (12d1p);
    \draw[->] (12a0p) --  (12d2p);
    \draw[->] (12a1p) --  (12d2p);
    \draw[->] (12a1p) --  (12d0p);
    \draw[->] (12a2p) --  (12d0p);
    \draw[->] (12a2p) --  (12d1p);
    
    \draw[->] (12d2p) --  (13d0p);
    \draw[->] (12d1p) --  (13d2p);
    \draw[->] (12d0p) --  (13d1p);
    
    \draw[->] (13d0p) to [out=330, in=30, looseness=1]  (13d1m);
    \draw[->] (13d0p) to [out=330, in=30, looseness=1]  (13d2m);
    \draw[->] (13d1p) to [out=330, in=30, looseness=1]  (13d0m);
    \draw[->] (13d1p) to [out=330, in=30, looseness=1] (13d2m);
    \draw[->] (13d2p) to [out=330, in=30, looseness=1]  (13d0m);
    \draw[->] (13d2p) to [out=330, in=30, looseness=1]  (13d1m);
    
    \draw[->] (13d2m) --  (12a0m);
    \draw[->] (13d1m) --  (12a2m);
    \draw[->] (13d0m) --  (12a1m);
    
    \draw[->] (12a0m) --  (12d1m);
    \draw[->] (12a0m) --  (12d2m);
    \draw[->] (12a1m) --  (12d0m);
    \draw[->] (12a1m) --  (12d2m);
    \draw[->] (12a2m) --  (12d0m);
    \draw[->] (12a2m) --  (12d1m);
    
    \draw[->] (12d2m) --  (11a0m);
    \draw[->] (12d1m) --  (11a2m);
    \draw[->] (12d0m) --  (11a1m);
    
    \draw[->] (11a0m) to [out=150, in=210, looseness=1] (11d1p);
    \draw[->] (11a0m) to [out=150, in=210, looseness=1] (11d2p);
    \draw[->] (11a1m) to [out=150, in=210, looseness=1] (11d0p);
    \draw[->] (11a1m) to [out=150, in=210, looseness=1] (11d2p);
    \draw[->] (11a2m) to [out=150, in=210, looseness=1] (11d0p);
    \draw[->] (11a2m) to [out=150, in=210, looseness=1] (11d1p);

    \end{tikzpicture}

%% file: graphics/figure18.tex
    \begin{tikzpicture}[node distance={5mm}, inner sep = 1pt, font=\scriptsize, minimum size=4mm, main/.style = {draw, circle}, time/.style = {}, invisi/.style = {}]

    \node[invisi] (blub) at (5.5,6.5) {}; 
    
    \node[invisi, font=\small, red] (b) at (5.5,6) {$\hat{v}$};
    \node[main, red] (11d2p_art) at (5.5,5.5) {$2$};
    \node[main, red] (11d1p_art) [below of=11d2p_art]  {$1$};
    \node[main, red] (11d0p_art) [below of=11d1p_art] {$0$};
    
    \node[invisi, font=\small, red] (b) at (4,6) {$\hat{v}'$};
    \node[main, red] (11d2p) at (4,5.5) {$2$};
    \node[main, red] (11d1p) [below of=11d2p]  {$1$};
    \node[main, red] (11d0p) [below of=11d1p] {$0$};
    
    \node[main] (12a2p) at (7,5.5) {$2$};
    \node[main] (12a1p) [below of=12a2p] {$1$};
    \node[main] (12a0p) [below of=12a1p] {$0$};
    
    \node[main] (12d2p) at (9,5.5) {$2$};
    \node[main] (12d1p) [below of=12d2p] {$1$};
    \node[main] (12d0p) [below of=12d1p] {$0$};
    
    \node[main] (13d2p) at (12,5.5) {$2$};
    \node[main] (13d1p) [below of=13d2p] {$1$};
    \node[main] (13d0p) [below of=13d1p] {$0$};
    
    \node[main] (13d2m) at (12,3) {$2$};
    \node[main] (13d1m) [below of=13d2m] {$1$};
    \node[main] (13d0m) [below of=13d1m]  {$0$};
        
    \node[main] (12a2m) at (9,3) {$2$};
    \node[main] (12a1m) [below of=12a2m] {$1$};
    \node[main] (12a0m) [below of=12a1m] {$0$};
    
    \node[main] (12d2m) at (7,3) {$2$};
    \node[main] (12d1m) [below of=12d2m] {$1$};
    \node[main] (12d0m) [below of=12d1m] {$0$};
    
    \node[main] (11a2m) at (4,3) {$2$};
    \node[main] (11a1m) [below of=11a2m] {$1$};
    \node[main] (11a0m) [below of=11a1m]{$0$};

    \draw[red, ->] (11d2p_art) --  (12a0p);
    \draw[red, ->] (11d1p_art) --  (12a2p);
    \draw[red, ->] (11d0p_art) --  (12a1p);
    
    \draw[->] (12a0p) --  (12d1p);
    \draw[->] (12a0p) --  (12d2p);
    \draw[->] (12a1p) --  (12d2p);
    \draw[->] (12a1p) --  (12d0p);
    \draw[->] (12a2p) --  (12d0p);
    \draw[->] (12a2p) --  (12d1p);
    
    \draw[->] (12d2p) --  (13d0p);
    \draw[->] (12d1p) --  (13d2p);
    \draw[->] (12d0p) --  (13d1p);
    
    \draw[->] (13d0p) to [out=330, in=30, looseness=1]  (13d1m);
    \draw[->] (13d0p) to [out=330, in=30, looseness=1]  (13d2m);
    \draw[->] (13d1p) to [out=330, in=30, looseness=1]  (13d0m);
    \draw[->] (13d1p) to [out=330, in=30, looseness=1] (13d2m);
    \draw[->] (13d2p) to [out=330, in=30, looseness=1]  (13d0m);
    \draw[->] (13d2p) to [out=330, in=30, looseness=1]  (13d1m);
    
    \draw[->] (13d2m) --  (12a0m);
    \draw[->] (13d1m) --  (12a2m);
    \draw[->] (13d0m) --  (12a1m);
    
    \draw[->] (12a0m) --  (12d1m);
    \draw[->] (12a0m) --  (12d2m);
    \draw[->] (12a1m) --  (12d0m);
    \draw[->] (12a1m) --  (12d2m);
    \draw[->] (12a2m) --  (12d0m);
    \draw[->] (12a2m) --  (12d1m);
    
    \draw[->] (12d2m) --  (11a0m);
    \draw[->] (12d1m) --  (11a2m);
    \draw[->] (12d0m) --  (11a1m);
    
    \draw[->] (11a0m) to [out=150, in=210, looseness=1] (11d1p);
    \draw[->] (11a0m) to [out=150, in=210, looseness=1] (11d2p);
    \draw[->] (11a1m) to [out=150, in=210, looseness=1] (11d0p);
    \draw[->] (11a1m) to [out=150, in=210, looseness=1] (11d2p);
    \draw[->] (11a2m) to [out=150, in=210, looseness=1] (11d0p);
    \draw[->] (11a2m) to [out=150, in=210, looseness=1] (11d1p);

    \end{tikzpicture}

%% file: graphics/figure20.tex
  
    \begin{tikzpicture}[node distance={15mm},  inner sep=1pt, minimum size=5mm, main/.style = {draw, circle}, time/.style = {}, invisi/.style = {}]
    \node[main] (1) at (0,0) { $2$};
    \node[main] (2) [left of=1] { $1$};
    \node[main] (3) [right of=1] {$3$ };
    \node[main] (4) [above of=1] {$4$};
    \node[main] (5) [below of=1] {$5$ };
    
    \draw[blue] (2) --  (1) -- (3);
    \draw[red] (4) --  (1) -- (5);
    
    \end{tikzpicture}

%% file: graphics/figure21.tex
\begin{tikzpicture}[node distance={15mm}, minimum size=5mm, font=\footnotesize, main/.style = {draw, }, time/.style = {}, invisi/.style = {}]
    \node[main] (KL) at (-2,0) {Krumme Lanke};
    \node[main] (RS) at (3,0) {Rathaus Steglitz};
    \node[main] (AM) at (7,0) {Alt Mariendorf};
    \node[main] (S) at (3,3) {Spichernstr};
    \node[main] (HT) at (7,3) {Hallesches Tor};
    \node[main] (WS) at (12,3) {Warschauer Str};
    \node[main] (L) at (5,6) {Leopoldplatz};
    \node[main] (AT) at (2,9) {Alt Tegel};
    \node[main] (OS) at (7,8) {Osloer Str};
    
    \draw[orange] (RS) to node[midway, right] {$9$} (S) to node[midway, right] {$11$} (L) to node[midway, right] {$2$}  (OS);
    \draw[violet] (AT) to node[midway, right] {$13$}   (L) to node[midway, right] {$12$}  (HT) to node[midway, right] {$13$} (AM);
    \draw[black!30!green] (KL) to node[midway, above] {$19$}   (S) to node[midway, above] {$12$}  (HT) to node[midway, above] {$9$}  (WS);
    
    \end{tikzpicture}

%% file: graphics/figure22.tex
\begin{tikzpicture}[node distance={15mm}, minimum size=5mm, font=\small, inner sep =1pt, main/.style = {draw, }, time/.style = {}, invisi/.style = {}]
    \node[main] (K) at (-2,-1) {K};
    \node[main] (Rz) at (3,-1) {Rz};
    \node[main] (Mf) at (8,-1) {Mf};
    \node[main] (R) at (13,-1) {R};
    \node[main] (Sn) at (3,3) {Sn};
    \node[main] (U) at (2,4) {U};
    \node[main] (H) at (8,4) {H};
    \node[main] (Mi) at (8,5) {Mi};
    \node[main] (UD) at (8,6) {UD};
    \node[main] (A) at (10,6) {A};
    \node[main] (PA) at (11,8) {PA};
    \node[main] (Hö) at (13,6) {Hö};
    \node[main] (Me) at (8,3) {Me};
    \node[main] (WA) at (12,4) {WA};
    \node[main] (Lp) at (5,8) {Lp};
    \node[main] (Tg) at (2,10) {Tg};
    \node[main] (WI) at (3,11) {WI};
    \node[main] (Ol) at (7,9) {Ol};
    \node[main] (Kf) at (3,4) {Kf};
    \node[main] (Z) at (3,5) {Z};
    \node[main] (Bm) at (1,5) {Bm};
    \node[main] (Rl) at (-1,5) {Rl};
    \node[main] (RSp) at (-2.5,6) {RSp};
    \node[main] (Wt) at (4,4) {Wt};
     \node[main] (N) at (5,4) {N};
     \node[main] (Bp) at (5,2) {Bp};
     \node[main] (Ip) at (5,0) {Ip};
     \node[main] (G) at (6,4) {G};
     \node[main] (HBF) at (5,6) {HBF};
     \node[main] (M) at (7,4) {M};
     \node[main] (Kb) at (10,4) {Kb};
     \node[main] (Hp) at (10,2) {Hp};
     \node[main] (HMS) at (10,0) {HMS};
     \node[main] (Fp) at (1,1) {Fp};
      \node[main] (Be) at (3,1) {Be};
    
    \draw[thick, orange] (Rz) to node[midway, right] {$7$} (Be) to node[midway, right] {$2$} (Sn) to node[midway, left] {$1$} (Kf) to node[midway, left] {$2$} (Z) to node[midway, right] {$9$} (Lp) to node[midway, below] {$2$}  (Ol);
    \draw[thick, violet] (Tg) to node[midway, right] {$13$}   (Lp) to node[midway, right] {$8$}  (UD) to node[midway, right] {$2$}  (Mi) to node[midway, right] {$2$}  (H) to node[midway, right] {$2$}  (Me) to node[midway, right] {$11$} (Mf);
    \draw[thick, black!30!green] (K) to node[midway, above] {$16$}  (Fp) to node[midway, above] {$3$}  (Sn) to node[midway, below] {$3$}  (Wt) to [out=10, in=170, looseness=1 ] node[near end, above] {$2$}  (N) to [out=10, in=170, looseness=1 ] node[near end, above] {$4$}   (G)  to [out=10, in=170, looseness=1 ]  node[midway, above] {$1$} (M) to  [out=10, in=170, looseness=1 ] node[midway, above] {$2$} (H) to  [out=10, in=170, looseness=1 ] node[midway, above] {$4$}  (Kb) to  [out=10, in=170, looseness=1 ] node[midway, above] {$5$}  (WA);
     \draw[thick, green] (U) to node[midway, above] {$1$}  (Kf) to node[midway, above] {$1$}  (Wt) to node[near start, above] {$2$}  (N) to node[near start, above] {$4$}   (G)  to node[midway, below] {$1$} (M) to node[midway, below] {$2$} (H) to node[midway, below] {$4$}  (Kb) to node[midway, below] {$5$}  (WA);
      \draw[thick, yellow] (N) to node[midway, right] {$3$}   (Bp)  to node[midway, right] {$3$} (Ip);
       \draw[thick, blue] (WI) to node[midway, above] {$11$}   (Ol)  to node[midway, right] {$11$}   (A) to node[midway, right] {$6$}   (Kb) to node[midway, right] {$4$} (Hp) to node[midway, left] {$4$} (HMS) ;
       \draw[thick, brown] (Hö) to node[midway, above] {$33$}   (A) to node[midway, above] {$4$}   (UD) to node[midway, above] {$4$}   (HBF); 
       \draw[thick, red] (PA) to node[midway, left] {$11$}   (A) to node[midway, above] {$9$}   (Mi) to node[midway, above] {$5$}   (G) to [out=190, in=350, looseness=1 ] node[midway, below] {$4$}   (N) to [out=190, in=350, looseness=1 ] node[midway, below] {$1$}   (Wt) to node[midway, above] {$3$}   (Z) to node[midway, above] {$4$}   (Bm) to node[midway, below] {$10$}   (Rl);
       \draw[thick, cyan] (RSp) to node[midway, above] {$17$}   (Bm) to node[midway, left] {$5$}   (Fp) to node[midway, below] {$3$}   (Be) to node[midway, below] {$2$}   (Bp) to node[midway, below] {$6$}   (M) to node[midway, below] {$2$}   (Me) to node[midway, above] {$5$}   (Hp) to node[midway, right] {$17$}   (R);
    
    \end{tikzpicture}